\documentclass[12pt,smallextended]{svjour3}
\smartqed
\usepackage{lineno}
\modulolinenumbers[5]

\journalname{Journal of XXXX}
\usepackage{amssymb,amsmath,amsfonts}
\usepackage{epsfig}
\usepackage[bookmarksnumbered, colorlinks, urlcolor=blue, linkcolor=blue, citecolor=red, plainpages=false, pdfstartview=FitW]{hyperref}
\usepackage{mathptmx}
\usepackage{diagbox}
\usepackage{graphicx}
\usepackage{color}
\usepackage{caption}
\usepackage{wrapfig}
\usepackage{graphics}
\usepackage{enumerate}
\usepackage{amsxtra}
\usepackage{latexsym}
\usepackage{multirow}
\usepackage{slashbox}
\usepackage{subfigure}
\usepackage{epstopdf}
\usepackage{pdfpages}
\usepackage{algorithm}
\usepackage{algorithmic}

\usepackage{url}
\newcommand{\thmlist}{
\begin{list}{Step 1}
{\setlength{\leftmargin}{0.6 in}\setlength{\labelwidth} {0.5 in}}}

\newcommand{\alglist}{
\begin{list}{Step 1}
{\setlength{\leftmargin}{1.1 in} \setlength{\labelwidth}{1.0 in}}}

 \renewcommand{\proof} {\noindent {\bf Proof.} \quad}
 \newcommand{\eproof} {$\quad \square$}
 
 \newtheorem{assumption}{Assumption}

\renewcommand{\subtitle}[1]{\color{blue}}

\begin{document}

%\begin{frontmatter}

\title{The regularization continuation method for optimization problems
with nonlinear equality constraints}
\titlerunning{The regularization continuation method for nonlinear optimization}
\author{Xin-long Luo\textsuperscript{$\ast$}  \and Hang Xiao \and Sen Zhang}
\authorrunning{Luo \and Xiao \and Zhang}

\institute{
     Xin-long Luo
     \at
     Corresponding author. School of Artificial Intelligence, \\
     Beijing University of Posts and Telecommunications, P. O. Box 101, \\
     Xitucheng Road  No. 10, Haidian District, 100876, Beijing China\\
     \email{luoxinlong@bupt.edu.cn}            %  \\
     \and
     Hang Xiao
     \at
     School of Artificial Intelligence, \\
     Beijing University of Posts and Telecommunications, P. O. Box 101, \\
     Xitucheng Road  No. 10, Haidian District, 100876, Beijing China \\
     \email{xiaohang0210@bupt.edu.cn}
     \and
     Sen Zhang
     \at
     School of Artificial Intelligence, \\
     Beijing University of Posts and Telecommunications, P. O. Box 101, \\
     Xitucheng Road  No. 10, Haidian District, 100876, Beijing China \\
     \email{senzhang@bupt.edu.cn}
}

\date{Received: date / Accepted: date}
\maketitle

\begin{abstract}
This paper considers the regularization continuation method and the
trust-region updating strategy for the nonlinearly equality-constrained optimization
problem. Namely, it uses the inverse of the regularization quasi-Newton matrix
as the pre-conditioner to improve its computational efficiency in the well-posed
phase, and it adopts the inverse of the regularization two-sided projection
of the Hessian as the pre-conditioner to improve its robustness in the
ill-conditioned phase. Since it only solves a linear system of equations at
every iteration and the sequential quadratic programming (SQP) needs to solve a
quadratic programming subproblem at every iteration, it is faster than SQP.
Numerical results also show that it is more robust and faster than
SQP (the built-in subroutine fmincon.m of the MATLAB2020a environment
and the subroutine SNOPT executed in GAMS v28.2 (2019) environment).
The computational time of the new method is about one third of that of fmincon.m
for the large-scale problem. Finally, the global convergence analysis of the new
method is also given.
\end{abstract}



\vskip 2mm

%\textbf{AMS subject classifications.} 65K05 \and 65L05 \and 65L20

%\end{frontmatter}

% \linenumbers
% main text

\section{Introduction} \label{SUBINT}

% \vskip 2mm

In this article, we consider the optimization problem with nonlinear equality
constraints  as follows:
\begin{align}
  &\min_{x \in \Re^n} \; f(x)  \nonumber \\
  &\text{subject to} \; \; c(x) = 0,   \label{NLEQOPT}
\end{align}
where $f: \Re^{n} \to \Re$ and $c: \Re^{n} \to \Re^{m} \; (m \le n)$.
This problem has many applications in engineering
fields such as the visual-inertial navigation of an unmanned aerial
vehicle maintaining the horizontal flight \cite{CMFO2009,LLS2022},
constrained sparse regression \cite{BF2010}, sparse signal recovery
\cite{FS2009,VLLW2006}, image restoration and de-noising
\cite{FB2010,NWY2010,ST2010}, the Dantzig selector \cite{LPZ2012},
support vector machines \cite{FCG2010}, and sparse principal component analysis
(PCA) methods \cite{DGJL2007,WTH2009,WY2013}. And there are many
practical methods \cite{LJYY2019,LY2011,NW1999} and many efficient solvers
(the built-in subroutine fmincon \cite{Schittkowski1986} of MATLAB 2020a
\cite{MATLAB} and the subroutine SNOPT \cite{GMS2005,GMS2006} in the
GAMS environment \cite{GAMS}) to solve it based on the sequential quadratic
programming (SQP) method \cite{Han1977,Powell1978a,Wilson1963}.

\vskip 2mm

For the constrained optimization problem \eqref{NLEQOPT}, the continuation
method \cite{AG2003,CKK2003,Goh2011,KLQCRW2008,Pan1992,Tanabe1980} is another
method other than the traditional optimization method such as SQP
\cite{Han1977,Powell1978a,Powell1978b,Wilson1963} or the penalty function
method \cite{FM1990}. The advantage of the continuation method over the SQP method
is that the continuation method is capable of finding many local optimal points
of the non-convex optimization problem by tracking its trajectory, and it is
even possible to find the global optimal solution
\cite{BB1989,LXZ2022,Schropp2000,Yamashita1980}. However, the computational efficiency
of the classical continuation methods is lower than that of the traditional
optimization method such as SQP. Recently, Luo et al intensively investigate
the continuation methods for nonlinear equations \cite{LXL2022,LX2021},
linear programming problems \cite{LY2022}, linear complementarity problems \cite{LZX2022},
unconstrained optimization problems \cite{LXLZ2021,LXZ2022}, and
linearly constrained optimization problems \cite{LLS2022,LX2022}.
And they obtain the more robust and faster continuation methods than the
traditional optimization methods.

\vskip 2mm

Here, we extend their ideas to the optimization problem with nonlinear equality
constraints. In order to improve the computational efficiency of the continuation
method for the large-scale optimization problem further,
we consider the regularization BFGS (a special quasi-Newton method)
preconditioned technique for the regularization continuation
method and use an adaptive time step control based on the trust-region
updating strategy in this article. Moreover, in order to improve its robustness,
we replace the inverse of the regularization BFGS matrix with
the inverse of the regularization two-sided projection
of the Hessian matrix in the ill-posed phase.

\vskip 2mm

The rest of the paper is organized as follows. In Section 2, we give the
regularization continuation method with the switching preconditioned technique
and the trust-region updating strategy for the problem \eqref{NLEQOPT}. In
Section 3, we analyze the global convergence of this new method. In Section 4,
we report some promising numerical results of the new method, in comparison to
the state-of-art optimization methods (SQP, the built-in subroutine fmincon
\cite{Schittkowski1986} of MATLAB 2020a \cite{MATLAB} and the subroutine
SNOPT \cite{GMS2005,GMS2006} in the GAMS environment \cite{GAMS}) for some
large-scale problems. Finally, we give some discussions and conclusions in Section 5.

%\vskip 2mm

\section{The regularization continuation method}

%\vskip 2mm

In this section, we give the regularization continuation method with
the switching preconditioned technique and the adaptive time-stepping scheme
based on the trust-region updating strategy \cite{CGT2000} for the optimization
problem with nonlinear equality constraints. Firstly, we consider the
regularized projection gradient flow based on the KKT conditions of
equality-constrained optimization problems. Then, we construct the regularization
continuation method with the trust-region updating strategy to follow this special
ordinary differential equations (ODEs). The new method uses the
Broyden-Fletcher-Goldfarb-Shanno (BFGS) updating formula
\cite{Broyden1970,Fletcher1970,Goldfarb1970,Shanno1970} as the preconditioned
technique to improve its computational efficiency in the well-posed phase, and
it adopts the inverse of the regularized two-sided projection Hessian
as the pre-conditioner to improve its robustness in the ill-posed phase.
Furthermore, we use the generalized continuation Newton method \cite{LX2021} to
find an initial feasible point.

\vskip 2mm

\subsection{The regularized projection Newton flow}

\vskip 2mm

For the nonlinearly constrained optimization problem \eqref{NLEQOPT}, it is well known
that its optimal solution $x^{\ast}$ needs to satisfy the Karush-Kuhn-Tucker
conditions (p. 328, \cite{NW1999}) as follows:
\begin{align}
  \nabla_{x} L(x, \, \lambda) &= \nabla f(x) + A(x)^{T} \lambda = 0,
    \label{FOKKTG} \\
   c(x) & = 0,             \label{FOKKTC}
\end{align}
where $A(x)^{T} = [\nabla c_{1}(x), \, \ldots, \, \nabla c_{m}(x)]$ and the
Lagrangian function $L(x, \, \lambda)$ is defined by
\begin{align}
  L(x, \, \lambda) = f(x) + \lambda^{T}c(x).
      \label{LAGFUN}
\end{align}
Similarly to the method of the negative gradient flow for the unconstrained
optimization problem \cite{HM1996}, from the first-order necessary conditions
\eqref{FOKKTG}-\eqref{FOKKTC}, we can construct a dynamical system of
differential-algebraic equations for problem \eqref{NLEQOPT}
\cite{CL2011,LL2010,Luo2012,LLW2013,LLS2022,LX2022,Schropp2003} as follows:
\begin{align}
    & \frac{dx}{dt} = - \nabla L_{x}(x, \, \lambda)
      = -\left(\nabla f(x) + A(x)^{T} \lambda \right),  \label{DAGF} \\
    & c(x) = 0.                      \label{LACON}
\end{align}

\vskip 2mm

By differentiating the algebraic constraint \eqref{LACON} with respect to $t$
and substituting it into the differential equation \eqref{DAGF}, we obtain
\begin{align}
  & \frac{dc(x)}{dt} = A(x)\frac{dx}{dt}
  = - A(x) \left(\nabla f(x) + A(x)^{T}\lambda \right) \nonumber \\
  & \hskip 2mm = - A(x) \nabla f(x) - A(x)A(x)^{T} \lambda = 0.    \label{DIFALGC}
\end{align}
If we assume that matrix $A(x)$ has full row rank further, from equation
\eqref{DIFALGC}, we obtain
\begin{align}
   \lambda = - \left(A(x)A(x)^{T} \right)^{-1} A(x) \nabla f(x). \label{LAMBDA}
\end{align}
By substituting $\lambda$ of equation \eqref{LAMBDA} into equation \eqref{DAGF},
we obtain the projected gradient flow \cite{Tanabe1980} for the constrained
optimization problem \eqref{NLEQOPT} as follows:
\begin{align}
  \frac{dx}{dt} = - \left( I - A(x)^{T} \left(A(x)A(x)^{T}\right)^{-1}A(x)\right)
  \nabla f(x) = - P(x)g(x), \label{ODGF}
\end{align}
where $g(x) = \nabla f(x)$ and the projection matrix $P(x)$ is defined by
\begin{align}
  P(x)  = I - A(x)^{T} \left(A(x)A(x)^{T}\right)^{-1}A(x).  \label{PROMAT}
\end{align}

\vskip 2mm

It is not difficult to verify $(P(x))^{2} = P(x)$. That is to say, the projection matrix
$P(x)$ is symmetric and its eigenvalues are 0 or 1. From Theorem 2.3.1 in p. 73
of \cite{GV2013}, we know that its matrix 2-norm is
\begin{align}
     \|P(x)\| = 1. \label{MATNP}
\end{align}
We denote $P(x)^{\dagger}$ as the Moore-Penrose generalized inverse of the projection
matrix $P(x)$ (see p. 11, \cite{SY2006} or \cite{Moore1920,Penrose1955}). Since the
projection matrix $P(x)$ is symmetric and $(P(x))^{2} = P(x)$, it is not difficult to
verify
\begin{align}
      P(x)^{\dagger} = P(x).     \label{GINVP}
\end{align}
Actually, from equation \eqref{GINVP}, we have $P(x)(P(x))^{\dagger}P(x) = P(x)$,
$P(x)^{\dagger}P(x)P(x)^{\dagger} = P(x)^{\dagger}$,
$\left(P(x)^{\dagger}P(x)\right)^{T} = P(x)^{\dagger}P(x)$
and $\left(P(x)P(x)^{\dagger}\right)^{T} = P(x)P(x)^{\dagger}$.

% \vskip 2mm

\begin{remark}
If $x(t)$ is the solution of the ODE \eqref{ODGF}, it is not difficult to verify
that $x(t)$ satisfies $A(x) (dx/dt) = 0$. That is to say, if the initial point
$x_{0}$ satisfies $c(x_{0}) = 0$, the solution $x(t)$ of the projected
gradient flow \eqref{ODGF} also satisfies $c(x) = 0, \; \forall t \ge 0$. This
property is very useful such that we can construct an ODE method to follow the
trajectory of the ODE \eqref{ODGF} and obtain its steady-state solution $x^{\ast}$. \qed
\end{remark}

% \vskip 2mm

If we assume that $x(t)$ is the solution of the ODE \eqref{ODGF}, by using the
property $(P(x))^{2} = P(x)$, we obtain
\begin{align}
   & \frac{df(x)}{dt} = \left(\nabla f(x)\right)^{T} \frac{dx}{dt}
    = - (\nabla f(x))^{T} P(x) \nabla f(x) \nonumber \\
   & \hskip 2mm =  - g(x)^{T} (P(x))^{2} g(x)
    = - \|P(x)g(x)\|^{2} \le 0.  \nonumber
\end{align}
That is to say, $f(x)$ is monotonically decreasing along the solution curve $x(t)$
of the dynamical system \eqref{ODGF}. Furthermore, the solution $x(t)$ converges
to $x^{\ast}$ when $f(x)$ is lower bounded and $t$ tends to infinity
\cite{HM1996,Schropp2000,Tanabe1980}, where $x^{\ast}$ satisfies the first-order
Karush-Kuhn-Tucker conditions \eqref{FOKKTG}-\eqref{FOKKTC}. Thus, we can follow
the trajectory $x(t)$ of the ODE \eqref{ODGF} to obtain its steady-state solution
$x^{\ast}$, which is also one stationary point of the original optimization problem
\eqref{NLEQOPT}.

\vskip 2mm

However, since the right-hand-side function $P(x)g(x)$ of the ODE \eqref{ODGF}
is rank-deficient, we will confront the numerical difficulties when we use
the explicit ODE method to follow the projected gradient flow \eqref{ODGF}
\cite{AP1998,BCP1996,BJ1998}.  In order to mitigate the stiffness of the ODE
\eqref{ODGF}, we use the generalized inverse $(P(x)\nabla^{2} f(x)P(x))^{\dagger}$
of the two-sided projection matrix $P(x) \nabla^{2} f(x) P(x)$ as the pre-conditioner
for the ODE \eqref{ODGF}, which is used similarly to the system of nonlinear
equations \cite{LXL2022}, the unconstrained optimization problem
\cite{HM1996,LXLZ2021,LXZ2022}, the linear programming problem \cite{LY2022},
the linear complementarity problem \cite{LZX2022}, the underdetermined system
of nonlinear equations \cite{LLS2022,LX2022} and the linearly constrained optimization
problem \cite{LX2022}.

\vskip 2mm

By using the generalized inverse $\left(P(x)\nabla^{2}f(x)P(x)\right)^{\dagger}$ of
$P(x)\nabla^{2}f(x)P(x)$ as the pre-conditioner for the ODE \eqref{ODGF}, we
obtain the projected Newton flow
\begin{align}
    \frac{dx(t)}{dt} = - \left(P(x)\nabla^{2}f(x)P(x)\right)^{\dagger}P(x)g(x).
    \label{PCODGF}
\end{align}
Since $P(x)\nabla^{2} f(x)P(x)$ is singular, we reformulate equation \eqref{PCODGF} as
\begin{align}
    \left(P(x) \nabla^{2}f(x)P(x)\right) \frac{dx(t)}{dt} = - P(x)g(x). \label{DODGF}
\end{align}

\vskip 2mm

Although the projected Newton flow \eqref{DODGF} mitigates the
stiffness of the ODE such that we can adopt the explicit ODE method
to integrate it on the infinite interval, there are two disadvantages. One
is that the two-side projection matrix $P(x)\nabla^{2}f(x)P(x)$ may be not
positive semi-definite. Consequently, it can not ensure the objective function
$f(x)$ is monotonically decreasing along the solution $x(t)$ of the ODE \eqref{DODGF}.
The other is that the solution $x(t)$ of the ODE \eqref{DODGF} can not ensure
to satisfy the feasibility, i.e., $dc(x)/dt = A(x)dx(t)/dt = 0$.  In order to
overcome these two shortcomings, we use the  regularization technique
\cite{Hansen1994,Tikhonov1943,TA1977} for the projected Newton flow
\eqref{DODGF} as follows:
\begin{align}
    \left(\sigma(x) I + P(x) \nabla^{2}f(x)P(x)\right) \frac{dx(t)}{dt}
    = - P(x)g(x),     \label{TRODGF}
\end{align}
where the regularization parameter $\sigma(x)$ satisfies
$\sigma(x) + \mu_{min}\left(P(x) \nabla^{2}f(x)P(x)\right)$ $ \ge \sigma_{min} > 0$.
Here, $\mu_{min}(B)$ represents the smallest eigenvalue of matrix $B$.

% \vskip 2mm

\begin{remark} \label{RMLCC}
If we assume that $x(t)$ is the solution of the ODE \eqref{TRODGF}, from the
property $A(x)P(x) = 0$, we have
\begin{align}
    A(x)\left(\sigma(x) I + P(x) \nabla^{2}f(x) P(x)\right) \frac{dx(t)}{dt}
    = - A(x)P(x)g(x) = 0. \nonumber
\end{align}
Consequently, we obtain $\sigma(x) A(x)dx(t)/dt = 0$, i.e. $A(x)dx(t)/dt = 0$.
By integrating it, we obtain $c(x) = c(x_{0}) = 0$. That is to say, the solution
$x(t)$ of the ODE \eqref{TRODGF} preserves the feasibility $c(x) = 0$ when
$c(x_{0}) = 0$. \qed
\end{remark}

% \vskip 2mm

\begin{remark}
From $(P(x))^{2} = P(x)$, we know that the solution $x(t)$ of the ODE
\eqref{TRODGF} satisfies $P(x)dx(t)/dt = dx(t)/dt$. Consequently, from equation
\eqref{TRODGF} and the assumption $\sigma(x) +
\lambda_{min}\left(P(x) \nabla^{2}f(x)P(x)\right) \ge \sigma_{min} > 0$, we obtain
\begin{align}
     & \frac{df(x(t))}{dt} = (\nabla f(x))^{T}\frac{dx(t)}{dt}
     = (\nabla f(x))^{T}P(x) \frac{dx(t)}{dt} = (P(x)g(x))^{T}\frac{dx(t)}{dt}
     \nonumber \\
     & \hskip 2mm  =  - (P(x)g(x))^{T}
     \left(\sigma(x)I + P(x)\nabla^{2}f(x)P(x)\right)^{-1}(P(x)g(x)) \le 0. \nonumber
\end{align}
That is to say, $f(x)$ is monotonically decreasing
along the solution $x(t)$ of the ODE \eqref{TRODGF}. Furthermore, the solution $x(t)$
converges to $x^{\ast}$ when $f(x)$ is lower bounded and $\|P(x)\nabla^{2}f(x)P(x)\| \le M$
\cite{HM1996,LQQ2004,Schropp2000,Tanabe1980}, where $M$ is a positive constant
and $x^{\ast}$ is the stationary point of the  regularization projected
Newton flow \eqref{TRODGF}.  Thus, we can follow the trajectory $x(t)$ of the ODE
\eqref{TRODGF} to obtain its stationary point $x^{\ast}$, which is also one stationary
point of the original optimization problem \eqref{NLEQOPT}. \qed
\end{remark}

% \vskip 2mm

\subsection{The regularization continuation method} \label{SUBTRCM}

% \vskip 2mm

The solution curve $x(t)$ of the ODE \eqref{TRODGF} can  not be efficiently
followed by the general ODE method \cite{AP1998,BCP1996,BJ1998,JT1995}
such as backward differentiation formulas (BDFs, the subroutine ode15s.m of
the MATLAB R2020a environment \cite{MATLAB}). Thus, we need to construct
the particular method for this problem. We apply the first-order explicit Euler
method \cite{SGT2003} to the ODE \eqref{TRODGF}, then we obtain the regularized
projection Newton method as follows:
\begin{align}
    & \left(\sigma_{k} I + P(x_{k}) \nabla^{2}f(x_{k})P(x_{k})\right) d_{k}
    = - P(x_{k})g(x_{k}),     \label{TRNEWTON} \\
    & \hskip 2mm x_{k+1} = x_{k} + \alpha_{k} d_{k},
    \label{XK1LST}
\end{align}
where $\alpha_k$ is the time step. When $\alpha_{k} = 1$, the regularized
projection Newton method \eqref{TRNEWTON}-\eqref{XK1LST} is similar to the
Levenberg-Marquardt method \cite{Levenberg1944,LLT2007,Marquardt1963}.

\vskip 2mm

Since the time step $\alpha_{k}$ of the regularized projection Newton method
\eqref{TRNEWTON}-\eqref{XK1LST} is restricted by the numerical stability
\cite{SGT2003}. That is to say, for the linear test equation
$dx/dt = - \lambda x$, its time step $\alpha_{k}$ is restricted by the
stable region $|1-\alpha_{k}{\lambda}/(\sigma_{k} + \lambda)| \le 1$.
In order to avoid this disadvantage, similarly to the processing technique
of the nonlinear equations \cite{LXL2022,LY2022,LX2021}, the unconstrained
optimization problem \cite{LXLZ2021,LX2022} and the linearly constrained optimization
problem \cite{LX2022}, we replace $\alpha_{k}$ with $\Delta t_{k}/(1+\Delta t_{k})$
in equation \eqref{XK1LST} and let $\sigma_{k} = \sigma_{0}/{\Delta t_{k}}$ in
equation \eqref{TRNEWTON}. Then, we obtain the  regularization continuation
method:
\begin{align}
    & \left(\frac{\sigma_{0}}{\Delta t_{k}} I + B_{k}\right) d_{k}
    = - P_{k}g_{k}, \; s_{k}^{p} = \frac{\Delta t_{k}}{1 + \Delta t_{k}}(P_{k}d_{k}),
    \label{TRCM} \\
    &  \hskip 2mm x_{k+1}^{p} = x_{k} + s_{k}^{p},     \label{XK1P}
\end{align}
where $\Delta t_{k}$ is the time step, $P_{k} = P(x_k), \; g_{k}
= \nabla f(x_k)$ and $B_{k} = P(x_{k}) \nabla^{2}f(x_{k})P(x_{k})$ or its
quasi-Newton approximation.

\vskip 2mm

We denote $A_{k} = A(x_{k})$. The matrix $A_{k}A_{k}^{T}$ may be ill-conditioned.
Thus, the Cholesky factorization method may fail to solve the projection matrix
defined by equation \eqref{PROMAT} for the large-scale problem. Therefore, we use
the QR decomposition (pp. 247-248, \cite{GV2013}) to solve it as follows:
\begin{align}
   A_{k}^{T} = Q_{k}R_{k}, \label{QRAKT}
\end{align}
where $Q_{k} \in \Re^{n \times m}$ satisfies $Q_{k}^{T}Q_{k} = I$ and
$R_{k} \in \Re^{m \times m}$ is an upper triangle matrix. Consequently, the
projection matrix $P_{k}$ defined by equation \eqref{PROMAT} can be simplified
as follows:
\begin{align}
     P_{k} = I - Q_{k}Q_{k}^{T}. \label{SPROMAT}
\end{align}

\vskip 2mm

From equations \eqref{TRCM}-\eqref{XK1P} and the property $A(x_{k})P(x_{k}) = 0$
of the projection matrix $P(x_{k})$, it is not difficult to verify
$A(x_{k})s_{k}^{p} = 0$. Thus, when the constraint is linear, i.e.
$c(x) = Ax - b = 0$, $x_{k+1}^{p}$ also satisfies the linear constraint
$Ax_{k+1}^{p} - b = 0$ if $Ax_{k} - b = 0$ \cite{LX2022}. However, for
the nonlinear constraint $c(x) = 0$, this property is not true. In other words,
we may not have $c\left(x_{k+1}^{p}\right) = 0$ when $c(x_{k}) = 0$, where
$x_{k+1}^{p}$ is the solution of equations \eqref{TRCM}-\eqref{XK1P}. Therefore, we
need to compute the correction step $s_{k}^{c}$ such that $x_{k+1} = x_{k+1}^{p} + s_{k}^{c}$
is pulled back to the constraint $c(x) = 0$.
In order to save the computational time, we use the following generalized Newton iteration with the shortest
increment \cite{LX2021} to achieve this aim:
\begin{align}
    & s_{k}^{c} = -(A_{k+1}^{p})^{\dagger} c(x_{k+1}^{p}), \;
    (A_{k+1}^{p})^{\dagger} = (A_{k+1}^{p})^{T}\left(A_{k+1}^{p}(A_{k+1}^{p})^{T}\right)^{-1},
    \nonumber \\
    & s_{k} \triangleq s_{k}^{p} + s_{k}^{c}, \; x_{k+1} = x_{k+1}^{p} + s_{k}^{c} = x_{k} + s_{k},
      \label{GNEWTONM}
\end{align}
where $A_{k+1}^{p}$ equals $A(x_{k})$ or $A(x_{k+1}^{p})$. We first try $A_{k+1}^{p} = A(x_{k})$
and solve equation \eqref{GNEWTONM}. If $\|c(x_{k+1})\|_{\infty} > \epsilon_{0}$, we try
$A_{k+1}^{p} = A(x_{k+1}^{p})$ and solve equation \eqref{GNEWTONM} again, where $\epsilon_{0}$
is a given small constant that is less than the termination tolerance $\epsilon$
such as $\epsilon_{0} = \frac{1}{10} \epsilon$ and $\epsilon = 10^{-6}$.

\vskip 2mm

When $A_{k+1}^{p} = A(x_{k}) = A_k$, we have obtained its QR decomposition \eqref{QRAKT} and
let $Q_{k}^{p} = Q_{k}, \; R_{k}^{p} = R_{k}$. If $A_{k+1}^{p} = A(x_{k+1}^{p})$,
we also use the QR decomposition (pp. 247-248, \cite{GV2013}) to factorize it as
$(A_{k+1}^{p})^{T} = Q_{k}^{p}R_{k}^{p}$, where $Q_{k}^{p} \in \Re^{n \times m}$
satisfies $(Q_{k}^{p})^{T}Q_{k}^{p} = I$ and $R_{k}^{p} \in \Re^{m \times m}$ is
an upper triangle matrix. Then, we solve the linear system  \eqref{GNEWTONM}
as follows:
\begin{align}
     (R_{k}^{p})^{T} d_{k}^{c} = - c(x_{k+1}^{p}),  \;
     s_{k}^{c} = Q_{k}^{p} d_{k}^{c},  \;    x_{k+1} = x_{k+1}^{p} + s_{k}^{c}.  \label{CSTEP}
\end{align}

\vskip 2mm

We define the quadratic model $q_{k}(s)$ as follows:
\begin{align}
     q_k(s) = f(x_{k}) + s^{T}g_{k} + \frac{1}{2}s^{T}B_{k}s,
     \label{QOAM}
\end{align}
where $g_{k} = \nabla f(x_{k})$ and $B_{k} = P(x_{k})\nabla^{2} f(x_{k})P(x_{k})$
or its quasi-Newton approximation. Similarly to the trust-region subproblem of
the null space method (pp. 571-574, \cite{SY2006}), we decompose the predicted
reduction into two parts:
\begin{align}
   \textit{Hpred}_{k} = q_{k}(0) - q_{k}(s_{k}^{p}),
    \; \textit{Vpred}_{k} = q_{k}(s_{k}^{p}) - q_{k}(s_{k}^{p}+s_{k}^{c}). \label{HVPREDK}
\end{align}
Thus, in order to ensure the global convergence of the regularization continuation
method \eqref{TRCM}-\eqref{CSTEP}, the accepted prediction step $s_{k}^{p}$ and
the correction step $s_{k}^{c}$ need to satisfy the following condition:
\begin{align}
   \|s_{k}^{c}\| \le \theta_{1} \|s_{k}^{p}\|,    \label{VPRLEHPR}
\end{align}
where $\theta_{1}$ is a large positive constant such as $10^{6}$.
\begin{remark} The time step $\Delta t_k$ of the  regularization continuation
method \eqref{TRCM}-\eqref{GNEWTONM} is not restricted by the numerical stability.
Therefore, the large time step $\Delta t_{k}$ can be adopted in the steady-state
phase such that the  regularization continuation method \eqref{TRCM}-\eqref{GNEWTONM}
mimics the projected Newton method near the stationary point $x^{\ast}$ and it
has the fast rate of convergence. The most of all, the new step
$\alpha_{k} = \Delta t_{k}/(\Delta t_{k} + 1)$ is favourable to adopt the
trust-region updating strategy for adjusting the time step $\Delta t_{k}$.
Consequently, the regularization continuation method \eqref{TRCM}-\eqref{GNEWTONM}
accurately follows the trajectory of the  regularization flow \eqref{TRODGF} in
the transient-state phase and achieves the fast rate of convergence near the
stationary point $x^{\ast}$. \qed
\end{remark}
\subsection{The adaptive step control}

% \vskip 2mm

Another issue is how to adaptively adjust the time step $\Delta t_k$
at every iteration. We borrow the adjustment technique of the trust-region
method due to its robustness and its fast rate of convergence
\cite{CGT2000,Yuan2015}. According to the structure-preserving property of
the  regularization continuation method \eqref{TRCM}-\eqref{GNEWTONM},
$x_{k+1}$ will be expected to preserve the feasibility by computing the
correction step $s_{k}^{c}$ of equation \eqref{CSTEP}. That is to say,
$x_{k+1}$ satisfies $c(x_{k+1}) \approx 0$. Therefore, we use the objective
function $f(x)$ instead of the nonsmooth penalty function $f(x) + \sigma \|c(x)|_{1}$ as
the merit function. Similarly to the stepping-time scheme of the ODE method for
the unconstrained optimization problem \cite{Higham1999,LLT2007,LXLZ2021}
and the linearly constrained optimization problem \cite{LLS2022,LX2022},
we also need to construct a local approximation model of $f(x)$ around $x_{k}$.
Here, we adopt the quadratic function defined by equation \eqref{QOAM} as its
approximation model.

\vskip 2mm

We define the ratio of the actual reduction $Ared_{k}$ to the predicted
reduction $Pred_{k}$ as follows:
\begin{align}
    Ared_{k} \triangleq f(x_{k}) - f(x_{k} + s_{k}), \;
    Pred_{k} \triangleq q_k(0) - q_k(s_{k}),
    \; \rho_k = \frac{Ared_{k}}{Pred_{k}}.     \label{MRHOK}
\end{align}
Thus, we accept the trial step $s_{k}$ and let $x_{k+1} = x_{k}+s_{k}$,
when $\rho_{k} \ge \eta_{a}$, $\|c(x_{k}+s_{k})\|_{\infty} \le \epsilon_{0}$
and the predicted reduction $Pred_{k} = q_{k}(0) - q_{k}(s_{k})$ satisfies the
Armijo sufficient descent condition:
\begin{align}
     Pred_{k} \ge \eta_{q} \|s_{k}^{p}\| \|p_{g_k}\|,  \label{AMGESKPGK}
\end{align}
where $p_{g_k} = P(x_k)g(x_{k}) = P(x_k)\nabla f(x_k)$ and $s_{k}^{p}$ is
computed by equation \eqref{TRCM}, and  $\eta_{a}, \; \eta_{q}$ are the
small positive constants such as $\eta_{a} = \eta_{q} = 10^{-6}$,
and $\epsilon_{0}$ is less than the tolerance error $\epsilon$ such as
$\epsilon_{0} = \frac{1}{10} \epsilon$ and $\epsilon = 10^{-6}$. Otherwise, we
discard it and let $x_{k+1} = x_{k}$.

\vskip 2mm

Consequently, we reduce or enlarged the time step $\Delta t_{k+1}$ at every
iteration according to the ratio $\rho_{k}$ defined by equation \eqref{MRHOK}.
A particular adjustment strategy is given as follows:
\begin{align}
     \Delta t_{k+1} =
        \begin{cases}
          \gamma_1 \Delta t_k, \; {\text{if} \; \left(\text{the trial step $s_{k}$
          is accepted and} \; \rho_k  \ge \eta_2 \right),} \\
          \Delta t_k, \; {\text{else if} \; \left(\text{the trial step $s_{k}$ is
          accepted and} \; \eta_1 < \rho_k  < \eta_2 \right),}\\
          \gamma_2 \Delta t_k, \; {\text{others},}
        \end{cases} \label{ADTK1}
\end{align}
where the constants are selected as $\eta_1 = 0.25, \; \gamma_1 = 2,
\; \eta_2 = 0.75, \; \gamma_2 = 0.5$ according to our numerical experiments.

\vskip 2mm

\begin{remark}
This new time-stepping scheme based on the trust-region updating strategy
has some advantages compared to the traditional line search strategy
\cite{Luo2005}. If we use the line search strategy and the damped projection
Newton method \eqref{TRNEWTON}-\eqref{XK1LST} to follow the trajectory $x(t)$
of the projected Newton flow \eqref{TRODGF}, in order to achieve the fast rate
of convergence in the steady-state phase, the time step $\alpha_{k}$ of the
damped projection Newton method is tried from 1 and reduced by half with many
times at every iteration. Since the quadratic model $f(x_{k})
+ \nabla f(x_{k})^{T}s_{k} + \frac{1}{2}s_{k}^{T}B_{k}s_{k}$ may
not approximate $f(x_{k}+s_{k})$ well in the transient-state phase, the time
step $\alpha_{k}$ will be small. Consequently, the line search strategy consumes
the unnecessary trial steps in the transient-state phase. However, the selection
scheme of the time step based on the trust-region strategy
\eqref{MRHOK}-\eqref{ADTK1} can overcome this shortcoming. \qed
\end{remark}

\subsection{The adaptively preconditioned technique}

% \vskip 2mm

For the large-scale problem, the numerical evaluation of the two-sided projection
matrix $P(x_{k})\nabla^{2}f(x_{k})P(x_{k})$ consumes much time.
In order to overcome this shortcoming, in the well-posed phase, we use
the Broyden-Fletcher-Goldfarb-Shanno (BFGS) quasi-Newton matrix
(see \cite{Broyden1970,Fletcher1970,Goldfarb1970,Mascarenhas2004,Shanno1970} or
pp. 222-230, \cite{NW1999}) to approximate the two-sided projection
matrix $P(x_{k})\nabla^{2} f(x_{k})P(x_k)$ of the regularization continuation
method \eqref{TRCM}-\eqref{CSTEP} as follows:
\begin{align}
    B_{k+1} = \begin{cases}
                 B_{k} + \frac{y_{k}y_{k}^{T}}{y_{k}^{T}s_{k}}
                 - \frac{B_{k}s_{k}s_{k}^{T}B_{k}}{s_{k}^{T}B_{k}s_{k}}, \;
                 \text{if} \; y_{k}^{T}s_{k} > 0,  \\
                 B_{k}, \; \text{otherwise},
              \end{cases}
         \label{BFGS}
\end{align}
where $y_{k} = P(x_{k+1}) g(x_{k+1}) - P(x_{k}) g(x_{k}),
\; s_{k} = x_{k+1} - x_{k}$ and $B_{0} = I$.

\vskip 2mm

The BFGS updating matrix $B_{k}$ has some nice properties. For example, $B_{k+1}$ is
symmetric positive definite when $B_{k}$ is symmetric positive definite
and $B_{k+1}$ is updated by the BFGS formula \eqref{BFGS} (p. 199, \cite{NW1999}).
For the large-scale problem, it consumes much time to solve the large-scale
system \eqref{TRCM} of linear equations. In order to save the computational time
of solving linear equations, we obtain the inverse $(\sigma_{k}I + B_{k})^{-1}$ of
the regularization matrix $(\sigma_{k}I + B_{k})$ by using the following
Sherman-Morrison-Woodbury formula (p. 17, \cite{SY2006}):
\begin{align}
    (B + SV^{T})^{-1} = B^{-1} - B^{-1}S (I + V^{T}B^{-1}S)^{-1} V^{T}B^{-1},
    \label{SMWEQ}
\end{align}
where $B$ is an $n \times n$ nonsingular matrix and $S, \, V$ are two $p \times n$
matrices. Generally speaking, $p$ is less more than $n$. We denote matrices $S_{k}$
and $V_{k}$ as follows:
\begin{align}
   S_{k} & = \begin{cases}
              \left[S_{k-1}, \,  \frac{y_{k}}{(s_{k}^{T}y_{k})^{1/2}}, \,
              \frac{m_{k}}{(s_{k}^{T}m_{k})^{1/2}} \, \right], \;
               \text{if} \; y_{k}^{T}s_{k} > 0, \\
               \left[S_{k-1}, \, 0, \, 0 \right],  \; \text{otherwise}
           \end{cases}        \label{DEFSK} \\
   V_{k} & = \begin{cases}
               \left[V_{k-1}, \,  \frac{y_{k}}{(s_{k}^{T}y_{k})^{1/2}}, \,
               -\frac{m_{k}} {(s_{k}^{T}m_{k})^{1/2}} \, \right], \;
               \text{if} \; y_{k}^{T}s_{k} > 0, \\
               \left[V_{k-1}, \, 0, \, 0 \right],  \; \text{otherwise},
           \end{cases} \label{DEFVK}
\end{align}
where $m_{k} = B_{k}s_{k} = s_{k} + S_{k-1}(V_{k-1}^{T}s_{k})$. Then, from
equations \eqref{DEFSK}-\eqref{DEFVK}, the BFGS quasi-Newton formula
\eqref{BFGS} can be written as
\begin{align}
    B_{k+1} = I + S_{k}V_{k}^{T}.          \label{LBFGSSV}
\end{align}
By using the Sherman-Morrison-Woodbury formula \eqref{SMWEQ}, from equation
\eqref{LBFGSSV}, we obtain
\begin{align}
   (\sigma_{k}I + B_{k})^{-1} = \frac{1}{1 + \sigma_{k}}
    \left(I - S_{k-1}\left((1 + \sigma_{k})I
    + V_{k-1}^{T}S_{k-1}\right)^{-1}V_{k-1}^{T}\right). \label{INVLBFGS}
\end{align}

\vskip 2mm

\begin{remark}
It is worthwhile to discuss whether matrix $I + V^{T}B^{-1}S$ in equation
\eqref{SMWEQ} is nonsingular or not. Actually, when matrix $B$ is symmetric
positive definite and matrix $B+ SV^{T}$ is positive definite (that is to say,
all eigenvalues of $B+SV^{T}$ are greater than 0), $I + V^{T}B^{-1}S$ is
nonsingular.
\end{remark}
\proof  Since $B$ is symmetric positive definite, it can be decomposed as
$B = LL^{T}$ by the Cholesky factorization (p. 163, \cite{GV2013}), where
$L$ is a nonsingular matrix. Thus, we have
$B + SV^{T} = L(I + (L^{-1}S)(V^{T}L^{-T}))L^{T}$.
By combining the positive definiteness of $B + SV^{T}$, we know that
$I + (L^{-1}S)(V^{T}L^{-T})$ is positive definite. Therefore, all eigenvalues
of $(L^{-1}S)(V^{T}L^{-T})$ are greater than -1. Furthermore, the eigenvalues
of $(L^{-1}S)(V^{T}L^{-T})$ equal those of $(V^{T}L^{-T})(L^{-1}S)$ except for
their zero eigenvalues. Consequently, all eigenvalues of $(V^{T}L^{-T})(L^{-1}S)$
are greater than -1. Namely, all eigenvalues of $V^{T}B^{-1}S$ are greater than
-1. Therefore, $I + V^{T}B^{-1}S$ is a nonsingular matrix. \eproof

\vskip 2mm

According to our numerical experiments \cite{LLS2022,LX2022}, the regularization
continuation method \eqref{TRCM}-\eqref{CSTEP} with the BFGS updating
formula \eqref{DEFSK}-\eqref{INVLBFGS} works well for most problems and the
objective function decreases very fast in the well-posed phase. However, for
the ill-posed problems, the regularization continuation method
\eqref{TRCM}-\eqref{CSTEP} with the BFGS updating formula
\eqref{DEFSK}-\eqref{INVLBFGS} will approach the stationary solution
$x^{\ast}$ very slow in the ill-posed phase. Furthermore, it fails to get close
to the stationary solution $x^{\ast}$ sometimes.

\vskip 2mm

In order to improve the robustness of the  regularization continuation method
\eqref{TRCM}-\eqref{CSTEP}, we set $B_{k} = P(x_k)\nabla^{2}f(x_{k})P(x_k)$
in the ill-posed phase. Now, the problem is how to automatically identify the
ill-posed phase and switch to $B_{k} = P(x_k)\nabla^{2} f(x_{k})P(x_k)$ from the
BFGS updating formula \eqref{BFGS}. Here, we adopt the simple criterion. That is to say, we regard
that the  regularization continuation method \eqref{TRCM}-\eqref{CSTEP} is in
the ill-posed phase once there exists the time step $\Delta t_{K} \le 10^{-3}$.

\vskip 2mm

In the ill-posed phase,  the computational time of the projected Hessian matrix
$P(x_k)\nabla^{2}f(x_{k})P(x_k)$ is heavy if we update it at every iteration.
In order to save the computational time of the Hessian evaluation
$\nabla^{2} f(x_{k})$, we set
$B_{k+1} = B_{k}$ when $q_{k}(0) - q_{k}(s_{k})$ approximates
$f(x_{k}) - f(x_{k}+s_{k})$ well, where the approximation model $q_{k}(s)$
is defined by equation \eqref{QOAM}. Otherwise, we update
$B_{k+1} =  P(x_{k+1})\nabla^{2}f(x_{k+1})P(x_{k+1})$ in the ill-posed
phase \cite{LX2021,LX2022}. In the ill-posed phase, a practice updating strategy
is give by
\begin{align}
     B_{k+1}
       = \begin{cases}
            B_{k},  \; \text{if} \; |1- \rho_{k}| \le \eta_{1}, \\
             P(x_{k+1})\nabla^{2}f(x_{k+1})P(x_{k+1}), \; \text{otherwise},
          \end{cases} \label{UPDJK1}
\end{align}
where $\rho_{k}$ is defined by equations \eqref{QOAM}-\eqref{MRHOK} and
$\eta_{1} = 0.25$.

\subsection{Finding an initial feasible point} \label{SUBSECFIP}

\vskip 2mm

For the regularization continuation method based on the ODE system \eqref{TRODGF},
one of the important issues is how to find an initial feasible point $x_{0}$.
That is to say, the initial point $x_{0}$ needs to satisfy $c(x_{0}) = 0$. Here,
we use the generalized continuation Newton method \cite{LX2021} to solve the
under-determined system $c(x) = 0$ due to its robustness and efficiency. For
convenience, we give the rough description of the generalized continuation Newton
method. Its detailed description and its convergence analysis can be found in
reference \cite{LX2021}.

\vskip 2mm

For the under-determined system $c(z) = 0$, we construct the generalized Newton flow
\begin{align}
    \frac{dz(\tau)}{d \tau} = - A(z)^{\dagger}c(z),  \; z(\tau_{0}) = z_{0},
    \label{GNFLOW}
\end{align}
where $A(z)^{\dagger}$ is the Moore-Penrose generalized inverse of the Jacobian
matrix $A(z)$ (p. 11, \cite{SY2006} or p. 290, \cite{GV2013}). Then, we
construct the generalized continuation Newton method \cite{LX2021} to follow the
generalized Newton flow \eqref{GNFLOW} and obtain its steady-state solution
$z^{\ast}$ as follows:
\begin{align}
    \Delta z_{j} = - \frac{\Delta \tau_{j}}
    {1 + \Delta \tau_{j}} A_{j}^{\dagger}c(z_j), \;  z_{j+1} = z_{j} + \Delta z_{j},
    \label{GCNM}
\end{align}
where $A_{j}$ equals $A(x_{j})$ or its approximation according to the following
ratio
\begin{align}
  r_{j} = \frac{\|c(z_{j})\|-\|c(z_{j}+ \Delta z_{j})\|}
      {\|c(z_{j})\| - \|c(z_{j}) + A_{j} \Delta z_{j}\|} =
  \frac{\|c(z_{j})\|-\|c(z_{j}+ \Delta z_{j})\|}
   {(\Delta \tau_{j}/(1+\Delta \tau_{j}))\|c(z_j)\|}.  \label{RATIOJ}
\end{align}
In other words, we update the Jacobian matrix $A_{j+1}$ according to the following
strategy:
\begin{align}
     A_{j+1} = \begin{cases}
                 A_{j},  \; \text{if} \; |1- r_{j}| \le \eta_{1}, \\
                 A(z_{j+1}), \; \text{otherwise}.
               \end{cases} \label{UPDAJ1}
\end{align}

\vskip 2mm

For a real-world problem, $A_{j}A_{j}^{T}$ may be ill-conditioned. Thus, the Cholesky
decomposition may fail to solve the linear system \eqref{GCNM} for the large-scale
problem. In order to improve its robustness, we use the QR decomposition to solve
it as follows:
\begin{align}
    A_{j}^{T} = U_{j}W_{j}, \; W_{j}^{T} \Delta z_{j}^{m} = - c(z_{j}), \;
    \Delta z_{j}^{N} = U_{j} \Delta z_{j}^{m}, \;
    \Delta z_{j} = \frac{\Delta \tau_{j}}{1 + \Delta \tau_{j}} \Delta z_{j}^{N},
    \label{GCNMQR}
\end{align}
where $U_{j} \in \Re^{n \times m}$ satisfies $U_{j}^{T}U_{j} = I$ and
$W_{j} \in \Re^{m \times m}$ is an upper triangle matrix. In practice, in order
to save the computational time of decomposing the matrix $A_{j+1}$ when $A_{j}$
performs well, i.e. $|1-r_{j}| \le \eta_{1}$, according to the updating formula
\eqref{UPDAJ1}, we set $W_{j+1} = W_{j}$ and $U_{j+1} = U_{j}$ in equation
\eqref{GCNMQR}.

\vskip 2mm

The adaptive time step $\Delta \tau_{j}$ of the generalized continuation
Newton method \eqref{GCNM} is adjusted by the following trust-region updating
strategy:
\begin{align}
   \Delta \tau_{j+1} =
     \begin{cases}
    \gamma_1 \Delta \tau_j, &{\text{if} \;  \left|1- r_j \right| \le \eta_1,}\\
    \Delta \tau_j, & {\text{else if} \; \eta_1 < \left|1 - r_j \right| < \eta_2,}\\
    \gamma_2 \Delta \tau_j, &{\text{others},}
    \end{cases} \label{TAUSK1}
\end{align}
where the constants are selected as $\gamma_{1} = 2, \; \gamma_{2} = 0.5, \;
\eta_{1} = 0.25, \; \eta_{2} = 0.75$, according to our numerical experiments.

\vskip 2mm

For a real-world problem, the analytical Jacobian $A(z_{j})$ may not be
offered. Thus, in practice, we replace the Jacobian matrix $A(z_{j})$
with its difference approximation as follows:
\begin{align}
     A(z_{k}) \approx
    \left[\frac{c(z_{j} + \epsilon e_{1}) - c(z_{j})}{\epsilon}, \,
    \ldots, \, \frac{c(z_{j} + \epsilon e_{n}) - c(z_{j})}{\epsilon}\right],
    \label{NUMJA}
\end{align}
where the elements of $e_{i}$ equal zeros except for the $i$-th element which
equals 1, and the parameter $\epsilon$ can be selected as $10^{-6}$ according to
our numerical experiments.

\vskip 2mm

According to the above descriptions, we give the detailed implementation of
the generalized continuation Newton method with the trust-region updating strategy
(GCNMtr) to find an initial feasible point $z^{\ast}$ of $c(z) = 0$ in  Algorithm
\ref{ALGGCNMTR}.

\vskip 2mm

\begin{algorithm}
   \renewcommand{\algorithmicrequire}{\textbf{Input:}}
   \renewcommand{\algorithmicensure}{\textbf{Output:}}
   \caption{Generalized continuation Newton methods and the trust-region
   updating strategy for the under-determined system  (The GCNMTr method)}
   \label{ALGGCNMTR}
   \begin{algorithmic}[1]
      \REQUIRE ~~ \\
      Function $c: \; \Re^{n} \to \Re^{m}, \; m \le n$,
      the initial point  $z_0$ (optional), and the tolerance $\epsilon$ (optional).
	  \ENSURE ~~ \\
      An approximation solution $z^{\ast}$ of nonlinear equations.
      \STATE Set the default $z_0 = \text{ones} (n, \, 1)$  and
      $\epsilon = 10^{-7}$, when $z_0$ or $\epsilon$ is not provided.
      \STATE Initialize the parameters: $\eta_{a} = 10^{-6}, \; \eta_1 = 0.25, \;
      \gamma_1 =2, \; \eta_2 = 0.75, \; \gamma_2 = 0.5, \; \text{maxit} = 400$.
      \STATE Set $\Delta \tau_0 = 10^{-2}$, flag\_success\_trialstep = 1, $\text{itc} = 0, \; j = 0$.
      \STATE Evaluate $c_{j} = c(z_{j})$.
      \STATE Set $r_{j-1}  = 0$.
      \WHILE{(itc $<$ maxit)}
         \IF{(flag\_success\_trialstep == 1)}
              \STATE Set itc = itc + 1.
              \STATE Compute $\text{Res}_{j} = \|c_{j}\|_{\infty}$.
              \IF{($\text{Res}_{j} < \epsilon$)}
                 \STATE break;
              \ENDIF
              \IF{$(|1-r_{j-1}| > \eta_1)$}
                 \STATE Evaluate $A_j = A(z_{j})$ from equation \eqref{NUMJA}.
                 \STATE Use the QR decomposition $[\bar{U}_{j},\bar{W}_{j}] = \text{qr}(A_{j}^{T})$
                    and set $U_{j} = \bar{U}_{j}(:, 1:m), \; W_{j} = \bar{W}_{j}(1:m, 1:m)$.
              \ELSE
                 \STATE Set $U_{j} = U_{j-1}$, $W_{j} = W_{j-1}$.
              \ENDIF
              \STATE By solving $W_{j}^{T} \Delta z_{j}^{m} = - c_{j}$ and
              $\Delta z_{j}^{N} = U_{j} \Delta z_{j}^{m}$, we obtain the Newton
              step $\Delta z_{j}^{N}$.
         \ENDIF
         \STATE Set $\Delta z_{j} = {\Delta \tau_{j}}/{(1+\Delta \tau_{j})} \, \Delta z_{j}^{N}, \;
         z_{j+1} = z_{j} + \Delta z_{j}$.
         \STATE Evaluate $c_{j+1} = c(z_{j+1})$.
         \IF {$\left(\|c_{j}\| < \|c_{j+1}\|\right)$}
           \STATE Set $r_{j} = -1$;
         \ELSE
           \STATE Compute the ratio $r_{j}$ from equation \eqref{RATIOJ}.
         \ENDIF
         \STATE Adjust the time step size $\Delta \tau_{j+1}$ according to the
         trust-region updating strategy \eqref{TAUSK1}.
         \IF{$(r_{j} \ge \eta_{a})$}
               \STATE Accept the trial point $z_{j+1}$. Set flag\_success\_trialstep = 1.
           \ELSE
               \STATE Set $z_{j+1}  = z_{j}$, $c_{j+1} = c_{j}$, $\Delta z_{j+1}^{N} = \Delta z_{j}^{N}$,
               flag\_success\_trialstep = 0.
           \ENDIF
         \STATE Set $j \longleftarrow j+1$.
      \ENDWHILE
   \end{algorithmic}
\end{algorithm}

\vskip 2mm

By combining the discussions of the previous sections and Algorithm
\ref{ALGGCNMTR}, we give the detailed implementation of the  regularization
continuation method and the trust-region updating strategy for the optimization
problem \eqref{NLEQOPT}  with nonlinear equality constraints in Algorithm
\ref{ALGRCM}.

\begin{algorithm}
	\renewcommand{\algorithmicrequire}{\textbf{Input:}}
	\renewcommand{\algorithmicensure}{\textbf{Output:}}
    \newcommand{\algorithmicbreak}{\textbf{break}}
    \newcommand{\BREAK}{\STATE \algorithmicbreak}
	\caption{The regularization continuation method for optimization problems
     with nonlinear equality constraints (the Rcm method)}
    \label{ALGRCM}	
	\begin{algorithmic}[1]
		\REQUIRE ~~ the objective function $f: \; \Re^{n} \to \Re$, the equality
        constraints $c(x)  = 0, \; c: \; \Re^{n} \to \Re^{m}$,
        the initial point $z_0$ (optional), the terminated parameter
        $\epsilon$ (optional).
		\ENSURE ~~
        the optimal approximation solution $x^{\ast}$.
        \STATE If $z_{0}$ or $\epsilon$ is not provided, we set
        $z_0 = \text{ones}(n, \, 1)$ or $\epsilon = 10^{-6}$.
        We let $\epsilon_{0} = \frac{1}{10} \epsilon$.
        \STATE Initialize parameters: $\eta_{a} = 10^{-6}, \; \eta_{m} = 10^{-6}$,
        $\eta_1 = 0.25, \; \gamma_1 =2, \; \eta_2 = 0.75, \; \gamma_2 = 0.5$,
        $\sigma_{0} = 10^{-5}, \Delta t_K = 10^{-3}$, max\_itc = 300. Set $\Delta t_0 = 10^{-2}$,
        flag\_illposed\_phase = 0, flag\_success\_trialstep = 1,
        $s_{-1} = 0, \; y_{-1} = 0, \; \rho_{-1} = 0, \; B_{0} = I, \; H_{0} = I$,
        $Q_{-1}^{b} = I, \; R_{-1}^{b} = I$, $S_{-1} = \text{zeros}(n, \, 1)$,
        $V_{-1} = \text{zeros}(n, \, 1)$, itc = 0.
        \STATE Use the generalized continuation Newton method (the GCNMTr
        method, Algorithm \ref{ALGGCNMTR}) to find an initial feasible point from
        $z_{0}$ and denote this initial feasible point as $x_{0}$.
        \STATE Set $k = 0$. Evaluate $f_0 = f(x_0)$ and $g_0 = \nabla f(x_0)$.
        Evaluate $A_{0} = A(x_{0})$ from equation \eqref{NUMJA}.
        \STATE Factorize $A_{0}^{T}$ by the QR decomposition
        $[\bar{Q}_{0},\bar{R}_{0}] = \text{qr}(A_{0}^{T})$.
        \STATE Set $Q_{0} = \bar{Q}_{0}(:, 1:m)$, $R_{0} = \bar{R}_{0}(1:m, 1:m)$.
        \STATE Compute the projected gradient $p_{g_{0}} = P_{0}g_{0}$, where
        $P_{0}$ is computed by equation \eqref{SPROMAT}.
        \WHILE{$\left(\left(\|p_{g_k}\|_{\infty}> \epsilon\right) \; \text{and} \;
         (\text{itc} < \text{max\_itc})\right)$}
           \STATE itc = itc + 1;
           \IF{$\Delta t_{k} < \Delta t_{K}$}
              \STATE Set flag\_illposed\_phase = 1.
           \ENDIF
           \IF{(flag\_illposed\_phase == 0)}
               \IF{(flag\_success\_trialstep == 1)}
                    \IF{$k > 0$}
                        \STATE Solve equations \eqref{DEFSK}-\eqref{DEFVK} to obtain
                               matrices $S_{k-1}$ and $V_{k-1}$.
                    \ELSE
                        \STATE Set $S_{-1} = \text{zeros}(n, \, 1)$
                               and $V_{-1} = \text{zeros}(n, \, 1)$.
                    \ENDIF
                   \STATE Solve equations \eqref{TRCM} and \eqref{INVLBFGS} to obtain \\
                   $d_{k} = -\left(\left(p_{g_k} - S_{k-1}\left(\left((1 + \sigma_{0}/\Delta t_{k})I
                   + V_{k-1}^{T}S_{k-1}\right)^{-1}(V_{k-1}^{T}p_{g_{k}})\right)\right)\right)
                   /({1 + \sigma_{0}/\Delta t_{k}}).$
               \ENDIF
           \ELSE
               \IF{$(|\rho_{k-1} - 1| > \eta_1)$}
                   \STATE Evaluate  $H_{k} = P_{k}\nabla^{2}f(x_{k})P_{k}$ from equation \eqref{NUMHESS}.
                   \STATE Set $B_{k} = ({\sigma_{0}}/{\Delta t_{k}})I + H_{k}$ and use the QR decomposition
                   $B_{k}= Q_{k}^{b} R_{k}^{b}$.
                \ELSE
                   \STATE $Q_{k}^{b} = Q_{k-1}^{b}, \; R_{k}^{b} = R_{k-1}^{b}$;
               \ENDIF
               \STATE Solve the linear system $R_{k}^{b}d_{k} = - (Q_{k}^{b})^{T}p_{g_k}$ to obtain $d_{k}$.
          \ENDIF
          \STATE Set $s_{k}^{p} = \frac{\Delta t_{k}}{1+\Delta t_{k}}(P_{k}d_{k})$ and
          $x_{k+1}^{p} = x_{k} + s_{k}^{p}$. Evaluate $c_{k+1}^{p} = c(x_{k+1}^{p})$.
          \STATE Solve the linear systems $(R_{k})^{T} d_{k}^{c} = - c_{k+1}^{p}, \;
          s_{k}^{c} = Q_{k} d_{k}^{c}$.
          \STATE Set $x_{k+1} = x_{k+1}^{p} + s_{k}^{c}, \; c_{k+1} = c(x_{k+1})$.
          \IF{($\|c_{k+1}\|_{\infty} > \epsilon_{0}$)}
             \STATE Evaluate $A_{k+1}^{p} = A(x_{k+1}^{p})$ from equation \eqref{NUMJA}.
             \STATE Solve the linear system \eqref{CSTEP} to obtain $x_{k+1}$ and evaluate $c_{k+1} = c(x_{k+1})$.
          \ENDIF
          \STATE Evaluate $f_{k+1} = f(x_{k+1})$ and compute the ratio $\rho_{k}$
          from equations \eqref{QOAM}-\eqref{MRHOK}.
          \IF{($\rho_k \ge \eta_{a}, \; \|c_{k+1}\|_{\infty} \le \epsilon_{0}$
            and $s_{k}$ satisfies the sufficient descent condition \eqref{AMGESKPGK})}
             \STATE Set flag\_success\_trialstep = 1. Evaluate $g_{k+1} = \nabla f(x_{k+1})$.
             \STATE Evaluate $A_{k+1} = A(x_{k+1})$ from equation \eqref{NUMJA}.
             \STATE Factorize $A_{k+1}^{T}$ by the QR decomposition $[\bar{Q}_{k+1},\bar{R}_{k+1}] = \text{qr}\left(A_{k+1}^{T}\right)$.
             \STATE Set $Q_{k+1} = \bar{Q}_{k+1}(:, 1:m)$, $R_{k+1} = \bar{R}_{k+1}(1:m, 1:m)$.
             \STATE Set $p_{g_{k+1}} = P_{k+1}g_{k+1}$,where $P_{k+1}$ is computed by equation \eqref{SPROMAT}.
             \STATE Set $s_{k} = x_{k+1} - x_{k}$ and $y_{k} = p_{g_{k+1}} - p_{g_{k}}$.
          \ELSE
            \STATE Set flag\_success\_trialstep = 0.
            \STATE $x_{k+1} = x_{k}, \; f_{k+1} = f_{k}$,
            $p_{g_{k+1}} = p_{g_{k}}, \; g_{k+1} = g_{k}, \; H_{k+1} = H_{k}, \;  d_{k+1} = d_{k}.$
          \ENDIF
          \STATE Adjust the time step $\Delta t_{k+1}$ based on the
          trust-region updating strategy \eqref{ADTK1}.
          \STATE Set $k \leftarrow k+1$.
        \ENDWHILE
	\end{algorithmic}
\end{algorithm}

\section{Algorithm Analysis}

In this section, we analyze the global convergence of the  regularization
continuation method \eqref{TRCM}-\eqref{CSTEP} and the adaptive time step
control for the nonlinearly equality-constrained optimization problem
(i.e. Algorithm \ref{ALGRCM}). Similarly to the result of the trust-region
method for the unconstrained optimization problem \cite{Powell1975} and the
continuation method for linearly constrained optimization problem
\cite{LLS2022}, we have the following lower bound estimation of
$q_{k}(0)-q_{k}(s_{k}^{p})$.

\vskip 2mm

\begin{lemma} \label{LEMBSOAM}
Assume that the quadratic model $q_{k}(s)$ is defined by equation \eqref{QOAM}
and $d_{k}$ is the solution of equation \eqref{TRCM}. Furthermore, we suppose that
the time step $\Delta t_k$ satisfies
\begin{align}
    \left(\frac{\sigma_{0}}{\Delta t_{k}} \, I + B_k\right) \succ 0 \;
    \text{and} \; \left(\frac{\sigma_{0}}{\Delta t_{k}}
    \, I + B_{k}- P_{k}^{T}B_{k}P_{k} \right) \succeq 0,  \label{PSDASS}
\end{align}
where $P_{k} = P(x_{k})$ is the projection matrix defined by equation \eqref{PROMAT}.
Then, we have
\begin{align}
    q_{k}(0) - q_{k}(P_{k}d_{k}) \ge \frac{1}{2} \left\|p_{g_{k}} \right\|
    \min \left\{\left\|P_{k}d_{k}\right\|, \; \frac{\|p_{g_k}\|}{3\|B_{k}\|}\right\},
    \label{PLBREDST}
\end{align}
where $p_{g_k} = P_{k}g_{k}$.
\end{lemma}
\proof Let $\tau_{k} = \sigma_{0}/\Delta t_{k}$. From equation \eqref{TRCM}, we obtain
\begin{align}
    q_{k}(0) &- q_{k}(P_{k}d_{k}) = - \frac{1}{2}\left(d_{k}^{T}P_{k}^{T}B_{k}P_{k}d_{k}\right)
     - (P_{k}g_{k})^{T}d_{k}  \nonumber \\
    & = - \frac{1}{2}\left(d_{k}^{T}P_{k}^{T}B_{k}P_{k}d_{k}\right)
    + p_{g_k}^{T} \left(\tau_{k}I + B_{k}\right)^{-1}p_{g_k}  \nonumber \\
    & =  \frac{1}{2}\left(p_{g_k}^{T} \left(\tau_{k}I + B_{k}\right)^{-1}p_{g_k}
    + d_{k}^{T}\left(-P_{k}^{T}B_{k}P_{k} + \tau_{k}I + B_{k}\right)d_{k}\right).
    \label{DKGKDK}
\end{align}
We denote $\mu_{\min}\left(B_{k}-P_{k}^{T}B_{k}P_{k}\right)$ as the smallest eigenvalue of
matrix $\left(B_{k}-P_{k}^{T}B_{k}P_{k}\right)$, and set
\begin{align}
    \tau_{lb} = \min\left\{0, \; \mu_{\min} \left(B_{k}-P_{k}^{T}B_{k}P_{k}\right)\right\}.
    \label{LOWBTAU}
\end{align}
From equations \eqref{PSDASS}, \eqref{DKGKDK}-\eqref{LOWBTAU} and the bound
on the eigenvalues of matrix $(\tau_{k}I + B_{k})^{-1}$, we obtain
\begin{align}
    & q_{k}(0) - q_{k}(P_{k}d_{k}) \ge \frac{1}{2}\left(p_{g_k}^{T}
    \left(\tau_{k}I + B_{k}\right)^{-1}p_{g_k}
    + \left(\tau_{k} + \tau_{lb}\right)\left\|d_{k}\right\|^{2}\right)
    \nonumber \\
    & \quad \ge \frac{1}{2}\left(\frac{\left\|p_{g_{k}}\right\|^{2}}{\tau_{k}+\left\|B_{k}\right\|}
    + \left(\tau_{k} + \tau_{lb}\right)\left\|d_{k}\right\|^{2} \right).
    \label{LBEQK}
\end{align}
In the above second inequality, we use the property $|\mu_{i}(B_{k})|\le \|B_{k}\|$,
where $\mu_{i}(B_{k})$ is an eigenvalue of matrix $B_{k}$.

\vskip 2mm

Now we consider the properties of the function
\begin{align}
   \varphi(\tau) \triangleq \tau \left\|d_k\right\|^{2}
   + {\left\|p_{g_k}\right\|^2}/\left(\tau - \tau_{lb}+\left\|B_k \right\|\right).
   \label{VLF}
\end{align}
From equation \eqref{VLF}, we have $\varphi^{''}(\tau) = 2\|p_{g_k}\|^2
/\left(\tau - \tau_{lb} + \|B_{k}\|\right)^3  \ge 0$ when
$(\tau - \tau_{lb}$ $+ \|B_{k}\|) > 0$. Thus, the function $\varphi(\tau)$
attains its minimum $\varphi(\tau_{\min})$ when
$\tau_{\min}$ satisfies $\varphi^{'}(\tau_{\min})=0$ and
$\tau \ge -(-\tau_{lb}+\|B_k\|)$. Namely, we have
\begin{align}
  \varphi(\tau_{\min}) = 2\|p_{g_k}\| \|d_k\|
  + \left(\tau_{lb}-\|B_k\|\right) \|d_k\|^{2},   \label{MINV}
\end{align}
where
\begin{align}
    \tau_{\min} = \frac{\|p_{g_k}\|}{\|d_k\|} + \tau_{lb}-\|B_k\|. \label{MINLD}
\end{align}

\vskip 2mm

We prove the property \eqref{PLBREDST} via separately considering $\tau_{\min} \ge 0$ or
$\tau_{\min} < 0$ as follows.

\vskip 2mm

(i) When $\left(\|p_{g_k} \|/\|d_k\| + \left(\tau_{lb}-\|B_k\|\right)\right) \ge 0$,
from equation \eqref{MINLD}, we have $\tau_{\min} \ge 0$. From the assumption
\eqref{PSDASS} and the definition \eqref{LOWBTAU} of $\tau_{lb}$, we have
$\tau_k \ge -\tau_{lb}$. Thus, from equations \eqref{LBEQK}--\eqref{MINLD},
we obtain
\begin{align}
    & q_k(0)   - q_k(P_{k}d_{k}) \ge \frac{1}{2}
    \left((\tau_{k}+\tau_{lb}) \|d_k\|^{2} + \frac{\|p_{g_k}\|^2}{\tau_{k} + \|B_k\|}\right)
    =  \frac{1}{2} \varphi(\tau_{k}+\tau_{lb})
    \ge \frac{1}{2} \varphi(\tau_{\min}) \nonumber \\
    & \quad = \frac{1}{2} \left(\|p_{g_k}\| \|d_k\|
    + \left(\|p_{g_k}\| \|d_k\|+\left(\tau_{lb}-\|B_k\|\right)\|d_k\|^{2}\right)\right)
   \ge \frac{1}{2}\|p_{g_k}\| \|d_k\|.   \label{MINQCOV}
\end{align}

\vskip 2mm

(ii) The other case is $\left(\|p_{g_k} \|/\|d_k\| + \left(\tau_{lb}-\|B_k\|\right)\right) < 0$.
In this case, from equation \eqref{MINLD}, we have $\tau_{\min} < 0$. It is not
difficult to verify that $\varphi(\tau)$ is monotonically increasing
when $\tau \ge 0$ and $\tau_{\min}<0$. From the definition \eqref{LOWBTAU} of $\tau_{lb}$ and
the property \eqref{MATNP}, we have
\begin{align}
    |\tau_{lb}| &\le \left|\mu_{\min} \left(B_{k}-P_{k}^{T}B_{k}P_{k}\right)\right|
    \le \|B_{k} - P_{k}^{T}B_{k}P_{k}\| \le \|B_{k}\|+\|P_{k}^{T}B_{k}P_{k}\| \nonumber \\
    & \le \|B_{k}\|+\|P_{k}^{T}\|\|B_{k}\|\|P_{k}\|= 2\|B_k\|. \nonumber
\end{align}
By using this property and the monotonicity of $\varphi(\tau)$, from equations
\eqref{LBEQK}-\eqref{VLF}, we obtain
\begin{align}
     & q_{k}(0)  - q_{k}(P_{k}d_{k})
    \ge \frac{1}{2}\left((\tau_{k}+\tau_{lb}) \|d_k\|^{2}
    + \frac{\|p_{g_k}\|^2}{\tau_k + \|B_k \|}\right) \nonumber \\
    & \quad = \frac{1}{2} \varphi(\tau_{k}+\tau_{lb}) \ge \frac{1}{2} \varphi(0)
    = \frac{1}{2(-\tau_{lb}+\|B_k\|)}\|p_{g_k}\|^2
     \ge \frac{1}{6\|B_k\|}\|p_{g_k}\|^2.   \label{MINQMI}
\end{align}

\vskip 2mm

From equations (\ref{MINQCOV})-(\ref{MINQMI}), we get
\begin{align}
     q_{k}(0) - q_{k}(P_{k}d_{k}) \ge \frac{1}{2} \|p_{g_k}\|
     \min\left\{\|d_k\|, \;  \frac{\|p_{g_k}\|}{3\|B_k\|}\right\}.
     \label{MINQATR}
\end{align}
By using the property \eqref{MATNP} of matrix $P_{k}$, we have
\begin{align}
 \|P_{k}d_{k}\| \le \|P_{k}\| \|d_{k}\| = \|d_{k}\|. \label{PROPD}
\end{align}
Therefore, from inequalities \eqref{MINQATR}-\eqref{PROPD}, we obtain the estimation
\eqref{PLBREDST}. \qed

\vskip 2mm

When $\Delta t_{k} \le \frac{\sigma_{0}}{2\|B_{k}\|}$, by using the property $\|P_{k}\| = 1$ of
the projection matrix $P_{k}$, we have
\begin{align}
    \left(\frac{\sigma_{0}}{\Delta t_{k}} \, I + B_k\right) \succ 0 \;
    \text{and} \; \left(\frac{\sigma_{0}}{\Delta t_{k}}
    \, I + B_{k}- P_{k}^{T}B_{k}P_{k} \right) \succeq 0. \label{SPDBPBCON}
\end{align}
Thus, from equation \eqref{TRCM} and Lemma \ref{LEMBSOAM}, when
$\Delta t_{k}  \le \frac{\sigma_{0}}{4\|B_{k}\|}$, we obtain the estimation of
the predicted reduction of the horizontal step as follows.

\vskip 2mm

\begin{lemma} \label{LEMHPRED}
Assume that the quadratic model $q_{k}(s)$ is defined by equation \eqref{QOAM}
and $s_{k}^{p}$ is solved by the regularization method \eqref{TRCM}. Then, when
$\Delta t_{k} \le \frac{\sigma_{0}}{4\|B_{k}\|}$, we have the lower bounded estimation
of $q_{0} - q_{k}(s_{k}^{p})$ as follows:
\begin{align}
    \textit{Hpred}_{k} = q_{k}(0) - q_{k}(s_{k}^{p})
    = q_{k}(0) - q_{k}\left(\frac{\Delta t_{k}}{1+\Delta t_{k}} P_{k}d_{k}\right)
    \ge \frac{1}{2} \|P_{k}g_{k}\| \|s_{k}^{p}\|. \label{HPREDKEST}
\end{align}
\end{lemma}
\proof We prove this result via considering two different cases, i.e.
$(s_{k}^{p})^{T}B_{k}s_{k}^{p} \le 0$ or $(s_{k}^{p})^{T}B_{k}s_{k}^{p} > 0$.
(i) When $(s_{k}^{p})^{T}B_{k}s_{k}^{p} \le 0$, from equation \eqref{QOAM} and
equation \eqref{TRCM}, we have
\begin{align}
   & q_{k}(0) - q_{k}(s_{k}^{p}) = - g_{k}^{T}s_{k}^{p}
    - \frac{1}{2}(s_{k}^{p})^{T}B_{k}s_{k}^{p} \ge - g_{k}^{T}s_{k}^{p}
    = - \frac{\Delta t_{k}}{1 + \Delta t_{k}} (P_{k}g_{k})^{T}d_{k} \nonumber \\
   &
   = \frac{\Delta t_{k}}{1 + \Delta t_{k}} (P_{k}g_{k})^{T}
   \left(\frac{\sigma_{0}}{\Delta t_{k}}I+B_{k}\right)^{-1}(P_{k}g_{k})
   \ge \frac{\Delta t_{k}}{1 + \Delta t_{k}}
   \frac{1}{\sigma_{0}/{\Delta t_{k}} + \|B_{k}\|}\|P_{k}g_{k}\|^{2}.
   \label{HPREDGEPG}
\end{align}
Furthermore, from equation \eqref{TRCM}, when
$\Delta t_{k} < \frac{\sigma_{0}}{\|B_{k}\|}$, we have
\begin{align}
    \|P_{k}g_{k}\| & =
    \left\|\left(\frac{\sigma_{0}}{\Delta t_{k}}I + B_{k}\right)d_{k} \right\|
    \ge (\sigma_{0}/{\Delta t_{k}} - \|B_{k}\|) \|d_{k}\| \nonumber \\
    &  \ge (\sigma_{0}/{\Delta t_{k}} - \|B_{k}\|) \|P_{k}d_{k}\|.
    \label{PKGKGEPKDK}
\end{align}
By substituting equation \eqref{PKGKGEPKDK} into equation \eqref{HPREDGEPG},
when $\Delta t_{k} \le \frac{\sigma_{0}}{4\|B_{k}\|}$, we obtain
\begin{align}
    & q_{k}(0) - q_{k}(s_{k}^{p}) \ge  \frac{\Delta t_{k}}{1 + \Delta t_{k}}
    \frac{\sigma_{0}/{\Delta t_{k}} - \|B_{k}\|}
    {\sigma_{0}/{\Delta t_{k}} + \|B_{k}\|}\|P_{k}g_{k}\| \|P_{k}d_{k}\|
    \nonumber \\
    &
    \ge \frac{\sigma_{0}/{\Delta t_{k}} - \|B_{k}\|}
    {\sigma_{0}/{\Delta t_{k}} + \|B_{k}\|}\|P_{k}g_{k}\| \|s_{k}^{p}\|
    \ge \frac{4\|B_{k}\| - \|B_{k}\|}
    {4\|B_{k}\| + \|B_{k}\|}\|P_{k}g_{k}\| \|s_{k}^{p}\|
    = \frac{3}{5}\|P_{k}g_{k}\| \|s_{k}^{p}\|.
    \label{HPREDNEG}
\end{align}
Thus, the result \eqref{HPREDKEST} is true for the case
$(s_{k}^{p})^{T}B_{k}s_{k}^{p} \le 0$.

\vskip 2mm (ii) When $(s_{k}^{p})^{T}B_{k}s_{k}^{p} \ge 0$, i.e.
$(P_{k}d_{k})^{T}B_{k}(P_{k}d_{k}) \ge 0$, from equation
\eqref{QOAM} and equation \eqref{TRCM}, we have
\begin{align}
   & q_{k}(0) - q_{k}(s_{k}^{p}) = - g_{k}^{T}s_{k}^{p}
    - \frac{1}{2}(s_{k}^{p})^{T}B_{k}s_{k}^{p}
    \nonumber \\
    &
    = - \frac{\Delta t_{k}}{1 + \Delta t_{k}} g_{k}^{T}(P_{k}d_{k})
    - \frac{1}{2} \left(\frac{\Delta t_{k}}{1 + \Delta t_{k}}\right)^{2}
    (P_{k}d_{k})^{T}B_{k}(P_{k}d_{k}) \nonumber \\
    &
    \ge \frac{\Delta t_{k}}{1 + \Delta t_{k}} \left(- g_{k}^{T}(P_{k}d_{k})
    - \frac{1}{2} (P_{k}d_{k})^{T}B_{k}(P_{k}d_{k}) \right)
    = \frac{\Delta t_{k}}{1 + \Delta t_{k}} (q_{k}(0) - q_{k}(P_{k}d_{k})).
   \label{HPREDGEPG2}
\end{align}
When $\Delta t_{k} \le \frac{\sigma_{0}}{4\|B_{k}\|}$, by substituting equation
\eqref{PLBREDST} and equation \eqref{PKGKGEPKDK} into equation \eqref{HPREDGEPG2},
we obtain
\begin{align}
    & q_{k}(0) - q_{k}(s_{k}^{p}) \ge \frac{1}{2} \frac{\Delta t_{k}}{1+\Delta t_{k}}
    \|P_{k}g_{k}\|
    \min \left\{\|P_{k}d_{k}\|, \; \frac{\|P_{k}g_{k}\|}{3\|B_{k}\|}\right\}
    \nonumber \\
    & \hskip 2mm
    \ge \frac{1}{2} \frac{\Delta t_{k}}{1+\Delta t_{k}}
    \|P_{k}g_{k}\|
    \min \left\{\|P_{k}d_{k}\|, \;
    \frac{(\sigma_{0}/{\Delta t_{k}} - \|B_{k}\|) \|P_{k}d_{k}\|}{3\|B_{k}\|}\right\}
    \nonumber \\
    & \hskip 2mm
    \ge \frac{1}{2} \frac{\Delta t_{k}}{1+\Delta t_{k}} \|P_{k}g_{k}\|\|P_{k}d_{k}\|
    = \frac{1}{2} \|P_{k}g_{k}\|\|s_{k}^{p}\|.
    \label{HPREDGEGK2}
\end{align}
Thus, the result \eqref{HPREDKEST} is also true for the case
$(s_{k}^{p})^{T}B_{k}s_{k}^{p} > 0$. \qed

\vskip 2mm

In the following analysis of Algorithm \ref{ALGRCM}, $f(x)$ and $c(x)$ are assumed
to satisfy Assumption \ref{ASSFUN}.

\vskip 2mm

\begin{assumption} \label{ASSFUN}
Assume that $f(x)$ is twice continuously differentiable and $c(x)$ is continuously
differentiable, that the sequence $\{x_{k}\}$ is generated by Algorithm \ref{ALGRCM}
in an open set $\mathrm{S}$, that $\nabla f(x), \; \nabla^{2} f(x), \; A(x)$ are
bounded above on $\mathrm{S}$. We also assume that $\|B_{k}\| \, (k = 0, \, 1, \, 2, \, \ldots)$
are uniformly bounded. Namely, there exists three positive constants
$M_{B}, \, M_{g}, \, M_{A}$ such that
\begin{align}
    & \left\|\nabla^{2}f(x)\right\| \le M_{B},  \;
    \|A(x)\| \le M_{A}, \; \|A(x)^{\dagger}\| \le M_{A}, \forall x \in \mathrm{S},  \nonumber \\
    & \text{and} \; \|B_{k}\| \le M_{B}, \; \|\nabla f(x_{k})\| \le M_{g},
    \; k = 0, \, 1, \, 2, \, \ldots \label{BOUNDHESS}
\end{align}
hold.
\end{assumption}

\vskip 2mm

By combining the property $\|P_{k}\| = 1$ of the projection matrix $P_{k}$, from the
assumption  \eqref{BOUNDHESS}, we obtain
\begin{align}
    & \left\|P_{k}\nabla^{2}f(x_{k})P_{k}\right\|
    \le \|P_{k}\| \left\|\nabla^{2}f(x_{k})\right\| \|P_{k}\|
    = \left\|\nabla^{2}f(x_{k})\right\| \le M_{B}, \; \nonumber \\
    & \|P_{k}B_{k}P_{k}\| \le \|P_{k}\| \|B_{k}\| \|P_{k}\|
     = \|B_{k}\| \le M_{B}, \; k = 0, \, 1, \, 2, \, \ldots. \label{BOUNDPHESS}
\end{align}
According to the property of the matrix norm, we know that the absolute eigenvalue
of $P_{k}\nabla^{2}f(x_{k})P_{k}$ is less than $M_{B}$. If we denote
$\mu\left(P_{k}\nabla^{2}f(x_{k})P_{k}\right)$ as the eigenvalue of
$P_{k}\nabla^{2}f(x_k)P_{k}$, we know that the eigenvalue of
$\left(\frac{\sigma_{0}}{\Delta t_{k}}I + P_{k}\nabla^{2}f(x_k)P_{k}\right)$
is $\frac{\sigma_{0}}{\Delta t_{k}} + \mu\left(P_{k}\nabla^{2}(f_{k})P_{k}\right)$.
Consequently, from equation \eqref{BOUNDPHESS}, we known that
\begin{align}
   \frac{\sigma_{0}}{\Delta t_{k}}I + P_{k}\nabla^{2}f(x_{k})P_{k} \succ 0,
   \;   \text{when} \;    \Delta t_{k} < \frac{\sigma_{0}}{M_{B}}. \label{SPDCON}
\end{align}

\vskip 2mm

\textcolor{blue}{In order to prove the global convergence of Algorithm \ref{ALGRCM}, we need to
prove the iteration sequence $\{x_{k}\}$ preserves the constraint feasible
$c(x_{k}) \approx 0$ and its convergence point $x^{\ast}$ satisfies the KKT
condition \eqref{FOKKTG}. In the following Lemma \ref{LEMCONC}, we prove that
the iteration sequence $\{x_{k}\}$ generated by Algorithm \ref{ALGRCM} will
preserve on the neighbourhood of the feasibility $c(x) = 0$ when the time
steps $\Delta t_{k} \; (k = 0, \, 1, \, \ldots)$ are sufficiently small.
The proof of the convergence point $x^{\ast}$ satisfying the KKT condition
\eqref{FOKKTG} is left to Theorem \ref{THEGLCON}.}

\vskip 2mm

\begin{lemma} \label{LEMCONC}
Assume that $c(x)$ is continuously differentiable and its Jacobian matrix $A(x)$
is Lipschitz continuous. Namely, there exists a positive constant $L_{A}$ such that
\begin{align}
     \left\|A(x) - A(y) \right\| \le L_{A} \|x - y\|, \;
     \forall x, \, y  \in \Re^{n}. \label{LIPCONA}
\end{align}
Furthermore, we suppose $s_{k}$ is solved by equations \eqref{TRCM}-\eqref{GNEWTONM}
and $B_{k},  \; \nabla f(x_{k}) \, (k=0, \, 1, \, 2, \, \ldots)$ are
uniformly bounded. Namely, they satisfy Assumption \ref{ASSFUN}.
Then, there exists a positive constant $\delta_{c}$ such that
\begin{align}
    \|c(x_{k}+s_{k})\| \le \epsilon_{0},      \label{CXKLEEPS}
\end{align}
when $\Delta t_{k} \le \delta_{c}$, $0< \epsilon_{0} \le \min\{1/(L_{A}M_{A}^{2}),\; \epsilon\}$
and $\|c(x_{k})\| \le \epsilon_{0}$.
\end{lemma}

\vskip 2mm

\proof When $A_{k+1}^{p} = A(x_{k})$ and $\|c(x_{k+1})\| \le \epsilon_{0}$, we
obtain the result \eqref{CXKLEEPS}. Therefore, we only need to consider the case
$A_{k+1}^{p} = A(x_{k+1}^{p})$ below. From the first-order Taylor expansion and
equation \eqref{GNEWTONM}, we have
\begin{align}
    & c(x_{k}+s_{k}) = c(x_{k}+s_{k}^{p} + s_{k}^{c}) = c(x_{k+1}^{p} + s_{k}^{c})
    = c(x_{k+1}^{p}) + \int_{0}^{1} A(x_{k+1}^{p}+ ts_{k}^{c})s_{k}^{c}dt
    \nonumber \\
    & \hskip 2mm  =  c(x_{k+1}^{p}) + A(x_{k+1}^{p})s_{k}^{c}
    +\int_{0}^{1} (A(x_{k+1}^{p}+ ts_{k}^{c}) - A(x_{k+1}^{p}))s_{k}^{c}dt \nonumber \\
    & \hskip 2mm
    =  \int_{0}^{1} (A(x_{k+1}^{p}+ ts_{k}^{c}) - A(x_{k+1}^{p}))s_{k}^{c}dt.
    \label{TEXPCXK1}
\end{align}
By substituting the Lipschitz continuity \eqref{LIPCONA} of $A(x)$ and the bounded
assumption $\|A(x)^{\dagger}\| \le M_{A}$ into equation \eqref{TEXPCXK1}, we obtain
\begin{align}
    & \|c(x_{k}+s_{k})\| \le
    \int_{0}^{1} \|A(x_{k+1}^{p}+ ts_{k}^{c}) - A(x_{k+1}^{p})\| \|s_{k}^{c}\|dt
    \le \frac{1}{2} L_{A} \|s_{k}^{c}\|^{2}
    \nonumber \\
    & \hskip 2mm
    = \frac{1}{2} L_{A}  \left\|A(x_{k+1}^{p})^{\dagger}c(x_{k+1}^{p})\right\|^{2}
    \le \frac{1}{2} L_{A}  \left\|A(x_{k+1}^{p})^{\dagger}\right\|^{2}
    \left\|c(x_{k+1}^{p})\right\|^{2}
    \nonumber \\
    & \hskip 2mm
    \le \frac{1}{2} L_{A}  M_{A}^{2} \left\|c(x_{k+1}^{p})\right\|^{2}.
    \label{CXK1LEUB}
\end{align}

\vskip 2mm

In order to estimate the upper bound of $\|c(x_{k}+s_{k})\|$, we need to estimate
the upper bound of $\|c(x_{k+1}^{p})\|$. From the property $A_{k} P_{k} = 0$ of
the projection matrix $P_{k}$, we have
\begin{align}
    A_{k}s_{k}^{p} = \frac{\Delta t_{k}}{1+\Delta t_{k}}A_{k}P_{k}d_{k} = 0.
    \label{AKSKPZ}
\end{align}
Thus, from the first-order Taylor expansion, the property \eqref{AKSKPZ}
and the Lipschitz continuity \eqref{LIPCONA} of $A(x)$, we have
\begin{align}
    & \|c(x_{k+1}^{p})\| =
    \left\|c(x_{k}) + \int_{0}^{1}A(x_{k}+ts_{k}^{p})s_{k}^{p}dt\right\|
    \nonumber \\
    & \hskip 2mm
    = \left\|c(x_{k}) + \int_{0}^{1}(A(x_{k}+ts_{k}^{p})- A(x_{k}))s_{k}^{p}dt\right\|
    \nonumber \\
    & \hskip 2mm
    \le \|c(x_{k})\| + \int_{0}^{1} \|A(x_{k}+ts_{k}^{p}) - A(x_{k})\| \|s_{k}^{p}\|dt
    \le \|c(x_{k})\| + \frac{1}{2}L_{A} \|s_{k}^{p}\|^{2}. \label{CXK1PUB}
\end{align}
Furthermore, from equation \eqref{TRCM} and the assumption \eqref{BOUNDHESS} of $g(x)$ and
$B_{k}$, when $\Delta t_{k} \le \sigma_{0}/(2M_{B})$, we obtain
\begin{align}
    & \|s_{k}^{p}\| = \frac{\Delta t_{k}}{1 + \Delta t_{k}}
    \left\| \left(\frac{\sigma_{0}}{\Delta t_{k}} I + B_{k} \right)^{-1}P_{k}g_{k} \right\|
    \le \frac{\Delta t_{k}}{1 + \Delta t_{k}} \frac{\|P_{k}g_{k}\|}
    {\sigma_{0}/{\Delta t_{k}} - \|B_{k}\|} \nonumber \\
    & \hskip 2mm
    \le \frac{\Delta t_{k}}{1 + \Delta t_{k}} \frac{\|P_{k}\| \|g_{k}\|}
    {\sigma_{0}/{\Delta t_{k}} - M_{B}}
    \le \frac{\Delta t_{k}}{1 + \Delta t_{k}} \frac{M_{g}}
    {2M_{B} - M_{B}} = \frac{\Delta t_{k}}{1 + \Delta t_{k}}\frac{M_{g}}{M_{B}}.
    \label{SKPUPB}
\end{align}
According to the inductive hypothesis, we have $\|c(x_{k})\| \le \epsilon_{0}$. By
substituting it and equation \eqref{SKPUPB} into equation \eqref{CXK1PUB},
when $\Delta t_{k} \le \sigma_{0}/(2M_{B})$, we obtain
\begin{align}
    \|c(x_{k+1}^{p})\| \le \epsilon_{0} + \frac{1}{2}L_{A}
    \left(\frac{\Delta t_{k}}{1 + \Delta t_{k}}\frac{M_g}{M_B}\right)^{2}.
    \label{CXK1PUPBE}
\end{align}

\vskip 2mm

Thus, from equation \eqref{CXK1LEUB} and equation \eqref{CXK1PUPBE}, when
$\Delta t_{k} \le \sigma_{0}/(2M_{B})$, we prove
$\|c(x_{k}+s_{k})\|\le \epsilon_{0}$ if only we prove
\begin{align}
    \frac{1}{2}L_{A}M_{A}^{2}\left(\epsilon_{0} + \frac{1}{2}L_{A}
    \left(\frac{\Delta t_{k}}{1 + \Delta t_{k}}\frac{M_g}{M_B}\right)^{2}\right)^{2}
    \le \epsilon_{0}. \nonumber
\end{align}
Namely, when $\Delta t_{k} \le \sigma_{0}/(2M_{B})$, we only need to prove
\begin{align}
    \frac{\Delta t_{k}}{1+\Delta t_{k}} \le \frac{2M_{B}^{2}}{L_{A}M_{g}^{2}}
    \left(-\epsilon_{0} + \sqrt{\frac{2\epsilon_{0}}{L_{A}M_{A}^{2}}}\right).
    \label{DELTAKUPB}
\end{align}
We select
\begin{align}
    \delta_{c} = \min\left\{\frac{\sigma_{0}}{4M_{B}}, \;
    \frac{\frac{2M_{B}^{2}}{L_{A}M_{g}^{2}}\left(-\epsilon_{0}
    + \sqrt{\frac{2\epsilon_{0}}{L_{A}M_{A}^{2}}}\right)}{
    \left|1 - \frac{2M_{B}^{2}}{L_{A}M_{g}^{2}}\left(-\epsilon_{0}
    + \sqrt{\frac{2\epsilon_{0}}{L_{A}M_{A}^{2}}}\right)\right|}\right\}.
    \label{DELTACUPB}
\end{align}
Then, when $\Delta t_{k} \le \delta_{c}$ and $\epsilon_{0} \le \min \{1/(L_{A}M_{A}^{2}),
\; \epsilon\}$, we know that equation \eqref{DELTAKUPB} is true. Namely, we have
$\|c(x_{k}+s_{k})\| \le \epsilon_{0}$ when $\Delta t_{k} \le \delta_{c}$ and
$\epsilon_{0} \le \min \{1/(L_{A}M_{A}^{2}), \; \epsilon\}$.  \qed

\vskip 2mm

\textcolor{blue}{In order to estimate the lower bound of the time steps
$\Delta t_{k} \, (k = 1, \, 2, \, \ldots)$, we need to prove the correction
step $s_{k}^{c}$ is less than the predicted step $s_{k}^{p}$ and the predicted
vertical reduction $\textit{Vpred}_{k}$ is less than the predicted horizontal
reduction $\textit{Hpred}_{k}$, where $s_{k}^{c}$ is solved by equation
\eqref{GNEWTONM}, and $s_{k}^{p}$ is solved by equation \eqref{TRCM}, and
$\textit{Vpred}_{k}$, $\textit{Hpred}_{k}$ are defined by equation
\eqref{HVPREDK}. Their proofs are given by the following Lemma \ref{LEMVHPRD}
under some assumptions.}

\vskip 2mm

\begin{lemma} \label{LEMVHPRD}
Assume that $c(x)$ is continuously differentiable and its Jacobian matrix $A(x)$
satisfies the Lipschitz continuity \eqref{LIPCONA}. Furthermore, we suppose $s_{k}^{p}$
is solved by equation \eqref{TRCM} and $B_{k} \; (k=0, \, 1, \, 2, \, \ldots)$ are
uniformly bounded. Namely, $B_{k} \; (k=0, \, 1, \, 2, \, \ldots)$ satisfy equation
\eqref{BOUNDHESS}. Assume that $\|P_{k}g_{k}\| \ge \epsilon \; (k = 0, \, 1, \, 2, \, \ldots)$.
Then, there exists five positive constants $\theta_{c}, \; \theta_{v}, \; \delta_{0}, \;
\delta_{v}$ and $\delta_{\epsilon_0}$ such that
\begin{align}
    \|s_{k}^{c}\| \le \theta_{c} \|s_{k}^{p}\|^{2} \; \text{and} \;
    |\textit{Vpred}_{k}| \le \theta_{v}
    \|P_{k}g_{k}\| \|s_{k}^{p}\|,  \label{VLHPREDK}
\end{align}
when $\delta_{0} \le \Delta t_{k} \le \delta_{v}$, \; $0< \epsilon_{0} \le \delta_{\epsilon_0}$
and $\|c(x_{k})\| \le \epsilon_{0}$, where $0 < \theta_{v} < 1/2$,
$s_{k}^{c}$ is solved by equation \eqref{GNEWTONM},
and $\textit{Vpred}_{k}$ is defined by equation \eqref{HVPREDK}.
\end{lemma}

\vskip 2mm

\proof From equation \eqref{TRCM}, when $\Delta t_{k} < \sigma_{0}/\|B_{k}\|$,
we have
\begin{align}
    P_{k}d_{k} = - P_{k}\left(\frac{\sigma_{0}}{\Delta t_{k}}I+B_{k}\right)^{-1}P_{k}g_{k}.
    \nonumber
\end{align}
By combining it with the Cauchy-Schwartz inequality $|x^{T}y| \le \|x\| \, \|y\|$ and
$P_{k}^{2} = P_{k}$, when $\Delta t_{k} < \sigma_{0}/\|B_{k}\|$,
we have
\begin{align}
    & \|P_{k}g_{k}\| \|P_{k}d_{k}\| \ge
    |g_{k}^{T} P_{k}^{2} d_{k}| =
    \left|(P_{k}g_{k})^{T}\left(\frac{\sigma_{0}}{\Delta t_{k}}I+B_{k}\right)^{-1}(P_{k}g_{k})\right|
    \nonumber \\
    & \hskip 2mm
    \ge \frac{1}{\sigma_{0}/{\Delta t_{k}} + \|B_k\|} \|P_{k}g_{k}\|^{2}. \label{PKDKGELB}
\end{align}
Namely, when $\Delta t_{k} < \sigma_{0}/\|B_{k}\|$, from equation \eqref{PKDKGELB},
we have
\begin{align}
   \|s_{k}^{p}\| = \left\|\frac{\Delta t_{k}}{1+\Delta t_{k}}P_{k}d_{k}\right\|
   = \frac{\Delta t_{k}}{1+\Delta t_{k}} \|P_{k}d_{k}\|
   \ge \frac{\Delta t_{k}}{1+\Delta t_{k}}
   \frac{1}{\sigma_{0}/{\Delta t_{k}} + \|B_k\|} \|P_{k}g_{k}\|.
   \label{SKPGELB}
\end{align}
Thus, according to the assumptions $\|B_k\| \le M_{B}$ and
$\|P_{k}g_{k}\| \ge \epsilon$, when $\delta_{0} \le \Delta t_{k}
<  \sigma_{0}/M_{B}$, from equation \eqref{SKPGELB}, we have
\begin{align}
   \|s_{k}^{p}\| \ge \frac{\delta_{0}}{1+\delta_{0}}
   \frac{1}{\sigma_{0}/\delta_{0} + M_{B}} \epsilon
   = \frac{\delta_{0}\epsilon}{(1+\delta_{0})(\sigma_{0}/\delta_{0} + M_{B})}.
   \label{SKPGELBA}
\end{align}

\vskip 2mm

From equation \eqref{GNEWTONM}, the Lipschitz continuity \eqref{LIPCONA} and
equation \eqref{CXK1PUB}, we have
\begin{align}
    \|s_{k}^{c}\| & = \|A(x_{k+1}^{p})^{\dagger}c(x_{k+1}^{p})\|
    \le M_{A} \|c(x_{k+1}^{p})\| = M_{A}\|c(x_{k}+s_{k}^{p})\|
    \nonumber \\
    &
    \le  M_{A} \left(\|c(x_{k})\| + \frac{1}{2}L_{A} \|s_{k}^{p}\|^{2}\right)
    \le  M_{A} \left(\epsilon_{0} + \frac{1}{2}L_{A} \|s_{k}^{p}\|^{2}\right).
    \label{SKCLEUB}
\end{align}
We select
\begin{align}
   \delta_{\epsilon_0} = \min \left\{\left(\frac{\delta_{0}\epsilon}{(1+\delta_{0})
   (\sigma_{0}/\delta_{0} + M_{B})}\right)^{2}, \;
   \frac{1}{L_{A}M_{A}^{2}} \right\} .
   \label{EPS0LEUB}
\end{align}
Then, from equations  \eqref{SKPGELBA}-\eqref{EPS0LEUB}, when
$\epsilon_{0} \le \delta_{\epsilon_0}$, we have
\begin{align}
   \|s_{k}^{c}\| \le M_{A} \left(1+\frac{1}{2}L_{A}\right)\|s_{k}^{p}\|^{2}
   = \theta_{c} \|s_{k}^{p}\|^{2},    \label{SKCLESKP2}
\end{align}
where we select $\theta_{c} = M_{A}(1+0.5L_{A})$. Thus, we prove the first
inequality of equation \eqref{VLHPREDK}.

\vskip 2mm

From equation \eqref{TRCM}, when $\Delta t_{k} \le \sigma_{0}/(2M_B)$,
we have
\begin{align}
    & \|s_{k}^{p}\| = \left\| \frac{\Delta t_{k}}
    {1 + \Delta t_{k}} P_{k}d_{k} \right\| \le \|P_{k}d_{k}\|
    \le \|d_{k}\|
    = \left\|\left(\frac{\sigma_{0}}{\Delta t_{k}} I + B_{k}\right)^{-1}P_{k}g_{k}\right\|
    \nonumber \\
    & \hskip 2mm
    \le \frac{1}{\sigma_{0}/{\Delta t_k} - M_{B}} \|P_{k}g_{k}\|
    \le \frac{M_g}{M_B}.
    \label{SKPLEUB}
\end{align}
Thus, from equation \eqref{HVPREDK}, equation \eqref{SKCLESKP2} and
the assumption \eqref{BOUNDHESS}, we have
\begin{align}
     & |\textit{Vpred}_{k}|  = |q_{k}(s_{k}^{p}) - q_{k}(s_{k}^{p} + s_{k}^{c})|
     = \left|(s_{k}^{c})^{T}\nabla q_{k}(s_{k}^{p}) +
      \frac{1}{2}(s_{k}^{c})^{T}\nabla^{2}q_{k}(s_{k}^{p})s_{k}^{c}\right|
      \nonumber \\
      & \hskip 2mm
      =  \left|(s_{k}^{c})^{T}\left(g_{k} + B_{k}s_{k}^{p}
      + \frac{1}{2}B_{k}s_{k}^{c}\right)\right|
      \le \|s_{k}^{c}\| \left\| g_{k} + B_{k}s_{k}^{p}
      + \frac{1}{2}B_{k}s_{k}^{c} \right\| \nonumber \\
      & \hskip 2mm
       \le \|s_{k}^{c}\| (\|g_k\| + \|B_{k}\| (\|s_{k}^{p}\|
      + \|s_{k}^{c}\|))
      \le \theta_{c}
      (M_{g} + M_{B}(\|s_{k}^{p}\|+\|s_{k}^{c}\|))\|s_{k}^{p}\|^{2}.
      \label{VPREDKLEUB}
\end{align}
By substituting equation \eqref{SKCLESKP2} and equation \eqref{SKPLEUB} into
equation \eqref{VPREDKLEUB}, when $\Delta t_{k} \le \sigma_{0}/(2M_B)$,
we obtain
\begin{align}
    |\textit{Vpred}_{k}| \le \frac{\theta_{c}(2M_{g}M_{B} + \theta_{c}M_{g}^{2})}
    {M_B} \|s_{k}^{p}\|^{2}
    \le \frac{\theta_{c}(2M_{g}M_{B} + \theta_{c}M_{g}^{2})}
    {M_B (\sigma_{0}/{\Delta t_k} - M_{B})} \|P_{k}g_{k}\| \|s_{k}^{p}\|.
    \label{VPREDKLPGSKP}
\end{align}
We select
\begin{align}
    \delta_{v} = \min \left\{\frac{\sigma_{0}}{4M_B}, \; \delta_{c}, \;
    \frac{\sigma_{0}M_{B}\theta_{v}}
    {(2M_{g}M_{B} + \theta_{c}M_{g}^{2})\theta_{c} + M_{B}^{2}\theta_{v}}\right\},
    \label{DELTAV}
\end{align}
where the positive constant $\delta_{c}$ is defined by equation \eqref{DELTACUPB}.
Thus, from equations \eqref{VPREDKLPGSKP}-\eqref{DELTAV}, when
$\Delta t_{k} \le \delta_{v}$, we have
$|\textit{Vpred}_{k}| \le \theta_{v} \|P_{k}g_{k}\| \|s_{k}^{p}\|$.  We select
the small positive constant $\delta_{0}$ to satisfy $\delta_{0}
\le \delta_{v}$. Then, by combining equation \eqref{SKCLESKP2}, when
$\delta_{0} \le \Delta t_{k} \le \delta_{v}$ and $\epsilon_{0} \le \delta_{\epsilon_0}$,
we obtain the result \eqref{VLHPREDK}.  \qed

\vskip 2mm

\begin{remark}
We select
\begin{align}
    \theta_{1} \ge \frac{\theta_{c}M_{g}}{M_{B}}. \label{THETA1}
\end{align}
Then, from equations \eqref{SKCLESKP2}-\eqref{THETA1}, when $\delta_{0} \le \Delta t_{k}
< \sigma_{0}/(2M_{B})$ and $\epsilon_{0} \le \delta_{\epsilon_0}$, we have
\begin{align}
    \|s_{k}^{c}\| \le \theta_{1} \|s_{k}^{p}\|. \label{SKCLESKP}
\end{align}
Namely, the condition \eqref{VPRLEHPR} is satisfied. \qed
\end{remark}

\vskip 2mm

\textcolor{blue}{In order to prove that $p_{g_k}$ converges to zero when $k$ tends
to infinity, we need to estimate the lower bound of time step sizes
$\Delta t_{k} \, (k = 1, \, 2, \, \ldots)$. By using Lemma \ref{LEMHPRED},
Lemma \ref{LEMCONC} and Lemma \ref{LEMVHPRD}, we can obtain the lower bound of time step sizes
$\Delta t_{k} \, (k = 1, \, 2, \, \ldots)$ as follows.}

\vskip 2mm

\begin{lemma} \label{LEMDTLB}
Assume that $f: \; \Re^{n} \rightarrow \Re$ is twice continuously differentiable
and $\nabla f(x)$ is Lipschitz continuous. Namely, there exists a positive
constant $L_{g}$ such that
\begin{align}
     \left\|\nabla  f(x) - \nabla  f(y) \right\| \le L_{g} \|x - y\|, \;
     \forall x, \, y  \in \Re^{n}. \label{LIPSCHCON}
\end{align}
Assume that $c(x)$ is continuously differentiable and its Jacobian matrix $A(x)$
satisfies the Lipschitz continuity \eqref{LIPCONA}. We suppose that the sequence
$\{x_{k}\}$ is generated by Algorithm \ref{ALGRCM} and Assumption \ref{ASSFUN}
holds. Furthermore, we assume
\begin{align}
    \|P_{k}g_{k}\| > \epsilon  \label{GKGEEPS}
\end{align}
holds for all $k = 0, \, 1, \, 2, \, \ldots$, where $g_{k} = g(x_k) = \nabla f(x_k), \;
P_{k} = P(x_{k})$ and the projection matrix $P(x)$ is defined by equation
\eqref{PROMAT}. Then, it exists a positive constant $\delta_{\Delta t}$ such that
\begin{align}
    \Delta t_{k} \ge \gamma_{2} \delta_{\Delta t} > 0, \; k = 0, \,  1, \, 2, \dots,
    \label{DTGEPN}
\end{align}
where $\Delta t_{k}$ is adaptively adjusted by the trust-region updating strategy
\eqref{MRHOK}-\eqref{ADTK1}.
\end{lemma}

\vskip 2mm

\proof According to Lemma \ref{LEMCONC}, by induction,
when $\delta_{0} \le \Delta t_{k} \le \delta_{v}$ and
$0< \epsilon_{0} \le \delta_{\epsilon_0}$, we have
\begin{align}
     \|c(x_{k})\| \le \epsilon_{0}, \; k = 0, \, 1, \, 2, \, \ldots,
     \label{CXKLEEPS0}
\end{align}
where $\delta_{\epsilon_0}$ is defined by equation \eqref{EPS0LEUB} and
$\delta_{v}$ is defined by equation \eqref{DELTAV}.

\vskip 2mm

From equation \eqref{MATNP}, we know $\|P_{k}\| = 1$. By using this property and the assumption
\eqref{BOUNDHESS}, we have
\begin{align}
     &\left|\mu_{\min}\left(B_{k}-P_{k}^{T}B_{k}P_{k}\right)\right| \le
     \left\|B_{k}-P_{k}^{T}B_{k}P_{k}\right\| \nonumber \\
     & \quad \le \|B_{k}\| + \|P_{k}^{T}\|\|B_{k}\|\|P_{k}\|
     = 2\|B_{k}\| \le 2M_{B}, \; k=0, \, 1, \, 2,\ldots,  \label{PROGKUPBD}
\end{align}
where $\mu_{\min}(B)$ represents the smallest eigenvalue of matrix $B$. Thus,
from equation \eqref{PROGKUPBD}, we obtain
\begin{align}
    & \mu_{\min}\left(\frac{\sigma_{0}}{\Delta t_{k}}I + B_{k}-P_{k}^{T}B_{k}P_{k}\right)
      = \frac{\sigma_{0}}{\Delta t_{k}} + \mu_{\min}\left(B_{k}-P_{k}^{T}B_{k}P_{k}\right)
     \nonumber \\
    & \quad \ge \frac{\sigma_{0}}{\Delta t_{k}} - 2M_{B}, \; k=0, \, 1, \, 2, \ldots.
    \label{SEPROGLBD}
\end{align}
Similarly, from the assumption \eqref{BOUNDHESS}, we have
\begin{align}
    \mu_{\min} \left(\frac{\sigma_{0}}{\Delta t_{k}} I + B_{k}\right)
     =  \frac{\sigma_{0}}{\Delta t_{k}}+ \mu_{\min} \left(B_{k}\right)
    \ge \frac{\sigma_{0}}{\Delta t_{k}} - M_{B}, \;  k = 0, \, 1, \, 2, \ldots.
    \label{SEGLBD}
\end{align}
Therefore, the positive definite conditions of equation \eqref{PSDASS} are
satisfied when $\Delta t_{k} \le \sigma_{0}/(2M_{B}) \, (k=1,\,2,\ldots)$.

\vskip 2mm

From the Lipschitz continuity \eqref{LIPSCHCON} of $\nabla f(x)$ and the
Cauchy-Schwartz inequality $|x^{T}y| \le \|x\| \|y\|$, we have
\begin{align}
    & \left|\int_{0}^{1} (\nabla f(x_{k}+t s_{k}) - \nabla f(x_{k}))^{T}s_{k}dt \right|
     \le \int_{0}^{1} \|\nabla f(x_{k}+t s_{k}) - \nabla f(x_{k})\| \|s_{k}\|dt \nonumber \\
    & \hskip 2mm \le \int_{0}^{1} L_{g} \|s_{k}\|^{2} t dt
    = \frac{1}{2}L_{g} \|s_{k}\|^{2}
    =  \frac{1}{2} L_{g} (\|s_{k}^{p} + s_{k}^{c}\|)^{2}
    \le \frac{1}{2} L_{g} (\|s_{k}^{p}\| + \|s_{k}^{c}\|)^{2}. \label{BMSOFK}
\end{align}
By substituting equation \eqref{VLHPREDK} and equation \eqref{SKPLEUB} into
equation \eqref{BMSOFK}, when $\delta_{0} \le \Delta t_{k} \le \delta_{v}$ and
$0< \epsilon_{0} \le \delta_{\epsilon_0}$, we obtain
\begin{align}
     & \left|\int_{0}^{1} (\nabla f(x_{k}+t s_{k}) - \nabla f(x_{k}))^{T}s_{k}dt \right|
     \le \frac{1}{2} L_{g}  (\|s_{k}^{p}\| + \theta_{c} \|s_{k}^{p}\|^{2})^{2}
     \nonumber \\
     & \hskip 2mm
     \le \frac{1}{2} L_{g}\left(1 + \frac{\theta_{c} M_g}{M_B}\right)^{2}\|s_{k}^{p}\|^{2},
     \label{BMSOFKP}
\end{align}
and
\begin{align}
    & \frac{1}{2} |s_{k}^{T}B_{k}s_{k}| \le \frac{1}{2} \|B_{k}\| \|s_{k}\|^{2}
    \le \frac{1}{2} M_{B} \|s_{k}^{p} + s_{k}^{c}\|^{2}
    \le \frac{1}{2} M_{B} (\|s_{k}^{p}\| + \|s_{k}^{c}\|)^{2} \nonumber \\
    & \le \frac{1}{2}M_{B}\left(\|s_{k}^{p}\| + \theta_{c} \|s_{k}^{p}\|^{2}\right)^{2}
    = \frac{1}{2}M_{B}(1+ \theta_{c}\|s_{k}^{p}\|)^{2} \|s_{k}^{p}\|^{2} \nonumber \\
    & \le \frac{1}{2}M_{B}\left(1 + \theta_{c} \frac{M_g}{M_B}\right)^{2}
    \|s_{k}^{p}\|^{2}. \label{SKBKSKUB}
\end{align}

\vskip 2mm

From the first-order Taylor expansion of $f(x_{k}+s_{k})$ and equation \eqref{BMSOFKP}, we have
\begin{align}
    & \left|f(x_{k})  - f(x_{k}+s_{k}) + \nabla f(x_{k})^{T}s_{k}\right|
    =  \left|\nabla f(x_{k})^{T}s_{k} - \int_{0}^{1} \nabla f(x_{k}
    + t s_{k})^{T} s_{k}dt\right| \nonumber \\
    & = \left|\int_{0}^{1}(\nabla f(x_{k}) - \nabla f(x_{k}+ts_{k}))^{T}s_{k}dt \right|
    \le \frac{1}{2} L_{g}\left(1 + \frac{\theta_{c} M_g}{M_B}\right)^{2}\|s_{k}^{p}\|^{2}.
     \label{SOTEFK}
\end{align}
From equation \eqref{HPREDKEST} and equation \eqref{VLHPREDK}, we have
\begin{align}
    & q_{k}(0) - q_{k}(s_k) = (q_{k}(0) - q_{k}(s_{k}^{p}))
    - (q_{k}(s_{k}^{p}) - q_{k}(s_{k}^{p} + s_{k}^{c})) \nonumber \\
    & \ge \frac{1}{2}\|P_{k}g_{k}\| \|s_{k}^{p}\| - \theta_{v} \|P_{k}g_{k}\| \|s_{k}^{p}\|
    = (0.5 - \theta_{v}) \|P_{k}g_{k}\| \|s_{k}^{p}\|. \label{PREESTLB}
\end{align}
Thus, from equations \eqref{QOAM}, \eqref{MRHOK}, \eqref{HPREDKEST}, \eqref{VLHPREDK}
and \eqref{SKBKSKUB}-\eqref{PREESTLB}, when $\delta_{0} \le \Delta t_{k} \le \delta_{v}$
and $0< \epsilon_{0} \le \delta_{\epsilon_0}$, we have
\begin{align}
     & |\rho_{k} - 1| =  \left|\frac{(f(x_{k}) - f(x_{k}+s_{k}))
     - (q_{k}(0) - q_{k}(s_{k}))}{q_{k}(0) - q_{k}(s_{k})}\right|
     \nonumber \\
     & \le \frac{0.5(L_{g}+M_{B})(1+ (\theta_{c}M_{g})/M_{B})^{2}\|s_{k}^{p}\|^{2}}
     {|q_{k}(0) - q_{k}(s_{k})|} \nonumber \\
     & \le \frac{0.5(L_{g}+M_{B})(1+ (\theta_{c}M_{g})/M_{B})^{2}\|s_{k}^{p}\|^{2}}
     {\left| |q_{k}(0) - q_{k}(s_{k}^{p})| - |q_{k}(s_{k}^{p}) - q_{k}(s_{k}^{p}+ s_{k}^{c})| \right|}
     \nonumber \\
    & \le \frac{0.5(L_{g}+M_{B})(1+ (\theta_{c}M_{g})/M_{B})^{2}\|s_{k}^{p}\|^{2}}
     {(0.5-\theta_{v})\|P_{k}g_{k}\| \|s_{k}^{p}\|}.
    \label{ESTRHOK}
\end{align}
By substituting equation \eqref{SKPLEUB} into equation \eqref{ESTRHOK}, we obtain
\begin{align}
    |\rho_{k} - 1| \le \frac{0.5(L_{g}+M_{B})(1+ (\theta_{c}M_{g})/M_{B})^{2}}
     {(0.5-\theta_{v})(\sigma_{0}/{\Delta t_k} - M_{B})}.
     \label{RHOKLEUB}
\end{align}
We denote
\begin{align}
    & \delta_{ub} = \frac{(0.5-\theta_{v})(1-\eta_{2})\sigma_{0}}
     {0.5(L_{g}+M_{B})\left(1+ (\theta_{c}M_{g})/M_B\right)^{2}
     +(0.5-\theta_{v})(1 - \eta_{2})M_{B}}, \nonumber \\
    & \delta_{\Delta t} \triangleq
     \min\left\{ \delta_{ub}, \; \frac{\sigma_{0}}{4M_{B}}, \; \delta_{v},
     \; \Delta t_{0} \right\}. \label{DELTAKUPBD}
\end{align}
Thus, from equations \eqref{RHOKLEUB}-\eqref{DELTAKUPBD}, when $\delta_{0} \le
\Delta t_{k} \le \delta_{\Delta t}$, we have $\rho_{k} \ge \eta_{2}$. We
select $\delta_{0} \le \gamma_{2} \delta_{\Delta t}$. We assume that $K$
is the first index such that $\Delta t_{K} \le \delta_{\Delta t}$. Then, from
equations \eqref{RHOKLEUB}-\eqref{DELTAKUPBD}, we know that
$\rho_{K} \ge \eta_{2}$. According to the time-stepping adjustment
formula \eqref{ADTK1} and Lemma \ref{LEMCONC}, when $0< \epsilon_{0} \le
\delta_{\epsilon_0}$, $x_{K} + s_{K}$ will be accepted and the time step
$\Delta t_{K+1}$ will be enlarged. Consequently,
$\Delta t_{k}\ge \gamma_{2}\delta_{\Delta t}$ holds for all
$k = 0, \, 1, \, 2, \ldots$. \qed

\vskip 2mm

\begin{remark}
We select $\eta_{q} \le 0.5 - \theta_{v}$. Then, when $\delta_{0} \le
\Delta t_{k} \le \delta_{v}$, from equation \eqref{PREESTLB}, we have
\begin{align}
    q_{k}(0) - q_{k}(s_k) \ge (0.5 - \theta_{v}) \|P_{k}g_{k}\| \|s_{k}^{p}\|
    \ge \eta_{q} \|P_{k}g_{k}\| \|s_{k}^{p}\|. \label{PREDKLB}
\end{align}
Namely, the sufficient descent condition \eqref{AMGESKPGK} is satisfied. \qed
\end{remark}

By using the results of Lemma \ref{LEMCONC} and Lemma \ref{LEMDTLB}, we prove
the global convergence of Algorithm \ref{ALGRCM} for the equality-constrained
optimization problem \eqref{NLEQOPT} as follows.

\vskip 2mm

\begin{theorem} \label{THEGLCON}
Assume that Assumption \ref{ASSFUN} holds, $\nabla f(x)$ satisfies the Lipschitz
continuity \eqref{LIPSCHCON}, and the Jacobian matrix $A(x)$ of $c(x)$ satisfies
the Lipschitz continuity \eqref{LIPCONA}. We suppose that the sequence
$\{x_{k}\}$ is generated by Algorithm \ref{ALGRCM} and $f(x_{k}) \, (k = 0, \, 1, \, \ldots)$
are bounded below. Then, for any positive constant $\epsilon$, we can select a sufficiently
small constant $\epsilon_{0} > 0$ such that
\begin{align}
    \|P_{K}g_{K}\| \le \epsilon \; \text{and} \; \|c(x_{K})\| \le \epsilon_{0}
    \label{PKGKLEEPS}
\end{align}
hold for an index $K$ when $\|c(x_{0})\| \le \epsilon_{0}$.
\end{theorem}
\proof We prove this result by contradiction. Assume that the result
\eqref{PKGKLEEPS} is not true. Then, we have
\begin{align}
  \|P_{k}g_{k}\| > \epsilon  \; \text{for all} \;  k = 0, \, 1, \, 2, \ldots.
  \label{PGKGECON}
\end{align}
Thus, from Lemma \ref{LEMCONC} and Lemma \ref{LEMDTLB}, we can select a
sufficiently small constant $\epsilon_{0} > 0$ and there exists an
infinite subsequence $\{x_{k_i}\}$ such that the trial steps
$s_{k_i}$ are accepted, i.e., $\rho_{k_i} \ge \eta_{a}, \, i = 1, \, 2, \ldots$.
Otherwise, all steps are rejected after a given iteration index, then the time
step will keep decreasing, which contradicts \eqref{DTGEPN}. Therefore,
from equation \eqref{MRHOK}, we have $\|c(x_{k_i})\| \le \epsilon_{0}$ and
\begin{align}
    f_{0} - \lim_{k \to \infty} f_{k} = \sum_{k = 0}^{\infty} (f_{k} - f_{k+1})
    \ge \eta_{a} \sum_{i = 0}^{\infty}
    \left(q_{k_{i}}(0) - q_{k_{i}}(s_{k_i})\right),
    \label{LIMSUMFK}
\end{align}
where $s_{k_{i}}$ is computed by equation \eqref{TRCM} and equation \eqref{GNEWTONM}.

\vskip 2mm

Since $f_{k} = f(x_{k}) \, (k = 0, \, 1, \, 2, \, \ldots)$ are bounded below and
the sequence $\{f_{k}\}$ decreases monotonically, we know that the limit
of the sequence $\{f(x_{k})\}$ exits and we denote it as $\lim_{k \to \infty} f_{k}
= f^{\ast}$. Thus, from equation \eqref{LIMSUMFK}, we have
\begin{align}
    \lim_{k_{i} \to \infty}
    \left(q_{k_{i}}(0) - q_{k_{i}}(s_{k_i})\right) = 0.
    \label{LIMQK}
\end{align}
By substituting equation \eqref{AMGESKPGK} into equation \eqref{LIMQK}, we
obtain
\begin{align}
     \lim_{i \to \infty} \|p_{g_{k_i}}\| \|s_{k_i}^{p}\| = 0.  \label{LIMPK}
\end{align}
According to the assumption $\|P_{k}g_{k}\| > \epsilon > 0$, from equation
\eqref{LIMPK}, we have
\begin{align}
    \lim_{i \to \infty} \|s_{k_i}^{p}\| = 0.   \label{LIMPDK}
\end{align}

\vskip 2mm

From Lemma \ref{LEMDTLB}, we know that it exists a positive constant
$\delta_{\Delta t}$ such that
\begin{align}
     \Delta t_{k} \ge \gamma_{2} \delta_{\Delta t} > 0, \;  k = 0, \, 1,\, 2, \ldots.
 \label{DELTKGE}
\end{align}
(i) When $B_{k_i}$ is updated by the BFGS quasi-Newton formula \eqref{BFGS}, we know
that $B_{k_i} \succ 0$. Thus, by substituting equation \eqref{DELTKGE} into
equation \eqref{SKPGELB}, we have
\begin{align}
     \|s_{k_i}^{p}\| & \ge \frac{\Delta t_{k_i}}{1+\Delta t_{k_i}}
   \frac{1}{\sigma_{0}/{\Delta t_{k_i}} + \|B_{k_i}\|} \|P_{k_i}g_{k_i}\|
   \nonumber \\
   & \ge \frac{\gamma_{2} \delta_{\Delta t}}{1+ \gamma_{2} \delta_{\Delta t}}
   \frac{1}{\sigma_{0}/{(\gamma_{2}\delta_{\Delta t})} + M_B} \|P_{k_i}g_{k_i}\|.
     \label{PDKEQ}
\end{align}
(ii) When $B_{k_i} = P_{k_i}\nabla^{2}f(x_{k_i})P_{k_i}$, from equation
\eqref{TRCM}, we have
\begin{align}
    & P_{k_i} \left(\frac{\sigma_{0}}{\Delta t_{k_i}}I
    + P_{k_i}\nabla^{2}f(x_{k_i})P_{k_i}\right)d_{k_i}
    = - P_{k_i}^{2}g_{k_i}
    \nonumber \\
    & \hskip 2mm
    = - P_{k_i}g_{k_i}
  = \left(\frac{\sigma_{0}}{\Delta t_{k_i}}I
    + P_{k_i}\nabla^{2}f(x_{k_i})P_{k_i}\right)d_{k_i},
    \nonumber
\end{align}
which gives $P_{k_i}d_{k_i} = d_{k_i}$. By combining it with
equation \eqref{TRCM} and equation \eqref{DELTKGE}, we have
\begin{align}
    & \frac{\gamma_{2} \delta_{\Delta t}}{1+ \gamma_{2} \delta_{\Delta t}}
    \|P_{k_i}g_{k_i}\|
    \le \left\|\frac{\Delta t_{k_i}}{1+\Delta t_{k_i}} P_{k_i}g_{k_i}\right\|
    \nonumber \\
    & \hskip 2mm
    = \left\|\frac{\Delta t_{k_i}}{1+\Delta t_{k_i}} \left(\frac{\sigma_{0}}
    {\Delta t_{k_i}}I  + P_{k_i} \nabla^{2}f(x_{k_i}) P_{k_i}\right) d_{k_i}\right\|
    \nonumber \\
    & \hskip 2mm
    = \left\|\frac{\Delta t_{k_i}}{1+\Delta t_{k_i}}
    \left(\frac{\sigma_{0}}{\Delta t_{k_i}}I  + P_{k_i} \nabla^{2}f(x_{k_i})
    P_{k_i}\right) P_{k_i}d_{k_i}\right\|
    \nonumber \\
    & \hskip 2mm
    = \left\|\left(\frac{\sigma_{0}}{\Delta t_{k_i}}I  + P_{k_i} \nabla^{2}f(x_{k_i})
    P_{k_i}\right) s_{k_i}^{p}\right\|
    \nonumber \\
    & \hskip 2mm
    \le \left(\frac{\sigma_{0}}{\Delta t_{k_i}} + \|B_{k_i}\|\right) \|s_{k_i}^{p}\|
    \le \left(\frac{\sigma_{0}}{\gamma_{2}\delta_{\Delta t}} + M_{B} \right) \|s_{k_i}^{p}\|.
    \label{PDKEQPH}
\end{align}
By substituting equations \eqref{PDKEQ}-\eqref{PDKEQPH} into equation \eqref{LIMPDK},
we obtain
\begin{align}
     \lim_{{i} \to \infty} \|P_{k_i}g_{k_i}\| = 0,  \label{PGKTOZ}
\end{align}
which contradicts the bounded assumption \eqref{PGKGECON} of
$P_{k}g_{k}  \, (k = 0, \, 1, \, 2, \, \ldots)$. Therefore, the result
\eqref{PKGKLEEPS} is true.  \qed

\vskip 2mm

\section{Numerical Experiments}

In this section, we conduct some numerical experiments to test the performance
of Algorithm \ref{ALGRCM} (the Rcm method). The codes are executed by a HP
notebook with the Intel quad-core CPU and 8Gb memory in the MATLAB R2020a
environment \cite{MATLAB}. For a real-world problem, the analytical Hessian matrix
$\nabla^{2}f(x_{k})$ may not be offered. Thus, in practice, we replace the
two sided projection $P(x_{k})\nabla^{2}f(x_{k})P(x_{k})$ of the Hessian matrix
with its difference approximation as follows:
\begin{align}
     P_{k}\nabla^{2}f(x_{k})P_{k} \approx
    \left[\frac{P_{k}g(x_{k} + \epsilon P_{k}e_{1}) - P_{k}g(x_{k})}{\epsilon}, \,
    \ldots, \, \frac{P_{k}g(x_{k} + \epsilon P_{k}e_{n}) - P_{k}g(x_{k})}{\epsilon}\right],
    \label{NUMHESS}
\end{align}
where the elements of $e_{i}$ equal zeros except for the $i$-th element which equals 1,
$P_{k} = P(x_{k})$ and the parameter $\epsilon$ can be selected as $10^{-6}$ according
to our numerical experiments.

\vskip 2mm

SQP \cite{FP1963,Goldfarb1970,NW1999,Wilson1963} is a representative
method for constrained optimization problems. And there are two representative
implementation codes of the SQP method. One is the built-in subroutine fmincon.m
(NLPQL) of the MATLAB2020a environment \cite{MATLAB,Schittkowski1986}. The other
is the subroutine SNOPT \cite{GMS2005,GMS2006} executed in GAMS v28.2 (2019)
environment \cite{GAMS}). fmincon and SNOPT are two robust and efficient solvers
for nonlinear equality-constrained optimization problems. Therefore, we select
these two typical solvers as the basis for comparison.

\vskip 2mm

We select $65$ nonlinear equality-constrained optimization problems from CUTEst
\cite{GOT2015} and construct $60$ test problems. Therefore, we use these $125$
problems to test their performance of these three algorithms (Rcm, fmincon and SNOPT).
For those $60$ constructed test problems, we use Ackley Function \cite{SB2013}
as their objective functions. Ackley Function can be written as:
\begin{align}
    f(x) = -a e^{-b \sqrt{\frac{1}{n} \sum_{i=1}^{n} x_{i}^{2}}}
    - e^{\frac{1}{n} \sum_{i=1}^{n} \cos(2\pi x_{i})} + a + e,    \label{ACKLEYFUN}
\end{align}
where $a=20, \; b=0.2, \; n=2000$. Then, we use the gradients
of $20$ unconstrained optimization problems from \cite{Andrei2008,Luksan1994,MGH1981,SB2013}
as their constrained functions. In order to test the performance of
the method for the different dimensions of constraints, we test the problems
with few constraints ($m << n$ such as $m = 10, \; n = 2000)$, medium constraints
($m \approx \frac{1}{2}n$ such as $m = 1000, \; n = 2000$) and many constraints
($m \approx n$ such as $m = 1999, \; n = 2000$).

\vskip 2mm

\begin{figure}
\centering
\begin{minipage}{.9\textwidth}
  \centering
  \includegraphics[width=.9\linewidth]{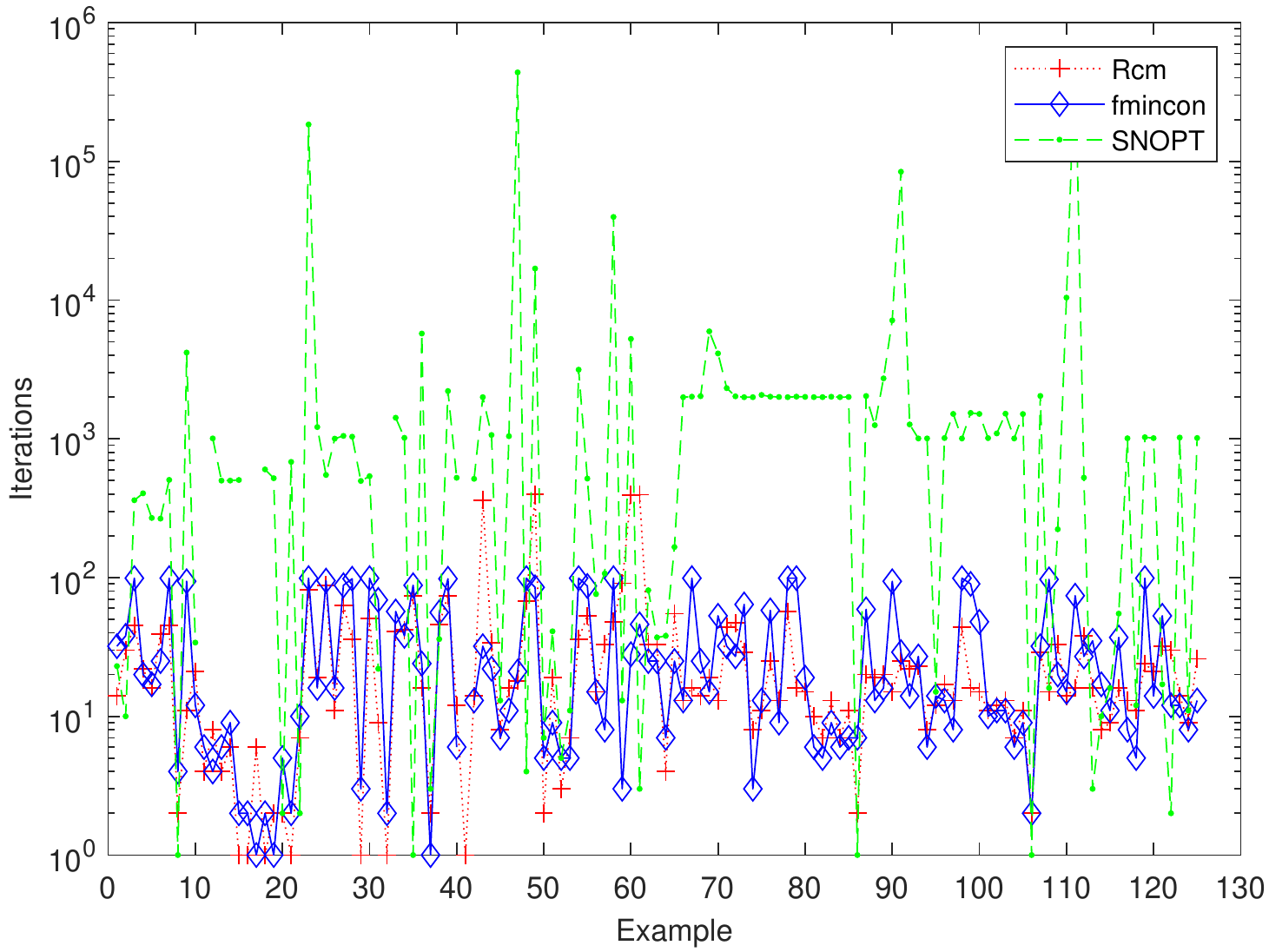}
  \captionof{figure}{Iterations of Rcm, fmincon and SNOPT for test problems.}
  \label{FIGITERCM}
\end{minipage}%

\vskip 2mm

\begin{minipage}{.9\textwidth}
  \centering
  \includegraphics[width=.9\linewidth]{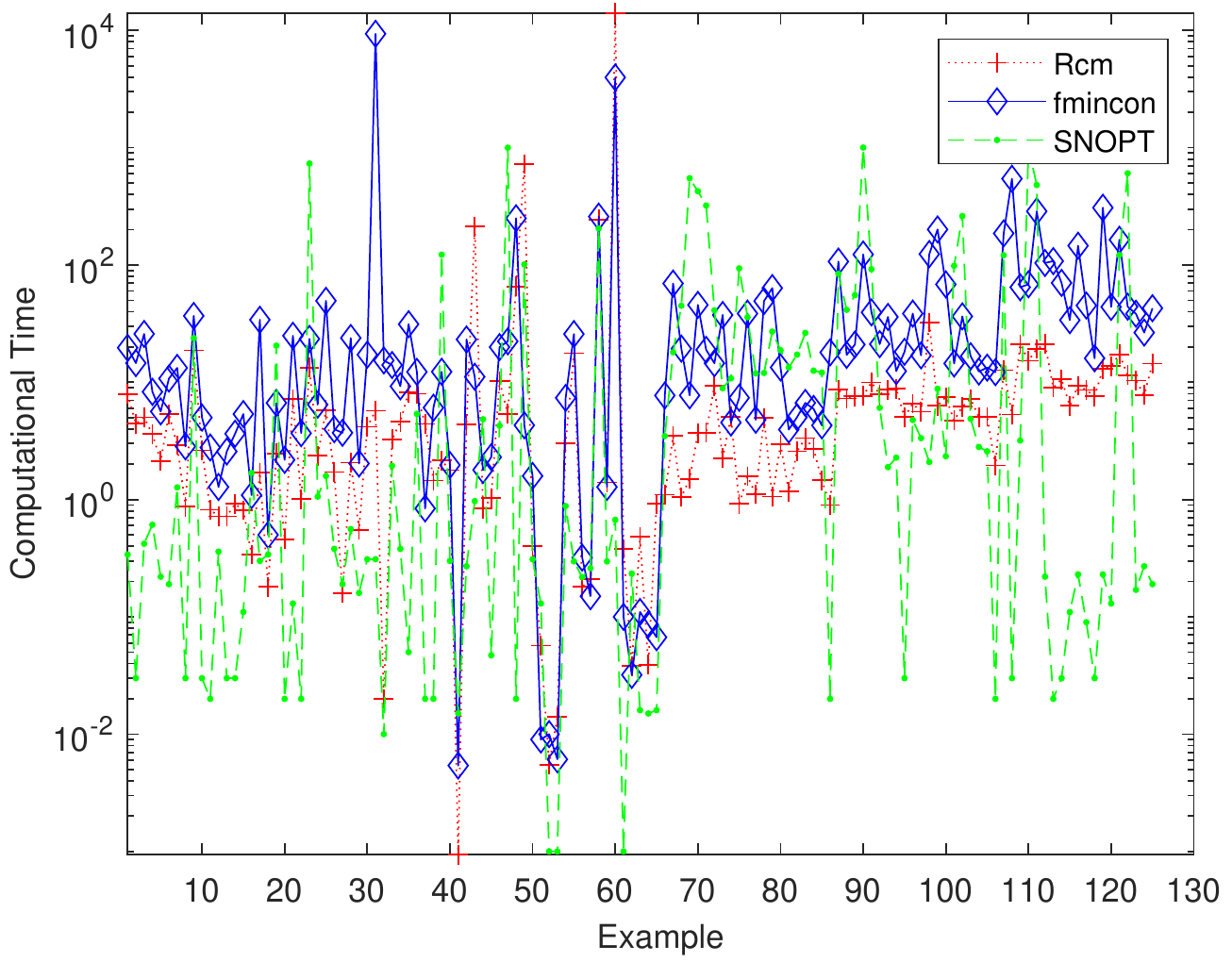}
  \captionof{figure}{CPU time of Rcm, fmincon and SNOPT for test problems.}
  \label{FIGCPURCM}
\end{minipage}

\vskip 2mm

\begin{minipage}{.9\textwidth}
  \centering
  \includegraphics[width=.9\linewidth]{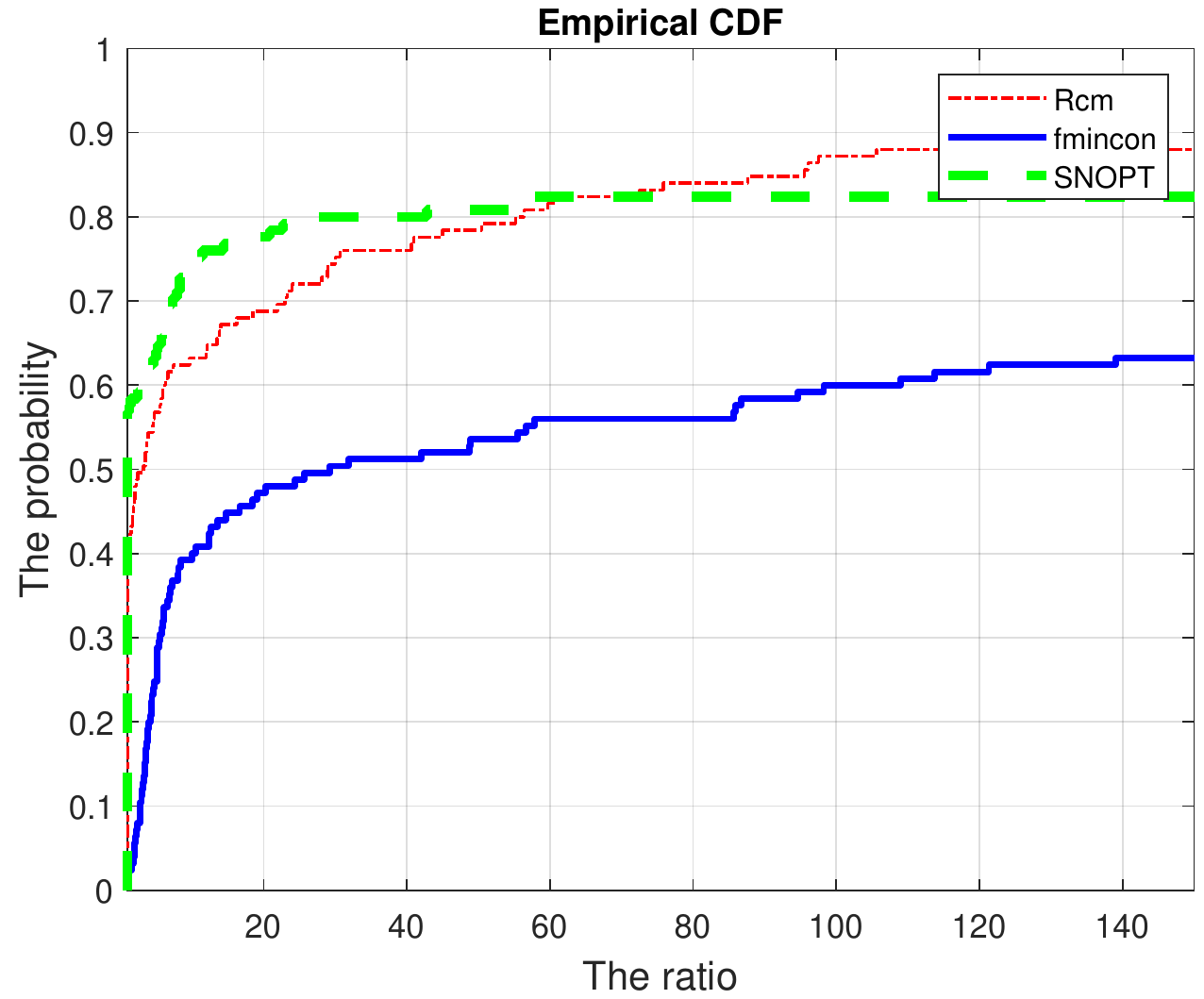}
  \captionof{figure}{Performance profile for nonlinear programming solvers.}
  \label{FIGPPRCM}
\end{minipage}

\end{figure}

The termination conditions of three compared algorithms (Rcm, fmincon and SNOPT)
are all set by
\begin{align}
    & \|\nabla_{x} L(x_{k}, \, \lambda_{k})\|_{\infty} \le 1.0 \times 10^{-6},
    \label{FOOPTTOL} \\
   & \|c(x_{k})\|_{\infty} \le 1.0 \times 10^{-6}, \;
    k = 1, \, 2, \, \ldots,   \label{FEATOL}
\end{align}
where the Lagrange function $L(x, \, \lambda)$ is defined by equation \eqref{LAGFUN}
and $\lambda$ is defined by equation \eqref{LAMBDA}. We test those $100$ problems
with $n \approx 1000$ to $n = 2000$. Numerical results are arranged in Tables
\ref{TABCOMRSA1}-\ref{TABCOMRSA5}. Numerical results of 65 constrained
optimization problems from CUTEst \cite{GOT2015} are arranged in Tables
\ref{TABCOMRSA1}-\ref{TABCOMRSA2}. And numerical results of 60 constructed constrained
optimization problems are arranged in Tables \ref{TABCOMRSA3}-\ref{TABCOMRSA5}.
The computational time and the number of iterations of
Rcm, fmincon and SNOPT are illustrated in Figure \ref{FIGITERCM}-\ref{FIGCPURCM},
respectively.

\vskip 2mm

Since Algorithm \ref{ALGRCM} (Rcm) use Algorithm \ref{ALGGCNMTR} (GCNMTr) to
find an initial feasible point, we count its iterations including the iterations
of GCNMTr (which are represented by $itc_{G}$) and the self iterations of Rcm
(which are represented by $itc_{R}$) in Tables \ref{TABCOMRSA1}-\ref{TABCOMRSA5}.
SNOPT gives the number of solving QP subproblems (which are represented by
``Major'') and uses SQOPT solver \cite{GMS2006} to find a solution of the QP
subproblem. SNOPT also gives the total iterations of solving QP subproblems,
which are represented by ``ICTS''. In order to ensure the evaluation
objectivity, we count the number of solving QP subproblems and the total
iterations of solving QP subproblems for SNOPT in Tables
\ref{TABCOMRSA1}-\ref{TABCOMRSA5}. fmincon only gives the number of solving
QP subproblems and we denote them as ``steps'' in Tables
\ref{TABCOMRSA1}-\ref{TABCOMRSA5}.

\vskip 2mm

\textcolor{blue}{In order to evaluate and compare those three methods (Rcm, fmincon and SNOPT)
fairly, we also adopt the performance profile as a evaluation tool \cite{DM2002}.
The performance profile for a solver is the (cumulative) distribution function
for a performance metric, which is the ratio of the computing time of the solver
versus the best time of all of the solvers as the performance metric. If the
solver fails to solve a problem, we let its ratio be a bigger number such as
$999$. Figure \ref{FIGPPRCM} is the performance profiles for three nonlinear
programming solvers (Rcm, fmincon and SNOPT). We also count the statistic
number of failed problems and fasted problems computed by Rcm, fminincon and
SNOPT. The statistical results are put in Table 7.  }

\vskip 2mm

From Tables \ref{TABCOMRSA1}-\ref{TABCOMRSA5} and Table \ref{TABPFRCM}, we find
that Rcm can solve most of test problems and few problems ($6.4\%$) can
not be solved. However, fmincon fails to solve CUTEst problems about
$\frac{28}{125} \; (22.4\%)$ and SNOPT fails to solve problems about
$\frac{21}{125} \; (16.8\%)$. From Tables \ref{TABCOMRSA1}-\ref{TABCOMRSA5} and
Figure \ref{FIGCPURCM}, Figure \ref{FIGPPRCM}, we find that the computational
time of Rcm is significantly less than that of fmincon
for most of test problems. From Tables \ref{TABCOMRSA1}-\ref{TABCOMRSA5} and
Figure \ref{FIGITERCM}, we also find that the iterations of Rcm is significantly
less than that of SNOPT for most of test problems. One of reasons is that Rcm
only needs to solve a linear system of equations with dimension $n$  at every
iteration. SQP needs to solve a linear system of equations with dimension $(m+n)$
when it solves a quadratic programming subproblem at every iteration
(pp. 531-532, \cite{NW1999}) and involves about $\frac{2}{3}(m+n)^{3}$ flops
(p. 116, \cite{GV2013}). The other reason is that Rcm uses the adaptive updating
technique of the Jacobian matrix $A(x)$ of constraints $c(x) = 0$ when it uses
Algorithm \ref{ALGGCNMTR} (GCNMTr) to find an initial feasible point. This
technique can significantly save the computational time when GCNMTr solves
an under-determined system of nonlinear equations \cite{LX2021}.
\textcolor{blue}{We also find that SNOPT is significantly faster than Rcm on
several test cases. The reason is due to the more gain of the complied language
than that of the interpreted language, since SNOPT is written by C++ language
(a complied language), and Rcm is written by MATLAB language (an interpreted
language).}

\vskip 2mm

\newpage

\vskip 2mm

\clearpage

\newpage

\begin{table}[!http]
  \newcommand{\tabincell}[2]{\begin{tabular}{@{}#1@{}}#2\end{tabular}}
  \scriptsize
  \centering
  \caption{Numerical results of Rcm, fmincon and SNOPT for CUTEst problems (Exam. 1-20).} \label{TABCOMRSA1}
  \resizebox{\textwidth}{!}{
  \begin{tabular}{|c|c|c|c|c|c|c|c|c|c|}
  \hline
  \multirow{2}{*}{Problems }
  & \multicolumn{2}{c|}{Rcm}
  & \multicolumn{2}{c|}{fmincon}
  & \multicolumn{2}{c|}{SNOPT}  \\ \cline{2-7}
        & \tabincell{c}{$itc_{G}+itc_{R}$ \\(time)} & \tabincell{c}{KKT \\ $||Cons(x^{it})||_\infty$}
        & \tabincell{c}{steps \\(time)}   & \tabincell{c}{KKT \\ $||Cons(x^{it})||_\infty$}
        & \tabincell{c}{Major+ICTS \\(time)} & \tabincell{c}{KKT \\ $||Cons(x^{it})||_\infty$} \\ \hline

  \tabincell{c}{Exam. 1 LUKVLE1 \cite{GOT2015}  \\ (n = 1000, m = 998)}
  & \tabincell{c}{4+10 \\ (7.93)} &\tabincell{c}{4.59e-06 \\ (2.66e-15)}
  & \tabincell{c}{{32} \\ {(19.68)}} &\tabincell{c}{{1.32e-04} \\ {(1.78e-15)} }
  & \tabincell{c}{17+23 \\ (0.34)} &\tabincell{c}{1.40e-07\\ (8.55e-13)} \\ \hline

  \tabincell{c}{Exam. 2 LUKVLE2 \cite{GOT2015}\\ (n = 1000, m = 993)}
  & \tabincell{c}{5+27 \\ (4.43)} &\tabincell{c}{3.96e-06 \\ (2.50e-16)}
  & \tabincell{c}{{38} \\ {(14.57)} } &\tabincell{c}{{2.67e-04} \\ {(5.82e-11)}  }
  & \tabincell{c}{\textcolor{red}{3+10} \\ \textcolor{red}{(0.03)}} &\tabincell{c}{\textcolor{red}{1.70e+01} \\ \textcolor{red}{(2.51e+01)} \\ \textcolor{red}{(failed)}} \\ \hline

  \tabincell{c}{Exam. 3 LUKVLE14  \cite{GOT2015}\\ (n = 1001, m = 666)}
  & \tabincell{c}{4+41 \\ (5.02)} &\tabincell{c}{7.41e-06 \\ (1.78e-15)}
  & \tabincell{c}{\textcolor{red}{99} \\ \textcolor{red}{(25.84)}} &\tabincell{c}{\textcolor{red}{1.92e-03} \\ \textcolor{red}{(7.52e-08)} \\ \textcolor{red}{(failed)}}
  & \tabincell{c}{25+361 \\ (0.42)} &\tabincell{c}{1.80e-06\\ (5.06e-13)} \\ \hline

  \tabincell{c}{Exam. 4 LUKVLE11 \cite{GOT2015} \\ (n = 1001, m = 666)}
  & \tabincell{c}{5+17 \\ (3.61)} &\tabincell{c}{4.41e-06 \\ (8.88e-16)}
  & \tabincell{c}{20 \\ (8.29)} &\tabincell{c}{4.00e-06 \\ (8.88e-16)}
  & \tabincell{c}{{35+405} \\ {(0.61)}} &\tabincell{c}{{1.20e-06} \\ {(3.84e-10)} } \\ \hline

  \tabincell{c}{Exam. 5 LUKVLE16 \cite{GOT2015}\\ (n = 1001, m = 750)}
  & \tabincell{c}{3+13 \\ (2.13)} &\tabincell{c}{8.47e-06 \\ (8.88e-16)}
  & \tabincell{c}{17 \\ (5.64)} &\tabincell{c}{3.32e-06 \\ (7.54e-12)}
  & \tabincell{c}{18+269 \\ (0.22)} &\tabincell{c}{2.90e-07\\ (3.41e-13)} \\ \hline

  \tabincell{c}{Exam. 6 LUKVLE17 \cite{GOT2015}\\ (n = 1001, m = 750)}
  & \tabincell{c}{4+35 \\ (5.34)} &\tabincell{c}{4.50e-06 \\ (2.22e-16)}
  & \tabincell{c}{25 \\ (10.76)} &\tabincell{c}{7.08e-06 \\ (4.46e-11)}
  & \tabincell{c}{15+266 \\ (0.19)} &\tabincell{c}{3.90e-07\\ (6.18e-09)} \\ \hline

  \tabincell{c}{Exam. 7 LUKVLE9 \cite{GOT2015}\\ (n = 1000, m = 6)}
  & \tabincell{c}{1+50 \\ (2.91)} &\tabincell{c}{9.90e-06 \\ (7.92e-07)}
  & \tabincell{c}{\textcolor{red}{99} \\ \textcolor{red}{(13.27)}} &\tabincell{c}{\textcolor{red}{3.78} \\ \textcolor{red}{(8.20e-02)}  \\ \textcolor{red}{(failed)}}
  & \tabincell{c}{\textcolor{red}{2+506} \\ \textcolor{red}{(1.27)}} &\tabincell{c}{\textcolor{red}{4.00e+00} \\ \textcolor{red}{(3.10e-02)}  \\ \textcolor{red}{(failed)}} \\ \hline

  \tabincell{c}{Exam. 8 BDVALUE \cite{GOT2015}\\ (n = 1000, m = 1000)}
  & \tabincell{c}{2+0 \\ (0.87)} &\tabincell{c}{0 \\ (2.08e-08)}
  & \tabincell{c}{4 \\ (2.84)} &\tabincell{c}{0 \\ (4.83e-16)}
  & \tabincell{c}{5+1 \\ (0.03)} &\tabincell{c}{9.10e-10 \\ (1.72e-10)} \\ \hline

  \tabincell{c}{Exam. 9 BDQRTIC \cite{GOT2015}\\ (n = 1000, m = 996)}
  & \tabincell{c}{\textcolor{red}{6+5} \\ \textcolor{red}{(18.68)}} &\tabincell{c}{\textcolor{red}{2.81e-02} \\ \textcolor{red}{(1.10)}  \\ \textcolor{red}{(failed)}}
  & \tabincell{c}{\textcolor{red}{94} \\ \textcolor{red}{(36.80)}} &\tabincell{c}{\textcolor{red}{1.36} \\ \textcolor{red}{(3.30)}  \\ \textcolor{red}{(failed)}}
    & \tabincell{c}{\textcolor{red}{222+4185} \\ \textcolor{red}{(23.86)}} &\tabincell{c}{\textcolor{red}{1.60e-06} \\ \textcolor{red}{(1.62e+00)}  \\ \textcolor{red}{(failed)}}
\\ \hline

  \tabincell{c}{Exam. 10 BDEXP \cite{GOT2015}\\(n = 1000,  m = 998)}
  & \tabincell{c}{21+0 \\ (2.63)} &\tabincell{c}{0 \\ (1.46e-07)}
  & \tabincell{c}{12 \\ (5.00)} &\tabincell{c}{0 \\ (9.53e-07)}
  & \tabincell{c}{15+34 \\ (0.03)} &\tabincell{c}{1.80e-06 \\ (3.10e-07)} \\ \hline

  \tabincell{c}{Exam. 11 BROYDNBD \cite{GOT2015}\\(n = 1000,  m = 1000)}
  & \tabincell{c}{4+0 \\ (0.82) } &\tabincell{c}{0 \\ (1.47e-09) }
  & \tabincell{c}{6 \\ (2.78)} &\tabincell{c}{ 0 \\ (8.88e-16)}
  & \tabincell{c}{7+0 \\ (0.02)} &\tabincell{c}{ 1.60e-13 \\ (4.96e-10) } \\ \hline

  \tabincell{c}{Exam. 12 YAO \cite{GOT2015}\\ (n = 1000,  m = 1000)}
  & \tabincell{c}{1+7 \\ (0.72)} &\tabincell{c}{1.33e-06 \\ (0)}
  & \tabincell{c}{4 \\ (1.28)} &\tabincell{c}{4.14e-08 \\ (0)}
  & \tabincell{c}{7+1007 \\ (0.36)} &\tabincell{c}{5.40e-07 \\ (0)} \\ \hline

  \tabincell{c}{Exam. 13 WOODSNE \cite{GOT2015}\\ (n = 1000,  m = 1000)}
  & \tabincell{c}{4+0 \\ (0.72)} &\tabincell{c}{0\\ (7.23e-09) \\ }
  & \tabincell{c}{6 \\ (2.57)} &\tabincell{c}{0 \\ (1.27e-12)}
  & \tabincell{c}{7+500 \\ (0.03)} &\tabincell{c}{8.00e-07\\ (1.16e-07) \\ } \\ \hline

  \tabincell{c}{Exam. 14 WOODS  \cite{GOT2015}\\ (n = 1000,  m = 1000)}
  & \tabincell{c}{6+0 \\ (0.92)} &\tabincell{c}{3.33e-16\\ (7.04e-08) \\ }
  & \tabincell{c}{9 \\ (3.64)} &\tabincell{c}{1.50e-08 \\ (8.96e-14)}
  & \tabincell{c}{10+500 \\ (0.03)} &\tabincell{c}{5.10e-10 \\ (1.68e-08) \\ } \\ \hline

  \tabincell{c}{Exam. 15 SOSQP1 \cite{GOT2015}\\ (n = 1000,  m = 501)}
  & \tabincell{c}{1+0 \\ (0.81)} &\tabincell{c}{6.66e-15 \\ (4.15e-10)}
  & \tabincell{c}{2 \\ (5.36)} &\tabincell{c}{1.10e-13 \\ (1.14e-13)}
  & \tabincell{c}{3+505 \\ (0.11)} &\tabincell{c}{9.80e-12 \\ (1.31e-11)} \\ \hline

  \tabincell{c}{Exam. 16 ARGLALE \cite{GOT2015}\\ (n = 1000,  m = 1000)}
  & \tabincell{c}{1+0 \\ (0.34)} &\tabincell{c}{0 \\ (1.65e-08)}
  & \tabincell{c}{2 \\ (1.09)} &\tabincell{c}{0 \\ (3.29e-11)}
  & \tabincell{c}{0+0 \\ (1.69)} &\tabincell{c}{0 \\ (1.60e-12)} \\ \hline

  \tabincell{c}{Exam. 17 AUG2D \cite{GOT2015}\\ (n = 1000,  m = 1000)}
  & \tabincell{c}{1+0 \\ (11.91)} &\tabincell{c}{0 \\ (0)}
  & \tabincell{c}{1 \\ (34.11)} &\tabincell{c}{0 \\ (0)}
  & \tabincell{c}{0+0 \\ (0.3)} &\tabincell{c}{0.00e+0 \\ (0)} \\ \hline

  \tabincell{c}{Exam. 18 BLOCKQP1  \cite{GOT2015}\\ (n = 501,  m = 1005)}
  & \tabincell{c}{1+0 \\ (0.18)} &\tabincell{c}{5.08e-15 \\ (9.99e-08)}
  & \tabincell{c}{2 \\ (0.50)} &\tabincell{c}{1.55e-14 \\ (4.44e-15)}
  & \tabincell{c}{\textcolor{red}{50+602} \\ \textcolor{red}{(0.34)}} &\tabincell{c}{\textcolor{red}{2.70e-12} \\ \textcolor{red}{(1.65e+07)}  \\ \textcolor{red}{(failed)}}
  \\ \hline

  \tabincell{c}{Exam. 19 BROWNAL \cite{GOT2015}\\ (n = 1000,  m = 1000)}
  & \tabincell{c}{2+0 \\ (2.48)} &\tabincell{c}{6.46e-14 \\ (8.82e-08)}
  & \tabincell{c}{1 \\ (6.67)} &\tabincell{c}{6.63e-06 \\ (1.50e-06)}
  & \tabincell{c}{22+520 \\ (20.53)} &\tabincell{c}{2.50e-08 \\ (2.36e-11)} \\ \hline

  \tabincell{c}{Exam. 20 BROYDEN3D  \cite{GOT2015}\\ (n = 1000,  m = 1000)}
  & \tabincell{c}{2+0 \\ (0.46)} &\tabincell{c}{0  \\ (3.72e-08)}
  & \tabincell{c}{5 \\ (2.18)} &\tabincell{c}{5.08e-11 \\ (1.11e-15)}
  & \tabincell{c}{5+2 \\ (0.02)} &\tabincell{c}{1.70e-08 \\ (1.56e-05)} \\ \hline

\end{tabular}}
\end{table}

\newpage

\vskip 2mm

\begin{table}[!http]
  \newcommand{\tabincell}[2]{\begin{tabular}{@{}#1@{}}#2\end{tabular}}
  \scriptsize
  \centering
  \caption{Numerical results of Rcm, fmincon and SNOPT for CUTEst problems (Exam. 21-40).} \label{TABCOMRSA2}
  \resizebox{\textwidth}{!}{
  \begin{tabular}{|c|c|c|c|c|c|c|c|c|c|}
  \hline
  \multirow{2}{*}{Problems }
  & \multicolumn{2}{c|}{Rcm}
  & \multicolumn{2}{c|}{fmincon}
  & \multicolumn{2}{c|}{SNOPT}  \\ \cline{2-7}
        & \tabincell{c}{$itc_{G}+itc_{R}$ \\(time)} & \tabincell{c}{KKT \\ $||Cons(x^{it})||_\infty$}
        & \tabincell{c}{steps \\(time)}   & \tabincell{c}{KKT \\ $||Cons(x^{it})||_\infty$}
        & \tabincell{c}{Major+ICTS \\(time)} & \tabincell{c}{KKT \\ $||Cons(x^{it})||_\infty$} \\ \hline

  \tabincell{c}{Exam. 21 CHENHARK \cite{GOT2015} \\ (n = 1000,  m = 1000)}
  & \tabincell{c}{1+0 \\ (7.18)} &\tabincell{c}{9.76e-15  \\ (8.25e-09)}
  & \tabincell{c}{2 \\ (24.94)} &\tabincell{c}{1.55e-06 \\ (1.55e-15)}
  & \tabincell{c}{36+682 \\ (0.13)} &\tabincell{c}{3.00e-10 \\ (4.75e-07)} \\ \hline

  \tabincell{c}{Exam. 22 CHNROSNB \cite{GOT2015} \\ (n = 1000,  m = 1000)}
  & \tabincell{c}{7+0 \\ (1.01)} &\tabincell{c}{1.55e-15  \\ (4.94e-08)}
  & \tabincell{c}{10 \\ (3.68)} &\tabincell{c}{1.85e-08 \\ (7.99e-15)}
  & \tabincell{c}{11+2 \\ (0.02)} &\tabincell{c}{1.40e-08\\ (7.85e-08)} \\ \hline

  \tabincell{c}{Exam. 23 DIXON3DQ \cite{GOT2015} \\ (n = 1000, m = 0)}
  & \tabincell{c}{0+83 \\ (13.34)} &\tabincell{c}{9.91e-06 \\ (0)}
  & \tabincell{c}{\textcolor{red}{99} \\ \textcolor{red}{(23.44)}} &\tabincell{c}{\textcolor{red}{6.66e-02} \\ \textcolor{red}{(0)}  \\ \textcolor{red}{(failed)}}
  & \tabincell{c}{168685+184416 \\ (733.17)} &\tabincell{c}{1.80e-06\\ (0)} \\ \hline

  \tabincell{c}{Exam. 24 DIXMAANA \cite{GOT2015} \\ (n = 1000,  m = 0)}
  & \tabincell{c}{0+16 \\ (2.37)} &\tabincell{c}{3.19e-06  \\ (0)}
  & \tabincell{c}{16 \\ (6.45)} &\tabincell{c}{3.81e-06 \\ (0)}
  & \tabincell{c}{12+1216 \\ (1.06)} &\tabincell{c}{2.60e-07 \\ (0)} \\ \hline

  \tabincell{c}{Exam. 25 ORTHGDS \cite{GOT2015} \\ (n = 1003, m = 1000)}
  & \tabincell{c}{1+40 \\ (5.79)} &\tabincell{c}{3.76e-08 \\ (1.32e-07)}
  & \tabincell{c}{\textcolor{red}{95} \\ \textcolor{red}{(49.52)}} &\tabincell{c}{\textcolor{red}{5.31e+02} \\ \textcolor{red}{(9.61e+02)}  \\ \textcolor{red}{(failed)}}
  & \tabincell{c}{41+549 \\ (1.59)} &\tabincell{c}{1.20e-06\\ (2.45e-11)} \\ \hline

  \tabincell{c}{Exam. 26 VARDIM \cite{GOT2015} \\ (n = 1000,  m = 0)}
  & \tabincell{c}{1+10 \\ (1.72)} &\tabincell{c}{1.60e-14 \\ (0)}
  & \tabincell{c}{16 \\ (3.98)} &\tabincell{c}{7.07e-06 \\ (0)}
  & \tabincell{c}{3+1003 \\ (0.38)} &\tabincell{c}{9.30e-06\\ (0)} \\ \hline

  \tabincell{c}{Exam. 27 SVANBERG \cite{GOT2015} \\ (n = 1000, m = 1000)}
  & \tabincell{c}{24+33 \\ (0.16)} &\tabincell{c}{6.04e-06 \\ (6.34e-07)}
  & \tabincell{c}{\textcolor{red}{88} \\ \textcolor{red}{(3.71)}} &\tabincell{c}{\textcolor{red}{8.64e+01} \\ \textcolor{red}{(2.12)}  \\ \textcolor{red}{(failed)}}
  & \tabincell{c}{\textcolor{red}{0+1050} \\ \textcolor{red}{(0.19)}} &\tabincell{c}{\textcolor{red}{2.00e+0} \\ \textcolor{red}{(2.10e+01)}  \\ \textcolor{red}{(failed)}} \\ \hline

  \tabincell{c}{Exam. 28 SROSENBR \cite{GOT2015} \\ (n = 1000, m = 0)}
  & \tabincell{c}{0+24 \\ (2.07)} &\tabincell{c}{1.88e-06 \\ (0)}
  & \tabincell{c}{\textcolor{red}{98} \\ \textcolor{red}{(23.87)}} &\tabincell{c}{\textcolor{red}{8.43e+02} \\ \textcolor{red}{(0)}  \\ \textcolor{red}{(failed)}}
  & \tabincell{c}{21+1036 \\ (0.56)} &\tabincell{c}{3.70e-08\\ (0)} \\ \hline

  \tabincell{c}{Exam. 29 SREADIN3 \cite{GOT2015} \\ (n = 1000,  m = 1000)}
  & \tabincell{c}{1+0 \\ (0.55)} &\tabincell{c}{2.08e-06  \\ (2.08e-08)}
  & \tabincell{c}{3 \\ (2.03)} &\tabincell{c}{2.08e-06 \\ (3.06e-12)}
  & \tabincell{c}{2+497 \\ (0.16)} &\tabincell{c}{1.90e-06 \\ (7.43e-06)} \\ \hline

  \tabincell{c}{Exam. 30 SOSQP2 \cite{GOT2015} \\ (n = 1000, m = 1000)}
  & \tabincell{c}{1+50 \\ (4.18)} &\tabincell{c}{3.84e-06 \\ (7.40e-07)}
  & \tabincell{c}{{99} \\ {(17.20)}} &\tabincell{c}{{1.03e-04} \\ {(1.33e-11)} }
  & \tabincell{c}{38+538 \\ (0.31)} &\tabincell{c}{2.00e-06\\ (2.81e-11)} \\ \hline

  \tabincell{c}{Exam. 31 SINROSNB \cite{GOT2015} \\ (n = 1000, m = 1000)}
  & \tabincell{c}{2+7 \\ (5.73)} &\tabincell{c}{5.17e-06 \\ (2.22e-16)}
  & \tabincell{c}{\textcolor{red}{69} \\ \textcolor{red}{(9.36e+03)}} &\tabincell{c}{\textcolor{red}{4.18e+04} \\ \textcolor{red}{(2.12e-02)}  \\ \textcolor{red}{(failed)}}
  & \tabincell{c}{18+22 \\ (0.31)} &\tabincell{c}{5.400e-07\\ (2.44e-13)} \\ \hline

  \tabincell{c}{Exam. 32 SIPOW1 \cite{GOT2015} \\ (n = 1000,  m = 2)}
  & \tabincell{c}{1+0 \\ (0.02)} &\tabincell{c}{8.70e-19 \\ (1.16e-08)}
  & \tabincell{c}{2 \\ (15.35)} &\tabincell{c}{2.83e-06 \\ (1.96e-18)}
  & \tabincell{c}{0+0 \\ (0.01)} &\tabincell{c}{0\\ (0)} \\ \hline

  \tabincell{c}{Exam. 33 SINQUAD \cite{GOT2015} \\ (n = 1000,  m = 0)}
  & \tabincell{c}{0+39 \\ (3.26)} &\tabincell{c}{8.65e-06 \\ (0)}
  & \tabincell{c}{57 \\ (14.03)} &\tabincell{c}{2.17e-06 \\ (0)}
  & \tabincell{c}{247+1419 \\ (1.94)} &\tabincell{c}{1.70e-06\\ (0)} \\ \hline

  \tabincell{c}{Exam. 34 SINEALI \cite{GOT2015} \\ (n = 1000, m = 0)}
  & \tabincell{c}{0+42 \\ (4.61)} &\tabincell{c}{8.80e-16 \\ (0)}
  & \tabincell{c}{\textcolor{red}{38} \\ \textcolor{red}{(9.24)}} &\tabincell{c}{\textcolor{red}{6.78e+03} \\ \textcolor{red}{(0)}  \\ \textcolor{red}{(failed)}}
  & \tabincell{c}{16+1017 \\ (0.38)} &\tabincell{c}{8.80e-07\\ (0)} \\ \hline

  \tabincell{c}{Exam. 35 SEMICON1 \cite{GOT2015} \\ (n = 1000,  m = 1000)}
  & \tabincell{c}{74+0 \\ (8.19)} &\tabincell{c}{0  \\ (1.63e-08)}
  & \tabincell{c}{88 \\ (31.28)} &\tabincell{c}{0 \\ (3.80e-13)}
  & \tabincell{c}{\textcolor{red}{9+1} \\ \textcolor{red}{(0.05)}} &\tabincell{c}{\textcolor{red}{1.00e-01} \\ \textcolor{red}{(6.79e+02)}  \\ \textcolor{red}{(failed)}}  \\ \hline

  \tabincell{c}{Exam. 36 READING1 \cite{GOT2015} \\ (n = 1000, m = 999)}
  & \tabincell{c}{1+19 \\ (7.98)} &\tabincell{c}{4.26e-06 \\ (5.28e-08)}
  & \tabincell{c}{\textcolor{red}{24} \\ \textcolor{red}{(12.27)}} &\tabincell{c}{\textcolor{red}{1.25e-03} \\ \textcolor{red}{(1.04e-02)}  \\ \textcolor{red}{(failed)}}
  & \tabincell{c}{\textcolor{red}{104+5740} \\ \textcolor{red}{(5.38)}} &\tabincell{c}{\textcolor{red}{1.10e+03} \\ \textcolor{red}{(3.11e+05)}  \\ \textcolor{red}{(failed)}} \\ \hline

  \tabincell{c}{Exam. 37 POWELL20 \cite{GOT2015} \\ (n = 1000,  m = 1000)}
  & \tabincell{c}{1+6 \\ (4.41)} &\tabincell{c}{3.10e-07  \\ (1.14e-13)}
  & \tabincell{c}{1 \\ (0.84)} &\tabincell{c}{3.46e-11 \\ (3.67e-12)}
  & \tabincell{c}{3+3 \\ (0.02)} &\tabincell{c}{2.27e-13 \\ ()} \\ \hline

  \tabincell{c}{Exam. 38 LUKVLE3 \cite{GOT2015} \\ (n = 1000, m = 2)}
  & \tabincell{c}{3+43 \\ (1.45)} &\tabincell{c}{9.52e-08 \\ (5.98e-07)}
  & \tabincell{c}{{56} \\ {(6.07)}} &\tabincell{c}{{2.67e-04}}
  & \tabincell{c}{24+36 \\ (0.02)} &\tabincell{c}{1.00e-06 \\ (3.52e-13)}  \\ \hline

  \tabincell{c}{Exam. 39 LUKVLE7 \cite{GOT2015} \\ (n = 1000, m = 4)}
  & \tabincell{c}{4+54 \\ (2.15)} &\tabincell{c}{1.10e-06 \\ (6.03e-07)}
  & \tabincell{c}{\textcolor{red}{98} \\ \textcolor{red}{(12.30)}} &\tabincell{c}{\textcolor{red}{1.08e+05} \\ \textcolor{red}{(2.66e+03)}  \\ \textcolor{red}{(failed)}}
  & \tabincell{c}{{1150+2210} \\ {(122.55)}} &\tabincell{c}{{4.10e-04} \\ {(6.85e-12)} } \\ \hline

  \tabincell{c}{Exam. 40 LUKVLE12 \cite{GOT2015} \\ (n = 1001,  m = 750)}
  & \tabincell{c}{2+11 \\ (1.78)} &\tabincell{c}{7.24e-06 \\ (8.88e-16)}
  & \tabincell{c}{6 \\ (1.97)} &\tabincell{c}{1.63e-06 \\ (3.41e-10)}
  & \tabincell{c}{{20+524} \\ {(0.30)}} &\tabincell{c}{{1.70e-07} \\ {(2.62e-12)} } \\ \hline

\end{tabular}}
\end{table}

\vskip 2mm

\begin{table}[!http]
  \newcommand{\tabincell}[2]{\begin{tabular}{@{}#1@{}}#2\end{tabular}}
  \scriptsize
  \centering
  \caption{Numerical results of Rcm, fmincon and SNOPT for CUTEst problems (Exam. 41-65).} \label{TABCOMRSA2}
  \resizebox{\textwidth}{!}{
  \begin{tabular}{|c|c|c|c|c|c|c|c|c|c|}
  \hline
  \multirow{2}{*}{Problems }
  & \multicolumn{2}{c|}{Rcm}
  & \multicolumn{2}{c|}{fmincon}
  & \multicolumn{2}{c|}{SNOPT}  \\ \cline{2-7}
        & \tabincell{c}{$itc_{G}+itc_{R}$ \\(time)} & \tabincell{c}{KKT \\ $||Cons(x^{it})||_\infty$}
        & \tabincell{c}{steps \\(time)}   & \tabincell{c}{KKT \\ $||Cons(x^{it})||_\infty$}
        & \tabincell{c}{Major+ICTS \\(time)} & \tabincell{c}{KKT \\ $||Cons(x^{it})||_\infty$} \\ \hline

  \tabincell{c}{Exam. 41 aircrfta \cite{GOT2015} \\ (n = 8,  m = 5)}
  & \tabincell{c}{1+0 \\ (0.00094)} &\tabincell{c}{0.00e+00  \\ (0.00e+00)}
  & \tabincell{c}{0 \\ (0.0054)} &\tabincell{c}{0.00e+00  \\ (0.00e+00)}
  & \tabincell{c}{0+0 \\ (0.015)} &\tabincell{c}{0.00e+00  \\ (0.00e+00)} \\ \hline

  \tabincell{c}{Exam. 42 orthrdm2 \cite{GOT2015} \\ (n = 2003,  m = 1000)}
  & \tabincell{c}{4+10 \\ (4.39)} &\tabincell{c}{7.66E-06  \\ (8.53E-14)}
  & \tabincell{c}{13 \\ (23.19)} &\tabincell{c}{2.34E-06  \\ (2.93E-10)}
  & \tabincell{c}{6+509 \\ (0.27)} &\tabincell{c}{2.60E-07  \\ (1.13E-06)} \\ \hline

  \tabincell{c}{Exam. 43 huestis \cite{GOT2015} \\ (n = 1000,  m = 2)}
  & \tabincell{c}{\textcolor{red}{1+360} \\ \textcolor{red}{(212.67)}} &\tabincell{c}{\textcolor{red}{472.58} \\ \textcolor{red}{(2.76E-08)} \\ \textcolor{red}{(failed)}}
  & \tabincell{c}{\textcolor{red}{32} \\ \textcolor{red}{(11.14)}} &\tabincell{c}{\textcolor{red}{1.24E+03} \\ \textcolor{red}{(6.67E-06)}  \\ \textcolor{red}{(failed)}}
  & \tabincell{c}{5+1992 \\ (0.97)} &\tabincell{c}{3.70E-14  \\ (9.09E-13)} \\ \hline

  \tabincell{c}{Exam. 44 gilbert \cite{GOT2015} \\ (n = 1000,  m = 1)}
  & \tabincell{c}{6+28 \\ (0.84)} &\tabincell{c}{8.58E-06  \\ (2.24E-14)}
  & \tabincell{c}{22 \\ (1.79)} &\tabincell{c}{8.58E-06  \\ (9.34E-08)}
  & \tabincell{c}{31+1035 \\ (4.84)} &\tabincell{c}{1.00E-06  \\ (3.87E-09)} \\ \hline

  \tabincell{c}{Exam. 45 genhs28 \cite{GOT2015} \\ (n = 1000,  m = 998)}
  & \tabincell{c}{1+7 \\ (1.03)} &\tabincell{c}{3.32E-06  \\ (2.22E-16)}
  & \tabincell{c}{7 \\ (2.30)} &\tabincell{c}{1.51E-06  \\ (2.22E-16)}
  & \tabincell{c}{4+9 \\ (0.047)} &\tabincell{c}{3.60E-12  \\ (4.55E-13)} \\ \hline

  \tabincell{c}{Exam. 46 dtoc6 \cite{GOT2015} \\ (n = 2001,  m = 1000)}
  & \tabincell{c}{3+13 \\ (10.24)} &\tabincell{c}{8.03E-06  \\ (1.63E-13)}
  & \tabincell{c}{11 \\ (19.95)} &\tabincell{c}{1.68E-06  \\ (2.32E-12)}
  & \tabincell{c}{20+1025 \\ (4.27)} &\tabincell{c}{1.70E-06  \\ (2.91E-07)} \\ \hline

  \tabincell{c}{Exam. 47 dtoc5 \cite{GOT2015} \\ (n = 2001,  m = 1000)}
  & \tabincell{c}{1+17 \\ (5.38)} &\tabincell{c}{9.72E-06  \\ (8.34E-12)}
  & \tabincell{c}{21 \\ (22.35)} &\tabincell{c}{6.31E-06  \\ (6.53E-09)}
  & \tabincell{c}{\textcolor{red}{155561+282454} \\ \textcolor{red}{(1000)}} &\tabincell{c}{\textcolor{red}{4.90E-03} \\ \textcolor{red}{(9.54E-01)}  \\ \textcolor{red}{(failed)}} \\ \hline

  \tabincell{c}{Exam. 48 orthregc \cite{GOT2015} \\ (n = 2005,  m = 1000)}
  & \tabincell{c}{8+60 \\ (65.61)} &\tabincell{c}{8.76E-06  \\ (8.07E-07)}
  & \tabincell{c}{\textcolor{red}{99} \\ \textcolor{red}{(250.32)}} &\tabincell{c}{\textcolor{red}{2.61E-01} \\ \textcolor{red}{(3.74E-04)}  \\ \textcolor{red}{(failed)}}
  & \tabincell{c}{1+3 \\ (0.02)} &\tabincell{c}{4.40E-16  \\ (2.42E-13)} \\ \hline

  \tabincell{c}{Exam. 49 catenary \cite{GOT2015} \\ (n = 501,  m = 166)}
  & \tabincell{c}{\textcolor{red}{1+397} \\ \textcolor{red}{(725.98)}} &\tabincell{c}{\textcolor{red}{2.95E+01} \\ \textcolor{red}{(4.30E+05)}  \\ \textcolor{red}{(failed)}}
  & \tabincell{c}{\textcolor{red}{85} \\ \textcolor{red}{(4.30)}} &\tabincell{c}{\textcolor{red}{2.09E+08} \\ \textcolor{red}{(1.78E+00)}  \\ \textcolor{red}{(failed)}}
  & \tabincell{c}{7092+9790 \\ (101.23)} &\tabincell{c}{2.00E-06  \\ (2.00E-11)} \\ \hline

  \tabincell{c}{Exam. 50 broydn3d \cite{GOT2015} \\ (n = 1000,  m = 1000)}
  & \tabincell{c}{2+0 \\ (0.40)} &\tabincell{c}{0 \\ (3.72E-08)}
    & \tabincell{c}{5 \\ (1.59)} &\tabincell{c}{5.08E-11 \\ (8.88E-16)}
  & \tabincell{c}{5+2 \\ (0.31)} &\tabincell{c}{1.70E-08  \\ (1.56E-05)} \\ \hline

  \tabincell{c}{Exam. 51 hs007 \cite{GOT2015} \\ (n = 2,  m = 1)}
  & \tabincell{c}{6+13 \\ (0.057)} &\tabincell{c}{2.39E-06  \\ (0.00E+00)}
  & \tabincell{c}{9 \\ (0.0090)} &\tabincell{c}{6.11E-08  \\ (3.47E-12)}
  & \tabincell{c}{20+21 \\ (0.13)} &\tabincell{c}{1.56E-05  \\ (3.70E-07)} \\ \hline

  \tabincell{c}{Exam. 52 hs008 \cite{GOT2015} \\ (n = 2,  m = 2)}
  & \tabincell{c}{3+0 \\ (0.0055)} &\tabincell{c}{0.00E+00  \\ (5.12E-08)}
  & \tabincell{c}{5 \\ (0.010)} &\tabincell{c}{1.32E-06  \\ (1.26E-10)}
  & \tabincell{c}{4+1 \\ (0.001)} &\tabincell{c}{4.00E-11  \\ (2.69E-07)} \\ \hline

  \tabincell{c}{Exam. 53 hs009 \cite{GOT2015} \\ (n = 2,  m = 1)}
  & \tabincell{c}{1+6 \\ (0.014)} &\tabincell{c}{3.84E-06  \\ (0.00E+00)}
  & \tabincell{c}{5 \\ (0.0061)} &\tabincell{c}{4.33E-06  \\ (1.78E-15)}
  & \tabincell{c}{5+6 \\ (0.001)} &\tabincell{c}{2.80E-12  \\ (2.27E-13)} \\ \hline

  \tabincell{c}{Exam. 54 gouldqp2 \cite{GOT2015} \\ (n = 699,  m = 349)}
  & \tabincell{c}{1+35 \\ (3.038)} &\tabincell{c}{8.01E-06  \\ (2.00E-11)}
  & \tabincell{c}{99 \\ (7.36)} &\tabincell{c}{2.19E-04 \\ (4.62E-12)}
  & \tabincell{c}{1400+1750 \\ (0.88)} &\tabincell{c}{1.90E-06  \\ (1.14E-13)} \\ \hline

  \tabincell{c}{Exam. 55 optcdeg2 \cite{GOT2015} \\ (n = 1202,  m = 800)}
  & \tabincell{c}{1+52 \\ (17.73)} &\tabincell{c}{8.79E-06  \\ (7.71E-07)}
  & \tabincell{c}{87 \\ (25.77)} &\tabincell{c}{8.62E-06 \\ (3.13E-10)}
  & \tabincell{c}{58+459 \\ (0.297)} &\tabincell{c}{1.30E-06  \\ (1.14E-13)} \\ \hline

  \tabincell{c}{Exam. 56 hs100lnp \cite{GOT2015} \\ (n = 7,  m = 2)}
  & \tabincell{c}{1+14 \\ (0.18)} &\tabincell{c}{4.26E-06 \\ (1.42E-14)}
  & \tabincell{c}{15 \\ (0.32)} &\tabincell{c}{7.32E-04 \\ (8.11E-08)}
  & \tabincell{c}{26+50 \\ (0.218)} &\tabincell{c}{1.10E-09  \\ (2.67E-12)} \\ \hline

  \tabincell{c}{Exam. 57 hs111 \cite{GOT2015} \\ (n = 10,  m = 3)}
  & \tabincell{c}{\textcolor{red}{4+29} \\ \textcolor{red}{(0.21)}} &\tabincell{c}{\textcolor{red}{1.63E-02} \\ \textcolor{red}{(4.38E-05)} \\ \textcolor{red}{(failed)}}
  & \tabincell{c}{\textcolor{red}{8} \\ \textcolor{red}{(0.15)}} &\tabincell{c}{\textcolor{red}{1.21E-01} \\ \textcolor{red}{(3.64E-03)}  \\ \textcolor{red}{(failed)}}
  & \tabincell{c}{50+57 \\ (0.26)} &\tabincell{c}{4.70E-07  \\ (1.40E-08)} \\ \hline

  \tabincell{c}{Exam. 58 reading1 \cite{GOT2015} \\ (n = 2002,  m = 1000)}
  & \tabincell{c}{\textcolor{red}{1+47} \\ \textcolor{red}{(243.60)}} &\tabincell{c}{\textcolor{red}{4.17E+09} \\ \textcolor{red}{(2.07E+13)} \\ \textcolor{red}{(failed)}}
  & \tabincell{c}{\textcolor{red}{99} \\ \textcolor{red}{(258.70)}} &\tabincell{c}{\textcolor{red}{3.02E+02} \\ \textcolor{red}{(4.34E+03)}  \\ \textcolor{red}{(failed)}}
  & \tabincell{c}{\textcolor{red}{827+38923} \\ \textcolor{red}{(205.094)}} &\tabincell{c}{\textcolor{red}{1.50E+04} \\ \textcolor{red}{(3.22E+07)} \\ \textcolor{red}{(failed)}} \\ \hline

  \tabincell{c}{Exam. 59 hs100lnp \cite{GOT2015} \\ (n = 4,  m = 1)}
  & \tabincell{c}{1+90 \\ (1.40)} &\tabincell{c}{9.91E-06 \\ (1.32E-09)}
  & \tabincell{c}{3 \\ (1.28)} &\tabincell{c}{4.60E-07 \\ (8.53E-14)}
  & \tabincell{c}{5+8 \\ (0.297)} &\tabincell{c}{4.30E-11 \\ (1.14E-13)} \\ \hline

  \tabincell{c}{Exam. 60 reading2 \cite{GOT2015} \\ (n = 3003,  m = 2002)}
  & \tabincell{c}{1+394 \\ (14094.66)} &\tabincell{c}{2.53E-05 \\ (5.45E-09)}
  & \tabincell{c}{28 \\ (3975.68)} &\tabincell{c}{2.53E-05 \\ (7.64E-03)}
  & \tabincell{c}{2600+2650 \\ (0.672)} &\tabincell{c}{0.00E+00 \\ (1.01E-10)} \\ \hline

  \tabincell{c}{Exam. 61 rk23 \cite{GOT2015} \\ (n = 17,  m = 11)}
  & \tabincell{c}{\textcolor{red}{2+394} \\ \textcolor{red}{(0.38)}} &\tabincell{c}{\textcolor{red}{1.00E+00} \\ \textcolor{red}{(2.57E-11)} \\ \textcolor{red}{(failed)}}
  & \tabincell{c}{\textcolor{red}{46} \\ \textcolor{red}{(0.1001)}} &\tabincell{c}{\textcolor{red}{3.05E+04} \\ \textcolor{red}{(1.94E+01)}  \\ \textcolor{red}{(failed)}}
  & \tabincell{c}{\textcolor{red}{0+3} \\ \textcolor{red}{(0.001)}} &\tabincell{c}{\textcolor{red}{1.50E-05} \\ \textcolor{red}{(5.00E-01)} \\ \textcolor{red}{(failed)}} \\ \hline

  \tabincell{c}{Exam. 62 hs046 \cite{GOT2015} \\ (n = 5,  m = 2)}
  & \tabincell{c}{1+32 \\ (0.038)} &\tabincell{c}{8.84E-06 \\ (6.23E-09)}
  & \tabincell{c}{25 \\ (0.032)} &\tabincell{c}{3.01E-06 \\ (1.42E-08)}
  & \tabincell{c}{39+42 \\ (0.234)} &\tabincell{c}{1.90E-06 \\ (2.11E-10)} \\ \hline

  \tabincell{c}{Exam. 63 hs099 \cite{GOT2015} \\ (n = 23,  m = 18)}
  & \tabincell{c}{1+32 \\ (0.4808)} &\tabincell{c}{8.84E-06 \\ (6.23E-09)}
  & \tabincell{c}{25 \\ (0.110)} &\tabincell{c}{3.01E-06 \\ (1.42E-08)}
  & \tabincell{c}{14+23 \\ (0.016)} &\tabincell{c}{1.30E-05 \\ (4.51E-04)} \\ \hline

  \tabincell{c}{Exam. 64 hs0991 \cite{GOT2015} \\ (n = 31,  m = 26)}
  & \tabincell{c}{\textcolor{red}{4+0} \\ \textcolor{red}{(0.039)}} &\tabincell{c}{\textcolor{red}{8.81E+02} \\ \textcolor{red}{(1.40E-08)} \\ \textcolor{red}{(failed)}}
  & \tabincell{c}{\textcolor{red}{7} \\ \textcolor{red}{(0.086)}} &\tabincell{c}{\textcolor{red}{1.06E+06} \\ \textcolor{red}{(9.90E+02)}  \\ \textcolor{red}{(failed)}}
  & \tabincell{c}{14+23 \\ (0.015)} &\tabincell{c}{9.00E-06 \\ (3.15E-04)} \\ \hline

  \tabincell{c}{Exam. 65 hs99exp \cite{GOT2015} \\ (n = 31,  m = 21)}
  & \tabincell{c}{\textcolor{red}{5+50} \\ \textcolor{red}{(0.92)}} &\tabincell{c}{\textcolor{red}{2.91E+05} \\ \textcolor{red}{(7.05E+06)} \\ \textcolor{red}{(failed)}}
  & \tabincell{c}{\textcolor{red}{25} \\ \textcolor{red}{(0.067)}} &\tabincell{c}{\textcolor{red}{4.03E+13} \\ \textcolor{red}{(8.85E+12)}  \\ \textcolor{red}{(failed)}}
  & \tabincell{c}{30+136 \\ (0.016)} &\tabincell{c}{3.80E-08 \\ (3.37E-08)} \\ \hline

\end{tabular}}
\end{table}

\vskip 2mm

\newpage

\begin{table}[!http]
  \newcommand{\tabincell}[2]{\begin{tabular}{@{}#1@{}}#2\end{tabular}}
  \scriptsize
  \centering
  \caption{Numerical results of Rcm, fmincon and SNOPT for large-scale problems with $n = 2000, m = 10$.}
  \label{TABCOMRSA3}
  \resizebox{\textwidth}{!}{
    \begin{tabular}{|c|c|c|c|c|c|c|c|c|c|}
  \hline
  \multirow{2}{*}{Problems }
  & \multicolumn{2}{c|}{Rcm}
  & \multicolumn{2}{c|}{fmincon}
  & \multicolumn{2}{c|}{SNOPT}  \\ \cline{2-7}
        & \tabincell{c}{$itc_{G}+itc_{R}$ \\(time)} & \tabincell{c}{KKT \\ $||Cons(x^{it})||_\infty$}
        & \tabincell{c}{steps \\(time)}   & \tabincell{c}{KKT \\ $||Cons(x^{it})||_\infty$}
        & \tabincell{c}{Major+ICTS \\(time)} & \tabincell{c}{KKT \\ $||Cons(x^{it})||_\infty$} \\ \hline

  \tabincell{c}{66. Trid Function \cite{GOT2015} \\ (n = 2000, m = 10)}
  & \tabincell{c}{1+12 \\ (1.10)} &\tabincell{c}{8.23e-06 \\ (3.55e-15)}
  & \tabincell{c}{13 \\ (7.72)} &\tabincell{c}{4.88e-06 \\ (1.30e-11)}
  & \tabincell{c}{4+1995 \\ (3.48)} &\tabincell{c}{7.20e-07 \\ (0)}  \\ \hline

  \tabincell{c}{67. Grewank Function \cite{GOT2015} \\ (n = 2000, m = 10)}
  & \tabincell{c}{6+10 \\ (3.59)} &\tabincell{c}{7.02e-06 \\ (4.98e-15)}
  & \tabincell{c}{\textcolor{red}{99} \\ \textcolor{red}{(69.52)}} &\tabincell{c}{\textcolor{red}{69.56} \\ \textcolor{red}{(1.12e-03)} \\ \textcolor{red}{(failed)}}
  & \tabincell{c}{12+2013 \\ (18.03)} &\tabincell{c}{2.000e-07 \\ (9.96e-12)}\\ \hline

  \tabincell{c}{68. Dixon Price Function \cite{GOT2015} \\ (n = 2000, m = 10)}
  & \tabincell{c}{5+9 \\ (1.05)} &\tabincell{c}{5.01e-06  \\ (7.72e-12)}
  & \tabincell{c}{25 \\ (19.37)} &\tabincell{c}{2.45e-06 \\ (3.17e-08)}
  & \tabincell{c}{31+2030 \\ (44.95)} &\tabincell{c}{9.20e-08 \\ (2.38e-10)} \\ \hline

  \tabincell{c}{69. Rosenbrock Function \cite{GOT2015} \\ (n = 2000, m = 10)}
  & \tabincell{c}{5+14 \\ (1.50)} &\tabincell{c}{5.21e-06  \\ (1.10e-08)}
  & \tabincell{c}{15 \\ (7.73)} &\tabincell{c}{1.44e-06 \\ (2.96e-06)}
  & \tabincell{c}{401+5952 \\ (549.72)} &\tabincell{c}{1.00e-06 \\ (5.97e-11)}  \\ \hline

  \tabincell{c}{70. Trigonometric Function  \cite{GOT2015} \\ (n = 2000, m = 10)}
  & \tabincell{c}{1+12 \\ (3.68)} &\tabincell{c}{8.58e-06  \\ (1.43e-12)}
  & \tabincell{c}{53 \\ (45.72)} &\tabincell{c}{5.98e-05 \\ (1.32e-10)}
  & \tabincell{c}{\textcolor{red}{97+4126} \\ \textcolor{red}{(425.56)}} &\tabincell{c}{\textcolor{red}{2.00e-06} \\ \textcolor{red}{(2.00e+03)}  \\ \textcolor{red}{(failed)}} \\ \hline

  \tabincell{c}{71. Singular Broyden Function \cite{GOT2015} \\ (n = 2000, m = 10)}
  & \tabincell{c}{10+34 \\ (3.72)} &\tabincell{c}{7.22e-06  \\ (8.45e-07)}
  & \tabincell{c}{32 \\ (18.91)} &\tabincell{c}{4.11e-06 \\ (4.91e-11)}
  & \tabincell{c}{\textcolor{red}{211+2320} \\ \textcolor{red}{(321.55)}} &\tabincell{c}{\textcolor{red}{1.500e-06} \\ \textcolor{red}{(2.11e+00)}  \\ \textcolor{red}{(failed)}} \\ \hline

  \tabincell{c}{72. Extended Powell Singular Function \cite{GOT2015} \\ (n = 2000, m = 10)}
  & \tabincell{c}{12+35 \\ (9.30)} &\tabincell{c}{9.28e-06  \\ (4.44e-16)}
  & \tabincell{c}{27 \\ (14.53)} &\tabincell{c}{7.14e-06 \\ (1.43e-13)}
  & \tabincell{c}{33+2025 \\ (41.11)} &\tabincell{c}{1.30e-06 \\ (6.11e-08)} \\ \hline

  \tabincell{c}{73. Tridiagonal System Function \cite{GOT2015} \\ (n = 2000, m = 10)}
  & \tabincell{c}{1+28 \\ (2.25)} &\tabincell{c}{5.99e-06  \\ (8.86e-15)}
  & \tabincell{c}{64 \\ (37.47)} &\tabincell{c}{8.08e-06 \\ (1.94e-10)}
  & \tabincell{c}{5+1996 \\ (8.89)} &\tabincell{c}{4.60e-07 \\ (0)} \\ \hline

  \tabincell{c}{74. Discrete Boundary-Value Function \cite{GOT2015} \\ (n = 2000, m = 10)}
  & \tabincell{c}{1+7 \\ (5.10)} &\tabincell{c}{4.05e-06  \\ (4.24e-22)}
  & \tabincell{c}{3 \\ (4.54)} &\tabincell{c}{4.01e-06 \\ (5.59e-14)}
  & \tabincell{c}{6+1997 \\ (10.80)} &\tabincell{c}{1.40e-07 \\ (1.88e-09)} \\ \hline

  \tabincell{c}{75. Broyden Tridiagonal Function \cite{GOT2015} \\ (n = 2000, m = 10)}
  & \tabincell{c}{2+9 \\ (0.92)} &\tabincell{c}{8.38e-07  \\ (2.15e-10)}
  & \tabincell{c}{13 \\ (7.38)} &\tabincell{c}{2.46e-06 \\ (1.48e-08)}
  & \tabincell{c}{\textcolor{red}{68+2072} \\ \textcolor{red}{(93.66)}} &\tabincell{c}{\textcolor{red}{1.40e-06} \\ \textcolor{red}{(1.35e+00)}  \\ \textcolor{red}{(failed)}}  \\ \hline

  \tabincell{c}{76. Extended Wood Function \cite{GOT2015} \\ (n = 2000, m = 10)}
  & \tabincell{c}{11+14 \\ (1.58)} &\tabincell{c}{1.20e-06  \\ (1.08e-06)}
  & \tabincell{c}{58 \\ (38.54)} &\tabincell{c}{6.27e-08 \\ (3.48e-09)}
  & \tabincell{c}{23+2013 \\ (35.87)} &\tabincell{c}{1.10e-07 \\ (9.54e-11)} \\ \hline

  \tabincell{c}{77. Extended Cliff Function \cite{GOT2015} \\ (n = 2000, m = 10)}
  & \tabincell{c}{6+7 \\ (1.11)} &\tabincell{c}{1.41e-06  \\ (6.66e-16)}
  & \tabincell{c}{9 \\ (4.79)} &\tabincell{c}{8.18e-06 \\ (1.35e-14)}
  & \tabincell{c}{8+2003 \\ (11.91)} &\tabincell{c}{2.00e-08 \\ (8.26e-06)} \\ \hline

  \tabincell{c}{78. Extended Hiebert Function \cite{GOT2015} \\ (n = 2000, m = 10)}
  & \tabincell{c}{1+56 \\ (4.95)} &\tabincell{c}{9.91e-07  \\ (2.71e-25)}
  & \tabincell{c}{\textcolor{red}{99} \\ \textcolor{red}{(49.24)}} &\tabincell{c}{\textcolor{red}{1.10e+03} \\ \textcolor{red}{(8.42e+04)}  \\ \textcolor{red}{(failed)}}
  & \tabincell{c}{7+1997 \\ (12.02)} &\tabincell{c}{3.30e-07 \\ (5.02e-09)}  \\ \hline

  \tabincell{c}{79. Extended Maratos Function \cite{GOT2015} \\ (n = 2000, m = 10)}
  & \tabincell{c}{10+6 \\ (1.06)} &\tabincell{c}{1.76e-06  \\ (1.42e-14)}
  & \tabincell{c}{\textcolor{red}{99} \\ \textcolor{red}{(63.59)}} &\tabincell{c}{\textcolor{red}{3.31e+03} \\ \textcolor{red}{(3.40e+02)}  \\ \textcolor{red}{(failed)}}
  & \tabincell{c}{19+2014 \\ (27.06)} &\tabincell{c}{5.20e-08 \\ (2.77e-05)}  \\ \hline

  \tabincell{c}{80. Extended Psc1 Function \cite{GOT2015} \\ (n = 2000, m = 10)}
  & \tabincell{c}{7+8 \\ (2.95)} &\tabincell{c}{5.93e-06  \\ (2.22e-16)}
  & \tabincell{c}{19 \\ (13.39)} &\tabincell{c}{3.77e-05 \\ (7.77e-16)}
  & \tabincell{c}{{11+2006} \\ {(18.78)}} &\tabincell{c}{{8.20e-09} \\ {(5.58e-08)}} \\ \hline

  \tabincell{c}{81. Extended Quadratic Penalty QP 1 \\ Function \cite{GOT2015}  (n = 2000, m = 10)}
  & \tabincell{c}{1+9 \\ (1.17)} &\tabincell{c}{2.61e-06  \\ (4.18e-41)}
  & \tabincell{c}{6 \\ (3.95)} &\tabincell{c}{3.17e-07 \\ (3.52e-10)}
  & \tabincell{c}{8+1998 \\ (13.36)} &\tabincell{c}{2.10e-07 \\ (0)} \\ \hline

  \tabincell{c}{82. Extended Quadratic Penalty QP 2  \\ Function \cite{GOT2015} (n = 2000, m = 10)}
  & \tabincell{c}{1+6 \\ (2.58)} &\tabincell{c}{2.32e-06  \\ (6.07e-39)}
  & \tabincell{c}{{5} \\ {(5.04)}} &\tabincell{c}{{1.73e-04} \\ {(1.60e+10)}}
  & \tabincell{c}{8+1998 \\ (17.22)} &\tabincell{c}{2.80e-08 \\ (0)}  \\ \hline

  \tabincell{c}{83. Extended TET Function \cite{GOT2015} \\ (n = 2000, m = 10)}
  & \tabincell{c}{7+6 \\ (3.34)} &\tabincell{c}{1.77e-06  \\ (0)}
  & \tabincell{c}{9 \\ (6.67)} &\tabincell{c}{5.15e-06 \\ (2.22e-16)}
  & \tabincell{c}{16+2011 \\ (26.36)} &\tabincell{c}{2.30e-09 \\ (1.27e-08)} \\ \hline

  \tabincell{c}{84. EG2 Function \cite{GOT2015} \\ (n = 2000, m = 10)}
  & \tabincell{c}{1+6 \\ (2.72)} &\tabincell{c}{5.10e-06  \\ (2.95e-07)}
  & \tabincell{c}{6 \\ (5.93)} &\tabincell{c}{7.18e-06 \\ (1.89e-07)}
  & \tabincell{c}{{6+1996} \\ {(12.55)}} &\tabincell{c}{{1.70e-07} \\ {(1.14e-05)}} \\ \hline

  \tabincell{c}{85. Extended BD1 Function \cite{GOT2015} \\ (n = 2000, m = 10)}
  & \tabincell{c}{5+6 \\ (1.46)} &\tabincell{c}{1.35e-06  \\ (0)}
  & \tabincell{c}{7 \\ (4.29)} &\tabincell{c}{8.82e-06 \\ (2.22e-15)}
  & \tabincell{c}{8+2003 \\ (12.09)} &\tabincell{c}{1.70e-06 \\ (9.65e-08)} \\ \hline

\end{tabular}}
\end{table}

\vskip 2mm

\begin{table}[!http]
  \newcommand{\tabincell}[2]{\begin{tabular}{@{}#1@{}}#2\end{tabular}}
  \scriptsize
  \centering
  \caption{Numerical results of Rcm, fmincon and SNOPT for large-scale problems with $n = 2000, m = 1000$.}
  \label{TABCOMRSA4}
  \resizebox{\textwidth}{!}{
    \begin{tabular}{|c|c|c|c|c|c|c|c|c|c|}
  \hline
  \multirow{2}{*}{Problems }
  & \multicolumn{2}{c|}{Rcm}
  & \multicolumn{2}{c|}{fmincon}
  & \multicolumn{2}{c|}{SNOPT}  \\ \cline{2-7}
        & \tabincell{c}{$itc_{G}+itc_{R}$ \\(time)} & \tabincell{c}{KKT \\ $||Cons(x^{it})||_\infty$}
        & \tabincell{c}{steps \\(time)}   & \tabincell{c}{KKT \\ $||Cons(x^{it})||_\infty$}
        & \tabincell{c}{Major+ICTS \\(time)} & \tabincell{c}{KKT \\ $||Cons(x^{it})||_\infty$} \\ \hline

  \tabincell{c}{86. Trid Function \cite{SB2013} \\ (n = 2000, m = 1000)}
  & \tabincell{c}{1+1 \\ (0.90)} &\tabincell{c}{2.65e-06  \\ (1.80e-06)}
  & \tabincell{c}{7 \\ (17.97)} &\tabincell{c}{1.98e-05 \\ (5.82e-11)}
  & \tabincell{c}{0+1 \\ (0.02)} &\tabincell{c}{1.40e-13 \\ (3.64e-12)} \\ \hline

  \tabincell{c}{87. Grewank Function \cite{SB2013} \\ (n = 2000, m = 1000)}
  & \tabincell{c}{14+6 \\ (8.66)} &\tabincell{c}{4.07e-06  \\ (3.30e-12)}
  & \tabincell{c}{59 \\ (1.07e+02)} &\tabincell{c}{5.21e-06 \\ (7.70e-164)}
  & \tabincell{c}{{32+2035} \\ {(84.20)}} &\tabincell{c}{{1.40e-06} \\ {(4.64e-13)}} \\ \hline

  \tabincell{c}{88. Dixon Price Function \cite{SB2013} \\ (n = 2000, m = 1000)}
  & \tabincell{c}{8+11 \\ (7.28)} &\tabincell{c}{2.34e-06  \\ (1.01e-07)}
  & \tabincell{c}{13 \\ (17.11)} &\tabincell{c}{7.09e-07 \\ (7.41e-07)}
  & \tabincell{c}{219+1254 \\ (41.47)} &\tabincell{c}{8.20e-07 \\ (1.82e-09)}  \\ \hline

  \tabincell{c}{89. Rosenbrock Function \cite{MGH1981} \\ (n = 2000, m = 1000)}
  & \tabincell{c}{6+14 \\ (7.60)} &\tabincell{c}{7.02e-06  \\ (1.70e-08)}
  & \tabincell{c}{16 \\ (20.98)} &\tabincell{c}{4.59e-06 \\ (2.11e-05)}
  & \tabincell{c}{276+2727 \\ (55.38)} &\tabincell{c}{7.00e-08 \\ (2.19e-08)}  \\ \hline

  \tabincell{c}{90. Trigonometric Function  \cite{MGH1981} \\ (n = 2000, m = 1000)}
  & \tabincell{c}{8+7 \\ (7.67)} &\tabincell{c}{9.74e-06  \\ (5.58e-11)}
  & \tabincell{c}{\textcolor{red}{94} \\ \textcolor{red}{(1.23e+02)}} &\tabincell{c}{\textcolor{red}{2.45e-02} \\ \textcolor{red}{(4.80e-11)}  \\ \textcolor{red}{(failed)}}
  & \tabincell{c}{\textcolor{red}{245+7140} \\ \textcolor{red}{(1003.45)}} &\tabincell{c}{\textcolor{red}{5.90e+02} \\ \textcolor{red}{(1.28e+03)}  \\ \textcolor{red}{(failed)}} \\ \hline

  \tabincell{c}{91. Singular Broyden Function \cite{Luksan1994} \\ (n = 2000, m = 1000)}
  & \tabincell{c}{10+15 \\ (9.97)} &\tabincell{c}{8.38e-06  \\ (1.92e-06)}
  & \tabincell{c}{29 \\ (39.45)} &\tabincell{c}{3.25e-05 \\ (3.24e-12)}
  & \tabincell{c}{\textcolor{red}{5413+84179} \\ \textcolor{red}{(91.94)}} &\tabincell{c}{\textcolor{red}{1.20e-03} \\ \textcolor{red}{(2.16e+00)}  \\ \textcolor{red}{(failed)}}  \\ \hline

  \tabincell{c}{92. Extended Powell Singular Function \cite{MGH1981} \\ (n = 2000, m = 1000)}
  & \tabincell{c}{12+10 \\ (7.89)} &\tabincell{c}{8.45e-06  \\ (2.05e-14)}
  & \tabincell{c}{16 \\ (21.15)} &\tabincell{c}{4.34e-06 \\ (2.95e-09)}
  & \tabincell{c}{18+1270 \\ (6.06)} &\tabincell{c}{1.80e-06 \\ (1.62e-05)} \\ \hline

  \tabincell{c}{93. Tridiagonal System Function \cite{Luksan1994} \\ (n = 2000, m = 1000)}
  & \tabincell{c}{1+22 \\ (8.64)} &\tabincell{c}{9.80e-06  \\ (1.55e-15)}
  & \tabincell{c}{27 \\ (36.08)} &\tabincell{c}{9.80e-06 \\ (1.34e-08)}
  & \tabincell{c}{{4+1006} \\ {(1.89)}} &\tabincell{c}{{1.20e-07} \\ {(0)}} \\ \hline

  \tabincell{c}{94. Discrete Boundary-Value Function \cite{Luksan1994} \\ (n = 2000, m = 1000)}
  & \tabincell{c}{1+7 \\ (8.78)} &\tabincell{c}{2.51e-06  \\ (6.94e-18)}
  & \tabincell{c}{6 \\ (12.44)} &\tabincell{c}{7.90e-07 \\ (2.76e-14)}
  & \tabincell{c}{{5+1007} \\ {(2.28)}} &\tabincell{c}{{5.20e-07} \\ {(2.96e-07)}} \\ \hline

  \tabincell{c}{95. Broyden Tridiagonal Function \cite{Luksan1994} \\ (n = 2000, m = 1000)}
  & \tabincell{c}{2+10 \\ (5.04)} &\tabincell{c}{4.43e-06  \\ (5.46e-09)}
  & \tabincell{c}{14 \\ (18.19)} &\tabincell{c}{4.43e-06 \\ (6.80e-09)}
  & \tabincell{c}{\textcolor{red}{11+15} \\ \textcolor{red}{(0.03)}} &\tabincell{c}{\textcolor{red}{9.40e-14} \\ \textcolor{red}{(1.41)}  \\ \textcolor{red}{(failed)}} \\ \hline

  \tabincell{c}{96. Extended Wood Function \cite{Andrei2008} \\ (n = 2000, m = 1000)}
  & \tabincell{c}{11+6 \\ (6.58)} &\tabincell{c}{1.29e-06  \\ (0)}
  & \tabincell{c}{58 \\ (38.54)} &\tabincell{c}{6.27e-08 \\ (3.48e-09)}
  & \tabincell{c}{13+1013 \\ (4.77)} &\tabincell{c}{1.40e-06 \\ (8.53e-05)} \\ \hline

  \tabincell{c}{97. Extended Cliff Function \cite{Andrei2008} \\ (n = 2000, m = 1000)}
  & \tabincell{c}{6+7 \\ (5.61)} &\tabincell{c}{1.41e-06  \\ (6.66e-16)}
  & \tabincell{c}{13 \\ (16.95)} &\tabincell{c}{4.01e-09 \\ (6.66e-16)}
  & \tabincell{c}{12+1512 \\ (3.34)} &\tabincell{c}{3.30e-09 \\ (5.82e-05)} \\ \hline

  \tabincell{c}{98. Extended Hiebert Function \cite{Andrei2008} \\ (n = 2000, m = 1000)}
  & \tabincell{c}{1+43 \\ (32.32)} &\tabincell{c}{1.55e-06  \\ (2.96e-23)}
  & \tabincell{c}{\textcolor{red}{99} \\ \textcolor{red}{(1.24e+02)}} &\tabincell{c}{\textcolor{red}{2.03e+04} \\ \textcolor{red}{(3.09e+04)}  \\ \textcolor{red}{(failed)}}
  & \tabincell{c}{6+1006 \\ (2.09)} &\tabincell{c}{1.60e-07 \\ (5.02e-09)} \\ \hline

  \tabincell{c}{99. Extended Maratos Function \cite{Andrei2008} \\ (n = 2000, m = 1000)}
  & \tabincell{c}{10+6 \\ (6.31)} &\tabincell{c}{2.00e-06  \\ (1.42e-14)}
  & \tabincell{c}{90 \\ (2.01e+02)} &\tabincell{c}{2.52e-07 \\ (1.42e-14)}
  & \tabincell{c}{37+1537 \\ (8.81)} &\tabincell{c}{2.80e-07 \\ (5.34e-05)} \\ \hline

  \tabincell{c}{100. Extended Psc1 Function \cite{Andrei2008} \\ (n = 2000, m = 1000)}
  & \tabincell{c}{7+8 \\ (7.49)} &\tabincell{c}{5.94e-06  \\ (5.55e-16)}
  & \tabincell{c}{48 \\ (68.38)} &\tabincell{c}{3.95e-05 \\ (8.21e-12)}
  & \tabincell{c}{12+1512 \\ (2.34)} &\tabincell{c}{6.40e-07 \\ (2.43e-13)} \\ \hline

  \tabincell{c}{101. Extended Quadratic Penalty QP 1 \\ Function \cite{Andrei2008}  (n = 2000, m = 1000)}
  & \tabincell{c}{6+5 \\ (4.72)} &\tabincell{c}{6.86e-06  \\ (8.31e-34)}
  & \tabincell{c}{10 \\ (14.46)} &\tabincell{c}{1.85e-06 \\ (8.40e-10)}
  & \tabincell{c}{12+1012 \\ (98.39)} &\tabincell{c}{1.20e-08 \\ (7.88e-06)} \\ \hline

  \tabincell{c}{102. Extended Quadratic Penalty QP 2  \\ Function \cite{Andrei2008} (n = 2000, m = 1000)}
  & \tabincell{c}{6+6 \\ (6.19)} &\tabincell{c}{2.20e-06  \\ (6.62e-35)}
  & \tabincell{c}{11 \\ (36.49)} &\tabincell{c}{1.22e-06 \\ (1.99e-12)}
  & \tabincell{c}{\textcolor{red}{89+1094} \\ \textcolor{red}{(261.47)}} &\tabincell{c}{\textcolor{red}{1.60e+01} \\ \textcolor{red}{(0.32)}  \\ \textcolor{red}{(failed)}} \\ \hline

  \tabincell{c}{103. Extended TET Function \cite{Andrei2008} \\ (n = 2000, m = 1000)}
  & \tabincell{c}{7+6 \\ (7.19)} &\tabincell{c}{4.67e-06  \\ (0)}
  & \tabincell{c}{11 \\ (16.70)} &\tabincell{c}{6.61e-07 \\ (0)}
  & \tabincell{c}{18+1518 \\ (4.88)} &\tabincell{c}{5.90e-07 \\ (7.02e-14)} \\ \hline

  \tabincell{c}{104. EG2 Function \cite{Andrei2008} \\ (n = 2000, m = 1000)}
  & \tabincell{c}{1+6 \\ (5.06)} &\tabincell{c}{5.11e-06  \\ (2.95e-07)}
  & \tabincell{c}{6 \\ (12.34)} &\tabincell{c}{7.72e-06 \\ (1.59e-07)}
  & \tabincell{c}{{6+1006} \\ {(2.81)}} &\tabincell{c}{{3.60e-08} \\ {(1.89e-08)}} \\ \hline

  \tabincell{c}{105. Extended BD1 Function \cite{Andrei2008} \\ (n = 2000, m = 1000)}
  & \tabincell{c}{5+6 \\ (5.08)} &\tabincell{c}{1.89e-06  \\ (0)}
  & \tabincell{c}{9 \\ (13.34)} &\tabincell{c}{2.50e-07 \\ (0)}
  & \tabincell{c}{9+1509 \\ (2.58)} &\tabincell{c}{2.30e-07 \\ (6.53e-05)} \\ \hline

\end{tabular}}
\end{table}

\vskip 2mm

\begin{table}[!http]
  \newcommand{\tabincell}[2]{\begin{tabular}{@{}#1@{}}#2\end{tabular}}
  \scriptsize
  \centering
  \caption{Numerical results of Rcm, fmincon and SNOPT for large-scale problems with $n = 2000, \; m = 1999$.}
  \label{TABCOMRSA5}
  \resizebox{\textwidth}{!}{
    \begin{tabular}{|c|c|c|c|c|c|c|c|c|c|}
  \hline
  \multirow{2}{*}{Problems }
  & \multicolumn{2}{c|}{Rcm}
  & \multicolumn{2}{c|}{fmincon}
  & \multicolumn{2}{c|}{SNOPT}  \\ \cline{2-7}
        & \tabincell{c}{$itc_{G}+itc_{R}$ \\(time)} & \tabincell{c}{KKT \\ $||Cons(x^{it})||_\infty$}
        & \tabincell{c}{steps \\(time)}   & \tabincell{c}{KKT \\ $||Cons(x^{it})||_\infty$}
        & \tabincell{c}{Major+ICTS \\(time)} & \tabincell{c}{KKT \\ $||Cons(x^{it})||_\infty$} \\ \hline

  \tabincell{c}{106. Trid Function \cite{SB2013} \\ (n = 2000, m = 1999)}
  & \tabincell{c}{1+1 \\ (1.95)} &\tabincell{c}{9.21e-06  \\ (1.67e-08)}
  & \tabincell{c}{2 \\ (12.19)} &\tabincell{c}{2.54e-06 \\ (2.33e-10)}
  & \tabincell{c}{0+1 \\ (0.02)} &\tabincell{c}{3.70e-12 \\ (7.28e-12)} \\ \hline

  \tabincell{c}{107. Grewank Function \cite{SB2013} \\ (n = 2000, m = 1999)}
  & \tabincell{c}{21+8 \\ (12.67)} &\tabincell{c}{9.48e-06  \\ (7.76e-17)}
  & \tabincell{c}{{32} \\ {(1.86e+02)}} &\tabincell{c}{{5.47e-04} \\ {(2.53e-15)}}
  & \tabincell{c}{30+2037 \\ (121.20)} &\tabincell{c}{1.70e-07 \\ (4.11e-10)} \\ \hline

  \tabincell{c}{108. Dixon Price Function \cite{SB2013} \\ (n = 2000, m = 1999)}
  & \tabincell{c}{8+6 \\ (5.32)} &\tabincell{c}{2.20e-06  \\ (1.42e-09)}
  & \tabincell{c}{\textcolor{red}{97} \\ \textcolor{red}{(5.44e+02)}} &\tabincell{c}{\textcolor{red}{1.25e+03} \\ \textcolor{red}{(1.95e+02)}  \\ \textcolor{red}{(failed)}}
  & \tabincell{c}{11+16 \\ (0.03)} &\tabincell{c}{1.10e-08 \\ (7.46e-07)} \\ \hline

  \tabincell{c}{109. Rosenbrock Function \cite{MGH1981} \\ (n = 2000, m = 1999)}
  & \tabincell{c}{6+27 \\ (21.10)} &\tabincell{c}{2.25e-06 \\ (6.20e-09)}
  & \tabincell{c}{20 \\ (64.61)} &\tabincell{c}{1.39e-06 \\ (1.07e-06)}
  & \tabincell{c}{253+223 \\ (3.19)} &\tabincell{c}{6.30e-05 \\ (7.58e-07)} \\ \hline

  \tabincell{c}{110. Trigonometric Function  \cite{MGH1981} \\ (n = 2000, m = 1999)}
  & \tabincell{c}{1+13 \\ (15.29)} &\tabincell{c}{1.24e-06  \\ (1.13e-15)}
  & \tabincell{c}{\textcolor{red}{15} \\ \textcolor{red}{(68.47)}} &\tabincell{c}{\textcolor{red}{1.41e-03} \\ \textcolor{red}{(2.89e-13)}  \\ \textcolor{red}{(failed)}}
  & \tabincell{c}{\textcolor{red}{96+10417} \\ \textcolor{red}{(1003.13)}} &\tabincell{c}{\textcolor{red}{1.60e+03} \\ \textcolor{red}{(9.71e+02)}  \\ \textcolor{red}{(failed)}} \\ \hline

  \tabincell{c}{111. Singular Broyden Function \cite{Luksan1994} \\ (n = 2000, m = 1999)}
  & \tabincell{c}{10+6 \\ (19.26)} &\tabincell{c}{1.84e-06  \\ (6.40e-07)}
  & \tabincell{c}{\textcolor{red}{74} \\ \textcolor{red}{(2.87e+02)}} &\tabincell{c}{\textcolor{red}{28.52} \\ \textcolor{red}{(9.41e-13)}  \\ \textcolor{red}{(failed)}}
  & \tabincell{c}{\textcolor{red}{49586+606083} \\ \textcolor{red}{(480.78)}} &\tabincell{c}{\textcolor{red}{3.70e-03} \\ \textcolor{red}{(2.16e+00)}  \\ \textcolor{red}{(failed)}} \\ \hline

  \tabincell{c}{112. Extended Powell Singular Function \cite{MGH1981} \\ (n = 2000, m = 1999)}
  & \tabincell{c}{12+26 \\ (21.14)} &\tabincell{c}{5.45e-06  \\ (6.38e-11)}
  & \tabincell{c}{27 \\ (1.04e+02)} &\tabincell{c}{6.46e-06 \\ (7.70e-12)}
  & \tabincell{c}{{23+525} \\{(0.22)}} &\tabincell{c}{{4.00e-07} \\ {(5.30e-08)}} \\ \hline

  \tabincell{c}{113. Tridiagonal System Function \cite{Luksan1994} \\ (n = 2000, m = 1999)}
  & \tabincell{c}{1+15 \\ (9.00)} &\tabincell{c}{9.28e-06  \\ (1.55e-15)}
  & \tabincell{c}{35 \\ (1.08e+02)} &\tabincell{c}{1.28e-02 \\ (4.88e-15)}
  & \tabincell{c}{1+3 \\ (0.02)} &\tabincell{c}{6.10e-18 \\ (0)} \\ \hline

  \tabincell{c}{114. Discrete Boundary-Value Function \cite{Luksan1994} \\ (n = 2000, m = 1999)}
  & \tabincell{c}{1+7 \\ (10.60)} &\tabincell{c}{7.57e-06  \\ (8.32e-17)}
  & \tabincell{c}{17 \\ (70.40)} &\tabincell{c}{1.16e-06 \\ (2.49e-16)}
  & \tabincell{c}{7+10 \\ (0.03)} &\tabincell{c}{1.60e-07 \\ (8.72e-07)} \\ \hline

  \tabincell{c}{115. Broyden Tridiagonal Function \cite{Luksan1994} \\ (n = 2000, m = 1999)}
  & \tabincell{c}{2+7 \\ (6.32)} &\tabincell{c}{2.47e-06  \\ (7.15e-10)}
  & \tabincell{c}{11 \\ (33.86)} &\tabincell{c}{2.77e-06 \\ (2.46e-06)}
  & \tabincell{c}{\textcolor{red}{11+16} \\ \textcolor{red}{(0.11)}} &\tabincell{c}{\textcolor{red}{1.20e-08} \\ \textcolor{red}{(1.41e+00)}  \\ \textcolor{red}{(failed)}} \\ \hline

  \tabincell{c}{116. Extended Wood Function \cite{Andrei2008} \\ (n = 2000, m = 1999)}
  & \tabincell{c}{11+5 \\ (9.35)} &\tabincell{c}{3.69e-06  \\ (3.79e-09)}
  & \tabincell{c}{37 \\ (1.46e+02)} &\tabincell{c}{2.56e-08 \\ (4.14e-11)}
  & \tabincell{c}{{54+55} \\ {(0.23)}} &\tabincell{c}{{7.00e-09} \\ {(3.30e-08)}} \\ \hline

  \tabincell{c}{117. Extended Cliff Function \cite{Andrei2008} \\ (n = 2000, m = 1999)}
  & \tabincell{c}{6+7 \\ (8.60)} &\tabincell{c}{1.41e-06  \\ (6.66e-16)}
  & \tabincell{c}{8 \\ (45.13)} &\tabincell{c}{1.29e-08 \\ (1.13e-13)}
  & \tabincell{c}{{9+1009} \\ {(0.09)}} &\tabincell{c}{{1.30e-08} \\ {(7.52e-07)}} \\ \hline

  \tabincell{c}{118. Extended Hiebert Function \cite{Andrei2008} \\ (n = 2000, m = 1999)}
  & \tabincell{c}{1+10 \\ (7.61)} &\tabincell{c}{8.09e-06  \\ (3.28e-42)}
  & \tabincell{c}{5 \\ (16.01)} &\tabincell{c}{1.90e-06 \\ (2.02e-13)}
  & \tabincell{c}{6+12 \\ (0.03)} &\tabincell{c}{1.40e-06 \\ (5.02e-09)} \\ \hline

  \tabincell{c}{119. Extended Maratos Function \cite{Andrei2008} \\ (n = 2000, m = 1999)}
  & \tabincell{c}{10+14 \\ (12.98)} &\tabincell{c}{7.99e-06  \\ (2.20e-06)}
  & \tabincell{c}{\textcolor{red}{99} \\ \textcolor{red}{(3.08e+02)}} &\tabincell{c}{\textcolor{red}{73.00} \\ \textcolor{red}{(49.68)}  \\ \textcolor{red}{(failed)}}
  & \tabincell{c}{{27+1027} \\ {(0.23)}} &\tabincell{c}{{1.50e-06} \\ {(6.32e-10)}} \\ \hline

  \tabincell{c}{120. Extended Psc1 Function \cite{Andrei2008} \\ (n = 2000, m = 1999)}
  & \tabincell{c}{7+14 \\ (13.73)} &\tabincell{c}{2.20e-06  \\ (3.33e-16)}
  & \tabincell{c}{14 \\ (43.92)} &\tabincell{c}{9.75e-06 \\ (4.49e-11)}
  & \tabincell{c}{{12+1012} \\ {(0.13)}} &\tabincell{c}{{7.50e-07} \\ {(7.61e-07)}} \\ \hline

  \tabincell{c}{121. Extended Quadratic Penalty QP 1 \\ Function \cite{Andrei2008}  (n = 2000, m = 1999)}
  & \tabincell{c}{13+19 \\ (17.29)} &\tabincell{c}{9.63e-06  \\ (7.22e-07)}
  & \tabincell{c}{\textcolor{red}{53} \\ \textcolor{red}{(1.64e+02)}} &\tabincell{c}{\textcolor{red}{8.95e-02} \\ \textcolor{red}{(1.74e-61)}  \\ \textcolor{red}{(failed)}}
  & \tabincell{c}{16+17 \\ (121.86)} &\tabincell{c}{1.30e-09 \\ (6.16e-09)} \\ \hline

  \tabincell{c}{122. Extended Quadratic Penalty QP 2  \\ Function \cite{Andrei2008} (n = 2000, m = 1999)}
  & \tabincell{c}{14+16 \\ (11.52)} &\tabincell{c}{2.42e-06  \\ (1.27e-08)}
  & \tabincell{c}{12 \\ (44.11)} &\tabincell{c}{4.27e-06 \\ (8.08e-10)}
  & \tabincell{c}{\textcolor{red}{1+2} \\ \textcolor{red}{(605.25)}} &\tabincell{c}{\textcolor{red}{3.00e-04} \\ \textcolor{red}{(1.59e+10)}  \\ \textcolor{red}{(failed)}} \\ \hline

  \tabincell{c}{123. Extended TET Function \cite{Andrei2008} \\ (n = 2000, m = 1999)}
  & \tabincell{c}{7+7 \\ (10.34)} &\tabincell{c}{5.45e-06  \\ (5.50e-12)}
  & \tabincell{c}{12 \\ (38.43)} &\tabincell{c}{4.53e-07 \\ (1.26e-07)}
  & \tabincell{c}{21+1021 \\ (0.17)} &\tabincell{c}{1.30e-06 \\ (9.35e-08)}  \\ \hline

  \tabincell{c}{124. EG2 Function \cite{Andrei2008} \\ (n = 2000, m = 1999)}
  & \tabincell{c}{3+6 \\ (7.80)} &\tabincell{c}{1.30e-06  \\ (3.68e-13)}
  & \tabincell{c}{8 \\ (26.53)} &\tabincell{c}{1.80e-06 \\ (3.68e-13)}
  & \tabincell{c}{10+11 \\ (0.27)} &\tabincell{c}{2.50e-08 \\ (2.99e-10)} \\ \hline

  \tabincell{c}{125. Extended BD1 Function \cite{Andrei2008} \\ (n = 2000, m = 1999)}
  & \tabincell{c}{5+21 \\ (14.41)} &\tabincell{c}{3.88e-06  \\ (2.03e-09)}
  & \tabincell{c}{13 \\ (42.83)} &\tabincell{c}{8.68e-08 \\ (8.60e-10)}
  & \tabincell{c}{14+1014 \\ (0.19)} &\tabincell{c}{6.50e-07 \\ (6.60e-06)} \\ \hline

\end{tabular}}
\end{table}

\vskip 2mm

\begin{table}[htbp]
  \renewcommand{\arraystretch}{1.5}
  \newcommand{\tabincell}[2]{\begin{tabular}{@{}#1@{}}#2\end{tabular}}
  \centering
  \caption{The number of failed problems computed by Rcm, fmincon and SNOPT.}
  \label{TABPFRCM}
  \resizebox{\textwidth}{!}{
    \begin{tabular}{|c|c|c|c|c|c|c|}
    \hline

    & Rcm & fmincon & SNOPT \cr\hline

    The number of failed problems
    & 8 & 28  & {21}   \cr \hline

    {The probability of failure}
    & $\frac{8}{125} \; (6.4\%)$ & $\frac{28}{125} \; (22.4\%)$ & $\frac{21}{125} \; (16.8\%)$  \cr \hline \hline

    The number of fasted problems
    & 48 & 3 & 72  \cr \hline

    The probability of fasted problem
    & $\frac{48}{125} \; (38.4 \%)$ & $\frac{3}{125} \; (2.4 \%)$ & $\frac{72}{125} \; (57.6 \%)$ \cr \hline

    \end{tabular}}
\end{table}
\vskip 2mm

\clearpage

\section{Conclusions}

\vskip 2mm

In this paper, we give the regularization continuation method with the
trust-region updating strategy (Rcm) for nonlinear equality-constrained optimization
problems. Namely, we use the inverse of the regularization quasi-Newton matrix
as the pre-conditioner to improve its computational efficiency in the well-posed
phase, and we use the inverse of the regularization two-sided projection of the
Hessian matrix as the pre-conditioner to improve its robustness in the
ill-conditioned phase. Since Rcm only solves a linear system of equations at
every iteration and SQP needs to solve a quadratic programming subproblem at
every iteration, Rcm is faster than SQP. Numerical results also show that Rcm
is more robust and faster than SQP (the built-in subroutine fmincon.m of the
MATLAB2020a environment \cite{MATLAB,Schittkowski1986} and the subroutine
SNOPT \cite{GMS2005,GMS2006} executed in GAMS v28.2 (2019) environment
\cite{GAMS}). The computational time of Rcm is about one third of that fmincon.m
for the large-scale problem. Therefore, Rcm is an alternating solver for
equality-constrained optimization problems and the regularization continuation
method is worthy to be explored further for inequality-constrained optimization
problems and the orthogonal matrix constraint problems \cite{WY2013}.

\vskip 2mm

\section*{Acknowledgments} This work was supported in part by Grant 61876199
from National Natural Science Foundation of China, Grant YBWL2011085 from
Huawei Technologies Co., Ltd., and Grant YJCB2011003HI from the Innovation
Research Program of Huawei Technologies Co., Ltd.. The authors are grateful
to two anonymous referees for their comments and suggestions which greatly
improve the presentation of this paper.

\vskip 2mm

\noindent \textbf{Conflicts of interest / Competing interests:} Not applicable.

\vskip 2mm

\noindent \textbf{Availability of data and material (data transparency):} If it is requested, we will
provide the test data.

\vskip 2mm

\noindent \textbf{Code availability (software application or custom code):} If it is requested, we will
provide the code.

\end{document}